\documentclass{article}
\usepackage{graphicx} 
\usepackage{amsmath}
\usepackage{amsthm} 
\usepackage{amsfonts} 
\usepackage{tikz} 
\usepackage{hyperref} 
\usepackage{mathtools}
\usepackage{amssymb} 
\usepackage{bbold}
\usepackage{geometry}
\usetikzlibrary{decorations.pathmorphing} 
\usepackage{enumitem} 
\usepackage{pgfplots} 
\usepackage{caption}
\usepackage{subcaption}
\usepackage[dvipsnames]{xcolor}
\usepackage{float} 
\usepackage{bm} 

\newcommand{\perron}{\mathcal{L}}
\newcommand{\D}{\operatorname{D}}
\newcommand{\Time}{\mathcal{T}}
\newcommand{\V}[1]{\left\|#1\right\|_{V}}
\newcommand{\Va}[1]{\left\|#1\right\|_{V_\varepsilon}}
\newcommand{\BV}[1]{\left\|#1\right\|_{BV}}
\newcommand{\BVa}[1]{\left\|#1\right\|_{BV_\varepsilon}}

\newcommand{\Lnorm}[1]{\left\|#1\right\|_{1}}

\newcommand{\N}{\mathbb{N}}

\newcommand{\f}{f^*}

\newcommand{\T}{\hat{T}}
\newcommand{\R}{\mathbb{R}}

\newcommand{\B}{\mathcal{B}}

\newcommand{\X}{\bm{X}}
\newcommand{\Te}{\T_{\varepsilon}}

\newtheorem{prop}{Proposition}[section]
\newtheorem{definition}{Definition}
\newtheorem{theorem}{Theorem}
\newtheorem{remark}{Remark}
\newtheorem{lemma}{Lemma}
\newtheorem{question}{Question}
\newenvironment{customproof}[1]
    {\par\vspace{1em}\noindent\textbf{#1.}\;}
    {\hfill$\qed$\par\vspace{1em}}

\title{ACIM instability of piecewise expanding maps through the lens of metastability}

\author{
Péter Bálint\thanks{
Department of Stochastics, Institute of Mathematics,
Budapest University of Technology and Economics,
M{\H u}egyetem rkp. 3., H-1111, Budapest, Hungary; HUN-REN--BME Stochastics Research Group,
Budapest University of Technology and Economics,
M{\H u}egyetem rkp. 3., H-1111 Budapest, Hungary;  Email:
\textit{balint.peter@ttk.bme.hu}. Support of the NKFIH Research Fund, grants 142169 and 144059, is thankfully acknowledged.}
\qquad
Ábel Komálovics\thanks{
 Department of Stochastics, Institute of Mathematics,
Budapest University of Technology and Economics,
M{\H u}egyetem rkp. 3., H-1111, Budapest, Hungary; Email: \textit{abelkomalovics@edu.bme.hu}. Project no. 294 has been implemented with the support provided by the Ministry of Culture and Innovation of Hungary from the National Research, Development and Innovation Fund, financed under the DKÖP-26-1-BME-8 funding scheme.}
}

\date{September 10, 2026}

\begin{document}

\maketitle

\begin{abstract}
    Motivated by Keller's W-shaped maps and its variants, we introduce a general class of families of expanding maps such that the perturbed maps have a shrinking almost invariant neighborhood about a fixed point of the limit map. The unique absolutely continuous invariant measures (ACIM) of the perturbed maps can converge to limit measures of various types. A local quantity is identified which determines if the limit is absolutely continuous, singular, or a non-trivial convex combination of these. Furthermore, for the case of a nontrivial convex combination, we prove that the dynamics of the perturbed system, when viewed on an appropriate slow time scale, converges to a jump Markov process, a convergence that extends to the diffusion coefficients for observables of bounded variation. Compared to analogous results on expanding maps with metastable behavior, a new feature of our setting is the emergence of a localized state of the Markov process which corresponds to the shrinking almost invariant interval about the fixed point. Our approach provides a general framework which, in particular, accommodates, to the best of our knowledge, all previously studied families that limit to Keller's W-shaped map.
\end{abstract}

\section*{Introduction}

The aim of the present paper is to explore new connections between two important phenomena of chaotic dynamical systems: statistical stability (more precisely, lack thereof) and metastable dynamics. Our exposition restricts to piecewise expanding maps of the interval, a context in which both of the above phenomena have been extensively studied. It is beyond the scope of the present paper to give a systematic review of the rich literature on these topics, our aim instead is to place our results in the relevant context. In short, the aim of our work is to show that, for a substantial class of examples, statistical instability can be coherently described in terms of the multiple time scale dynamics arising in an appropriate metastable system.

For piecewise expanding maps of the interval, \emph{statistical stability} typically reduces to ACIM stability (where the acronym ACIM refers to absolutely continuous invariant measure). Generally speaking, when modeling real world situations we often assume that there is an underlying chaotic dynamics, possessing a natural invariant measure, with respect to which we want to calculate the possibilities of certain events. Nonetheless, the dynamics can only be approximated, given measurement errors, and it would be convenient if the probabilities of events measured with respect to the approximate dynamics were close to the original probabilities. This motivates the study of \emph{ACIM stability}, which can be roughly described as follows.
Let $T_\varepsilon$, $\varepsilon\geq 0$ be a family of piecewise expanding maps of the interval $I$, all having unique ACIMs $\mu_\varepsilon$. Is it true that $T_\varepsilon \to T_0$ implies $\mu_\varepsilon\to \mu_0$?

This is not the case for all families. A well known example that demonstrates ACIM instability dates back to \cite{Keller82}, where Keller constructed a family of piecewise linear expanding maps of the interval $[0,1]$,  with unique invariant densities, which nevertheless converged to a Dirac point mass as $\varepsilon \to 0$. The limiting map of $\varepsilon=0$ is often called  Keller's W map, given the shape of its graph. The tip of the $W$ at $1/2$ corresponds to a critical fixed point (in other words, a turning point) on which the Dirac point mass is supported. Keller conjectured that such instability can only occur if there are shrinking neighborhoods of a periodic critical point of the limit map which are left invariant by the perturbed maps. During the early 2010's several papers \cite{singular limits, Ulams method, harmonic averages, various instabilities, spectrum instability, family} disproved this hypothesis, by constructing still piecewise linear variants of the original family from \cite{Keller82}, which all converged to Keller's W-map, yet exhibited somewhat different behavior.
Specifically, \cite{family} showed that with a simple, non-Markovian (yet piecewise linear) perturbation of this map, invariant densities can be obtained which converge to a nontrivial convex combination of the ACIM of the W-map and the Dirac mass on the fixed point at $1/2$. The perturbations did not necessarily leave any neighborhood of $1/2$ invariant. This family was reconsidered in \cite{Ulams method}, where the authors showed that although the densities do not converge in the above family, Ulam's method converges. \cite{spectrum instability} examined the instability of the spectrum of the Perron-Frobenius operator associated to another family of W-shaped maps. These maps were no longer symmetric, yet, could still produce arbitrary nontrivial convex coefficients in the limit measure. \cite{singular limits} showed that by introducing more parameters, Markovian families can be obtained for which the coefficients of the convex combination in the limit measure can be again tuned arbitrarily. This includes the trivial cases, when the limit measure is either the ACIM or the Dirac mass. \cite{various instabilities} introduced a family for which the type of the limiting measure was connected to the slopes of the maps neighboring the fixed point. If the harmonic averages of the two slopes were below $2$, the limit was the Dirac mass. If the harmonic average was exactly $2$, it could be tuned to produce any nontrivial convex coefficients. If the harmonic average was above $2$, the limit was the ACIM.


Although all these papers introduced novelty and brought forward the study of the problem of ACIM instability, they all used methods mostly suited for piecewise linear maps, and calculations that fit their parametrization of their specific family. This is a natural approach, yet it limits to a specific (in particular, piecewise linear) class of examples which can be studied with these ad hoc methods. In the present paper we show that these phenomena can be studied in a more systematic way, which not only allows for the description of a far larger (nonlinear) class of examples, but, we believe, also connects the description of ACIM instability to its true dynamical origins. As already observed in the above papers, in all of the studied families the fixed turning point of the limit map has shrinking neighborhoods that are \emph{almost invariant} under the perturbed maps. This results in the phenomenon where every perturbation has two partitioning sets that are almost invariant under the dynamics, one of which converges to a point. Furthermore, there is a \emph{flow of points} between the two sets, hence a positive mass may accumulate both in the shrinking set and its complement.

This phenomenon is reminiscent of the behavior of \textit{metastable systems}, which, in the specific setting of piecewise expanding interval maps, can be described as follows.
Let $T_0$ be a piecewise expanding $C^2$ map of the interval $I=[0,1]$, and let us assume $I$ can be partitioned into invariant nontrivial subintervals $I_j$, such that each $I_j$ supports a unique ACIM $\mu_j$ and $T_0|_{I_j}$ is mixing for each $\mu_j$. Take $C^2$-small perturbations  $T_\varepsilon$ of $T_0$ that have a unique mixing ACIM $\mu_\varepsilon$ supported on the entire $I$. Metastability has a rich literature, we refer to \cite{random metastable} and references therein for an overview. In the context of chaotic dynamical systems, a systematic approach was initiated in \cite{KL}, which pointed out some key spectral techniques for studying metastable phenomena, relating these also to the theory of open dynamical systems (or systems with holes). Several further aspects of deterministic metastable systems were studied in the works \cite{BV13, FP12, FS10, K12}. Particularly
important for us is \cite{metastable}, which studied metastability in the above described setting of piecewise expanding $C^2$ interval maps and showed, under quite general conditions, that $\mu_\varepsilon$ converges to a convex combination of the measures $\mu_j$ when the  perturbation tends to zero.  Furthermore, the convex coefficients depend on the rates with which points flow between the  subintervals. In \cite{difcoef10} the authors showed that the time-rescaled process of the wandering of points approximates a jump Markov-process described by the above rates. They also proved that for all observables the rescaled diffusion coefficient of the observable converges to the diffusion coefficient of the induced observable defined on the state space of the Markov-chain. Recently, \cite{random metastable} and \cite{random diffusion} generalized these two papers to random metastable systems.

The main difference in the setting of metastable systems and the families studied in \cite{Keller82} and the above mentioned follow-up papers is that in the latter, one \emph{almost invariant} interval shrinks to a point. This introduces highly nontrivial technical challenges in the analysis of these systems. Yet, in this paper we aim to argue that the phenomenon of ACIM instability can be described by the natural theory of metastability in a large class of systems. We attempt to present the results and their proofs to simultaneously show the resemblance to their counterparts in the literature of metastable systems \cite{difcoef10, metastable}. Our results thus provide
a more systematic description of the convex combination of the limiting invariant measure, for a much wider class of families, generalizing the results on W shaped maps. Furthermore, we prove that the process of the flow of points between the shrinking (localized) and the extended (delocalized) almost invariant components of the dynamics, when suitably rescaled, converges to a jump Markov process, resembling the analogous results of \cite{difcoef10,random diffusion} on metastable systems.

The paper is structured as follows.
In Section \ref{section: conditions}, we introduce the class of systems which are in the scope of this paper. Here, and in the following sections we only discuss maps that have one fixed turning point, and leave the discussion about multiple turning points to the Appendix.
In Section \ref{section: results} we state our three main results analogous to those of \cite{difcoef10, metastable}.
In Section \ref{section: tools} we introduce the tools we will need in order to study the systems our framework accomodates. Most importantly, the \textit{conjugation}, which will allow us to compare our systems to those of metastable nature.
In Sections \ref{section: proof of theo 1}, \ref{section: proof of theo 2}, and \ref{section: proof of theo 3} we present the proofs of the three main results, respectively. There is a fundamental behavior of the maps that we take advantage of. Namely, after a point leaves the local state, it can be \emph{frozen} for an extended amount of time. This will allow us to wait for curves and densities, suffering extreme evolution leaving the local state, to regularize.
Utilizing this freezing effect, we can derive two key tools. On the one hand, we introduce the $\varepsilon$-variance and the associated norm with respect to which we can construct a Lasota-Yorke estimate in Section \ref{section: e-variance}; on the other hand, we state and prove our highly non-trivial Growth lemma in Section \ref{section: growth lemma}. Throughout the exposition, we aim to emphasize both the similarities and the differences of the phenomena (and thus of the required techniques) in our framework and in that of the metastable systems studied in \cite{metastable, difcoef10}.
In section \ref{section: examples} a new family of maps is presented as an example illustrating the applicability of our results. We then discuss both the generality of our conditions and the practical aspects of their verification.
In Section \ref{section: conclusions} we make some conclusions and remarks regarding the framework and the results of the paper. We also pose two questions regarding ACIM instability, one connecting this phenomenon to the \emph{shape} of the maps in the family, and the other connecting it to the \emph{emergence of the jump Markov behavior} described in Theorem \ref{theo: markov}. These questions aim to refine Keller's original conjecture of \cite{Keller82}.
In Appendix \ref{section: multiple} we present the naturally adapted conditions for systems with multiple turning points, and state the generalized results. We also discuss the slight modifications in the proofs needed to fit the framework of multiple turning points.

Throughout this paper, we will use the following notations. For two functions $f,g: \mathbb{R}^+ \rightarrow \mathbb{R}$, $f(x) \ll g(x)$ holds if there is a $C>0$ and an $x_0\in\mathbb{R}$ such that $|f(x)| \leq C g(x) $ for all $x>x_0$.
Two functions are asymptotic, if $f(x) \ll g(x) \ll f(x)$. We will denote this by $f(x) \asymp g(x)$.
The functions $g$ and $f$ are asymptotically equivalent, $f(x)\sim g(x)$, if $\lim_{x\to \infty}f(x)/g(x)=1$.

\section{Assumptions}\label{section: conditions}

\subsection*{The limit map}

\begin{enumerate}[label=(C\arabic*), start=0]
    \item\label{C0} We assume that we have a piecewise $C^2$ expanding map $\T: I\to I$,
        with $I=[0, 1]$ and critical points $C_0=\{c_1, \dots, c_m\}$.
        It follows from this that $\T$ has a distorsion bound, meaning, that there is a $\tau>0$ such that for any $n\in \N^+$ if $\T^n$ is $C^2$ smooth on $[x, y]$, then
        \begin{align}\label{eq: distorsion}
            \frac{\D\T^n(x)}{\D\T^n(y)}\leq \tau.
        \end{align}
        We assume that $\T$ has a unique ACIM $\hat\mu$ with density $\hat\phi$, and it is mixing.

    \item\label{C1} We assume that $\T$ has a single turning point $a\in int(I)$, meaning
    \begin{align*}
        \T(a)=a, \quad \partial_\pm\T (a)=\lim_{x \to a\pm}\D\T (x), \quad sign(\partial_-\T (a)) \neq sign(\partial_+\T (a)),
    \end{align*}
    where $sign(y)$ denotes the sign of a value $y\neq 0$, and $\T$ is continuous at $a$. See Fig \ref{fig: perturbation} for reference. Notice that this implies $a=c_{i_0}$ for some $i_0 \leq m$. See Appendix \ref{section: multiple} about maps having multiple turning points.

\end{enumerate}
As we will utilize methods from $\cite{difcoef10, metastable}$, we adopt some of their technical assumptions, which are labeled in \cite{metastable} by (I2)-(I4). First, we define \textit{infinitesimal holes}, the finite set of points  $H=\T^{-1}(a)\setminus \{a\}$.

\begin{enumerate}[label=(C\arabic*), start=2]
    \item\label{C2} We assume that there is no return of the critical set to the infinitesimal holes, hence to $a$: $a\notin \T^k (C_0\setminus \{a\})$ for all $k\geq 1$.

    \item\label{C3} We assume that $\hat\phi$ is positive on $a$.
\end{enumerate}

Both of these properties are generic. Condition \ref{C2} implies that $\hat\phi$ is continuous at each infinitesimal hole; see \cite[Section 4.2]{metastable}. Condition \ref{C3} implies that $\hat  \phi$ is positive at at least one point of the infinitesimal holes $H$, which is important for Theorems \ref{theo: markov} and \ref{theo: diffusion} but not for Theorem \ref{theo: convex}, see Remark \ref{remark: C3}.

\begin{enumerate}[label=(C\arabic*), start=4]
    \item\label{C4} There are no periodic critical points except for $a$ and possibly $\partial I$.
\end{enumerate}
This condition will be used to construct a uniform Lasota-Yorke inequality in Section \ref{section: lasota yorke}.

\subsection*{The perturbation}

\begin{figure}[H]
\begin{center}
\begin{tabular}{cc}
\begin{tikzpicture}
    \begin{scope}[scale=1.5]
        \draw[black] (0, 0) -- (0, 4) -- (4, 4) -- (4, 0) -- (0, 0);

        \node at (0.3, 3.7) {$\T$};

        \draw[lightgray] (0, 0) -- (4, 4);


        \draw[] (2.15, -0.1) -- (2.15, 0.1);
        \node at (2.15, -0.2) {$a$};

        \draw[] (-0.1, 2.15) -- (0.1, 2.15);
        \node at (-0.3, 2.35) {$a$};

        \draw[thick, decorate, decoration={snake, amplitude=0.3mm, segment length=12mm}] (2.15, 2.15) -- (2.7, 1.3);
        \draw[thick, decorate, decoration={snake, amplitude=0.3mm, segment length=8mm}]  (2.15, 2.15) -- (1.5, 1);

        \draw[thick, decorate, decoration={snake, amplitude=0.2mm, segment length=15mm}] (0, 1) -- (1, 2.5);
        \draw[thick, decorate, decoration={snake, amplitude=0.4mm, segment length=12mm}] (1, 1.6) -- (1.5, 0);
        \draw[thick, decorate, decoration={snake, amplitude=0.1mm, segment length=10mm}] (2.65, 4) -- (3.5, 2.7);
        \draw[thick, decorate, decoration={snake, amplitude=0.4mm, segment length=15mm}] (3.5, 2) -- (4, 2.8);
    \end{scope}
\end{tikzpicture}
&
\begin{tikzpicture}
    \begin{scope}[scale=1.5]
        \draw[black] (0, 0) -- (0, 4) -- (4, 4) -- (4, 0) -- (0, 0);

        \node at (0.3, 3.7) {$\T_\varepsilon$};

        \draw[lightgray] (0, 0) -- (4, 4);


        \draw[] (-0.05, 2.5) -- (0.05, 2.5);
        \node at (-0.3, 2.5) {$a_2$};

        \draw[] (2-0.2, -0.05) -- (2-0.2, 0.05);
        \node at (1.8, -0.2) {$a_1$};

        \draw[] (-0.05, 1.8) -- (0.05, 1.8);
        \node at (-0.3, 1.8) {$a_1$};

        \draw[] (2.15, -0.1) -- (2.15, 0.1);
        \node at (2.15, -0.22) {$a_\varepsilon$};

        \draw[] (2+0.5, -0.05) -- (2+0.5, 0.05);
        \node at (2.5, -0.2) {$a_2$};

        \draw[red] (1.8, 1.8) -- (2.5, 1.8) -- (2.5, 2.5) -- (1.8, 2.5) -- (1.8, 1.8);

        \draw[thick, decorate, decoration={snake, amplitude=0.3mm, segment length=12mm}] (2.15, 2.65) -- (2.7, 1.3);
        \draw[thick, decorate, decoration={snake, amplitude=0.3mm, segment length=8mm}]  (2.15, 2.65) -- (1.5, 1);

        \draw[thick, decorate, decoration={snake, amplitude=0.4mm, segment length=15mm}] (0, 1) -- (1, 3);
        \draw[thick, decorate, decoration={snake, amplitude=0.5mm, segment length=12mm}] (1, 1.7) -- (1.5, 0);
        \draw[thick, decorate, decoration={snake, amplitude=0.3mm, segment length=10mm}] (2.65, 4) -- (3.5, 2.5);
        \draw[thick, decorate, decoration={snake, amplitude=0.4mm, segment length=15mm}] (3.5, 1.9) -- (4, 2.9);
    \end{scope}
\end{tikzpicture}
\end{tabular}
\end{center}
\caption{The map $\T$ on the left and the perturbed map $\T_\varepsilon$ on the right. The box $\B$ is highlighted with red.}
\label{fig: perturbation}
\end{figure}
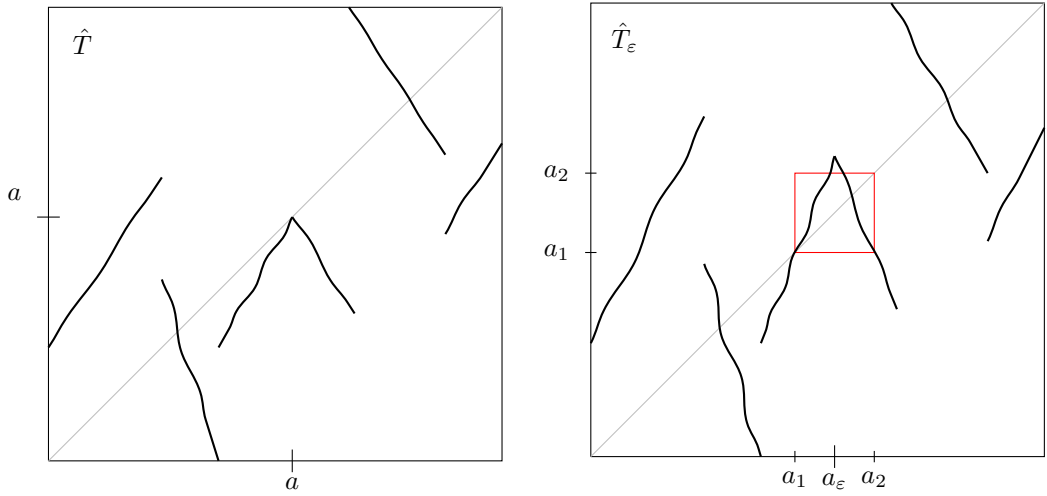

\begin{enumerate}[label=(C\arabic*), start=5]
    \item\label{C5} We consider a $C^2$ perturbation $\T_\varepsilon$,             $\varepsilon>0$ of $\T_0 = \T$.
        We assume that there is a family of piecewise expanding maps $\T_\varepsilon$ with critial points $C_\varepsilon=\{c_{\varepsilon, 1}, \dots, c_{\varepsilon, m} \}$ such that the functions $\varepsilon \mapsto c_{\varepsilon, i}$ are $C^2$ for all $i\leq m$.
        There exists a $\delta > 0$ such that $\T_{\varepsilon}|_{]c_{\varepsilon, i}, c_{\varepsilon, i+1}[}$ has a $C^2$ extension $\tilde T_{\varepsilon, i}: [c_i-\delta, c_{i+1}+\delta] \to \R$ such that $\|\tilde T_{\varepsilon, i} - \T_{0, i}\|_{C^2} \to 0$.
        We assume that for all $\varepsilon>0$,                       $\T_\varepsilon$ has a unique ACIM                            $\hat\mu_\varepsilon$ with density                            $\hat\phi_\varepsilon$, and it is mixing.

    \item\label{C6} We assume that $\T_\varepsilon$ is continuous at $a_\varepsilon = c_{\varepsilon, i_0}$ and satisfies
    \begin{align*}
        sign(\partial_-\T (a))\cdot \T_\varepsilon(a_\varepsilon)
        >
        sign(\partial_-\T (a))\cdot a_\varepsilon.
    \end{align*}

    This condition ensures that the graph of $\T_\varepsilon$ crosses the diagonal in the vicinity of the point $(a_\varepsilon, \T a_\varepsilon)$, see Fig \ref{fig: perturbation}.
\end{enumerate}

\begin{remark}
    Without loss of generality, we will assume that in (C1) $\partial_+\T (a) < 0 < \partial_-\T (a)$ holds, since in the other case one can study $U\circ \T_\varepsilon\circ U^{-1}$, where $U(x)=1-x$. Under this assumption, (C6) reduces to $\T_\varepsilon(a_\varepsilon)>a_\varepsilon.$
\end{remark}

Conditions (C2) and (C4) can be difficult to check for some systems. However, if $\T$ happens to be Markovian (has a finite Markov partition), the verification reduces to checking finitely many conditions. Notice that even if $\T$ is Markovian, this does not necessarily imply that $\T_\varepsilon$ is Markovian for any $\varepsilon>0$.

Notice that under these conditions, for sufficiently small $\varepsilon$, $\T_\varepsilon$ has a fixed point $a^{(\varepsilon)}_1 < a_\varepsilon$ which converges to $a$ as $\varepsilon\to 0$. Since $a$ is a turning point of $\T$, by assumption \ref{C6}, there is a another point $a^{(\varepsilon)}_2 > a_\varepsilon$ such that $\T_\varepsilon(a^{(\varepsilon)}_2)=\T_\varepsilon(a^{(\varepsilon)}_1)=a^{(\varepsilon)}_1$. We will call the interval $\mathcal{B}_\varepsilon=[a^{(\varepsilon)}_1, a^{(\varepsilon)}_2]$ the \textit{box}. For the ease of notation, we drop the $\varepsilon$ from these three, and we use $\B, a_1$ and $a_2$.

We can assume that $a_2 - a_1 \asymp \varepsilon $, and we know that the intervals $[a_1, a_\varepsilon]$ and $[a_\varepsilon, a_2]$ lie on separate branches of $\T_\varepsilon$.

\section{Results}\label{section: results}

As discussed in the introduction, one of the main points of our results is that although we start from a system that might experience instability, our setting shares some similarities with the metastable systems studied in the works \cite{difcoef10, metastable}. On the one hand, the presentation of the results aims to show these analogies. On the other hand, it is important to emphasize the new feature of our framework that there is a localized and a delocalized state.

\begin{remark}
    Throughout the paper, by the weak limit of measures we mean the weak limit w.r.t. the space of continuous functions $C(I)$ on $I$. Hence $\mu_\varepsilon \xrightarrow{w} \mu$ means that $\mu_\varepsilon(f) \to \mu(f)$ for all $f\in C(I)$.
\end{remark}

\subsection{Limit measure}

\begin{definition}
    We will refer to the sets
    \[
    \hat H_0(=\hat H_0(\varepsilon))=  (I \setminus \B) \cap \T_\varepsilon^{-1}\B, \qquad \hat H_1(=\hat H_1(\varepsilon))= \B \cap \T_\varepsilon^{-1}(I\setminus \B)
    \]
    as \textit{holes} (see Figure \ref{fig: holes}). The normalized Lebesgue measure on $\B$ is denoted as  $m_\B(A) = \frac{Leb(A\cap \B)}{Leb(\B)}$. Let us define the \textit{limit hole ratio} $l.h.r.$ as
\begin{equation*}
    l.h.r.
    =\lim_{\varepsilon\to 0}
    \frac{m_\B(\hat{H}_1)}{\hat{\mu}(\hat{H}_0)}.
\end{equation*}
\end{definition}

Consider the piecewise linear tent map  $F:[0,1]\to [0,1]$ with slopes of $s_1>0$ and $s_2<0$, that is
\[
F(t)=\begin{cases}
s_1\cdot t & \text{ if } 0\le t<t_0,\\
s_2\cdot(t-1) & \text{ if } t_0\le t\le 1.
\end{cases}
\]
with $t_0$ such that $F(t)$ is continuous. Then the \textit{top of the map} is at $h = F(t_0)=(|s_1|^{-1} + |s_2|^{-1})^{-1}$. This motivates introducing the following quantity:
\begin{align}\label{eq: h}
    h_\varepsilon = (|\partial_{-}\T_{\varepsilon}(a_\varepsilon)|^{-1} + |\partial_{+}\T_{\varepsilon}(a_\varepsilon)|^{-1})^{-1}, \quad \varepsilon \geq 0.
\end{align}
It turns out that the scaling of $h_\varepsilon$ determines if the ACIM of $\T_\varepsilon$ limits to the ACIM of $\hat T$, to a Dirac measure at $a$, or to a convex combination of these two, see Theorem~\ref{theo: convex} below. In the case of a nontrivial convex combination, the coefficients are given by the l.h.r. which in turn can be expressed in terms of some further quantities describing the perturbation, defined as:
\begin{align}\label{eq: cpm}
    c_\pm = \lim_{\varepsilon \to 0} \frac{\partial_{\pm}\T_{\varepsilon}(a_\varepsilon) - \partial_{\pm}\T(a)}{\varepsilon}
    =
    \partial_\varepsilon(\partial_{\pm}\T_{\varepsilon}(a_\varepsilon))|_{\varepsilon=0};
\end{align}

\begin{align*}
    c_a
    =&
    \lim_{\varepsilon\to 0}\frac{|\T_\varepsilon(a_\varepsilon)-a_\varepsilon|}{\varepsilon}.
\end{align*}

\begin{theorem}\label{theo: convex}
    If $\T_\varepsilon$ satisfies conditions (C0)-(C6), then                                                                                                                                                                                                                                                                                                                                                                                                                                                                                                                                                                              \begin{align*}
        \text{(I) If }& h_\varepsilon \leq 1 \; \forall \varepsilon>0, \text{ then } \hat\mu_\varepsilon \xrightarrow{w} \delta_a.
        \\
        \text{(II) If }& h_\varepsilon > 1 \; \forall \varepsilon\geq 0, \text{ then } \hat\phi_\varepsilon \xrightarrow{L^1} \hat\phi.
        \\
        \text{(III) If } h_\varepsilon > 1 \; \forall \varepsilon > 0& \text{ and } h_0 = 1,  \text{ then } \hat\mu_\varepsilon \xrightarrow{w} \alpha\hat\mu + (1-\alpha)\delta_a, \text{ where }
    \end{align*}
    \begin{align*}
        &\frac{\alpha}{1-\alpha}
        =
        l.h.r.
        =
        \frac{|c_-| \partial_{-}\T (a)^{-2} +  |c_+| \partial_{+}\T (a)^{-2}}
    {
    \sum_{x \in H}|\D \T x|^{-1}\hat\phi(x)
    \lim_{\varepsilon\to 0}\varepsilon^{-1}Leb(\B)}\text{ with}
     \\
      \lim_{\varepsilon\to 0}\varepsilon^{-1}Leb(\B)
     =&
     -c_a \left(
    (\partial_-\T(a) -  \mathbb{1}_{\partial_-\T(a)>0})(\partial_+\T(a) - \mathbb{1}_{\partial_+\T(a)>0})
    \right)^{-1}
    \left(
    |\partial_+\T(a)|
    +
    |\partial_-\T(a)|
    \right).
    \end{align*}
\end{theorem}

\begin{remark}
    Notice that if $\T$ satisfies the classic condition that the derivatives are uniformly separated from $2$, that is, $\inf_{I}|\D\T_\varepsilon|>\lambda>2$ for all $\varepsilon\geq 0$, then $h_0>1$, which means that by case (II), the family is ACIM stable.
\end{remark}

\begin{remark}\label{remark: C3}
    The only purpose of Condition \ref{C3} is to ensure that $l.h.r.<\infty$. However, the proof of Theorem \ref{theo: convex} does not utilize this assumption. If $\hat\phi(a)=0$, then one can see from the proof that in case (III), $\hat\phi_\varepsilon \xrightarrow{L^1} \hat\phi$.  Condition \ref{C3} will be important for Theorems \ref{theo: markov} and \ref{theo: diffusion}.
\end{remark}

\subsection{Jump Markov process}\label{section: process}

Beyond the description of the limit of the ACIM
$\hat{\mu}_{\varepsilon}$, in this section further results are included which express that, on a slow time scale,
the dynamics of $\T_{\varepsilon}$ converges to a jump Markov process. Theorems~\ref{theo: markov}
and~\ref{theo: diffusion} stated below are analogous to \cite[Theorems 1 and 2]{difcoef10}, yet, a specific feature of our setting is that there is
an extended and a localized metastable state, corresponding to $\hat{\mu}$ and $\delta_a$, respectively. Accordingly, these results are stated under the following additional condition:

\begin{enumerate}[label=(C\arabic*), start=7]
    \item\label{C7} We assume that case (III) of Theorem \ref{theo: convex} holds for the family $\T_\varepsilon$, $\varepsilon\geq 0$.
\end{enumerate}

If the family satisfies \ref{C7}, we can define the following Markov process, described in \cite{difcoef10} for metastable systems. Let the state space be $\mathcal{S}=\{0,1\}$, where $0$ corresponds to the extended state $I\setminus\B$, while $1$ corresponds to the localized sate $\B$. Define the jump rates as
\begin{align*}
    \beta_0
    =&
    \lim_{\varepsilon\to 0}\frac{\hat\mu(\hat H_0)}{\varepsilon}
    =
     \lim_{\varepsilon\to 0}\varepsilon^{-1}Leb(\B)
     \sum_{x \in H}|\D \T x|^{-1}\hat\phi(x),
    \\
    \beta_1
    =&
    \lim_{\varepsilon\to 0}\frac{m_\B(\hat H_1)}{\varepsilon}
    =
    |c_-| \partial_{-}\T (a)^{-2} +  |c_+| \partial_{+}\T (a)^{-2}.
\end{align*}

Then the process has a stationary distribution $\pi=(\pi_0, \pi_1)$ with $\pi_0/\pi_1(=\beta_1/\beta_0)=l.h.r.$ so that $\hat\mu_\varepsilon \xrightarrow{L^1} \pi_0 \hat\mu + \pi_1 \delta_a$.
Let $t_0 = 0$ and let $t_k^M$ be the $k$th time the process changes states. Let $\Time_k^M= t_k^M-t_{k-1}^M$ and let $z_k^M\in \mathcal{S}$ denote the state of the process after the $k$th transition. As in \cite{difcoef10}, for $r\in \mathcal{S}=\{0,1\}$ let $\mathbb{P}^r$ denote the probability measure on $\{0,1\}^{[0,\infty)}$ describing the process when started from state $r\in \mathcal{S}$, and evolved by the Markov dynamics. In particular, for any $k\ge 1$:
\begin{align*}
    &\text{d}\mathbb{P}^r(\Time_{2k-1}^M = t ) = \beta_{r} e^{-\beta_{r}t}dt, \qquad &\mathbb{P}^r(z_{2k-1}^M = 1-r) = 1;\\
    &\text{d}\mathbb{P}^r(\Time_{2k}^M = t ) = \beta_{1-r} e^{-\beta_{1-r}t}dt, \qquad &\mathbb{P}^r(z_{2k}^M = r) = 1.
\end{align*}

We can also describe a process generated by the system $\T_\varepsilon$. Let $z(x)=0$ for $x\in I\setminus\B$ and $z(x)=1$ for $x\in\B$. Let $t^\varepsilon_0 = 0$ and $t^\varepsilon_k = \inf \{t>t^\varepsilon_{k-1}\; | \; z(\T_\varepsilon^{t} x) \neq z(\T_\varepsilon^k x)\}$.

Theorem~\ref{theo: markov} establishes that in case (III), at the level of finite dimensional distributions, the dynamics of wandering between the box and the extended state can be approximated by a Markov-process.

\begin{theorem}\label{theo: markov}
    Let $\T_\varepsilon$ satisfy the conditions (C0)-(C7). Fix an integer $p\ge 1$ and two real parameters $0<s<S$. For any intervals $\Delta_k=[a_k, b_k]$, and numbers $r_k \in \mathcal{S}$, $k = 1, \dots, p$
    \begin{align*}
        \hat\mu(\varepsilon\Time^{\varepsilon}_{k}\in\Delta_k,\; z(t_k^{\varepsilon})
        =&
        r_k \; \forall k\leq p)
        \to \mathbb{P}^0(\Time_{k}^M\in\Delta_k,\; z_k^M=r_k \; \forall k\leq p),
        \\
        m_\B(\varepsilon\Time^{\varepsilon}_{k}\in\Delta_k,\; z(t_k^{\varepsilon})
        =&
        r_k \; \forall k\leq p)
        \to \mathbb{P}^1(\Time_{k}^M\in\Delta_k,\; z_k^M=r_k \; \forall k\leq p),
    \end{align*}

and the convergence is uniform for $\max_k b_k \leq S$ and $\min_k a_k\geq s$.
\end{theorem}

Our next aim is to show that this correspondence extends to the level of the asymptotic properties of the processes by studying the CLT in the perturbed system $\T_{\varepsilon}$. Our result, Theorem~\ref{theo: diffusion} is analogous to \cite[Theorem 2]{difcoef10}.

Let $X: I\to \R$ be an observable. For any fixed, sufficiently small $\varepsilon > 0$, the Ergodic theorem provides a law of large numbers, hence $N^{-1}\sum_{n=0}^{N-1}X\circ \T_\varepsilon^{n} \to \hat\mu_\varepsilon(X)$ as $N\to \infty$,  almost surely and in $L^1$ (w.r.t. $\hat\mu_\varepsilon$). If $X\in BV(I)$ is centered, then by \cite{CLT}, $N^{-1/2}\sum_{n=0}^{N-1}X\circ \T_\varepsilon^{n}$ converges to a normal distribution as $N\to \infty$. The variance $D^\varepsilon(X)$ of this limiting distribution is called the \textit{diffusion coefficient}.

For a $X\in BV(I)$ which is continuous at $a$, we can define an observable $\X: \mathcal{S}\to \R$ on the state space of the Markov-process described above by letting  $\X(0)=\hat\mu(X)$ and $\X(1) = \delta_a(X)=X(a)$. We may assume that $X$ is centered with respect to the limit measure, that is, $\pi_0\X(0)+\pi_1\X(1)=0$, that is, $\beta_1
\hat\mu(X)=-\beta_0 X(a)$. Let $\bm{D}(\X)$ denote the diffusion coefficient of $\X$.

\begin{theorem}\label{theo: diffusion}
    If $\T_\varepsilon$ satisfies conditions (C0)-(C7), then
    \begin{align*}
        \varepsilon D^{\varepsilon}(X) \to \textbf{D}(\X)
    \end{align*}
    holds for any centered observable $X\in BV(I)$ which is continuous at $a$.
\end{theorem}

\begin{remark}
    Note that $\bm{D}(\X)$ is easily computable. In fact, $\bm{D}(\X)=\langle \pi \X , G^{-1} \X \rangle$, where $\pi = (\pi_0, \pi_1)$ is the stationary distribution and $G$ is the infinitesimal generator of the Markov process described above. From this, one arrives at $\bm{D}(\X)=2\beta_0\beta_1\frac{(\X(0) - \X(1))^2}{(\beta_0 + \beta_1)^3}$.
\end{remark}

\section{The tools}\label{section: tools}

\subsection{The conjugation}\label{section: conjugation}

\begin{figure}[H]
\centering
\begin{tabular}{c c}
\begin{tikzpicture}
    \begin{scope}[scale=1.5]
        \draw[black] (0, 0) -- (0, 4) -- (4, 4) -- (4, 0) -- (0, 0);

        \node at (0.3, 3.7) {$\T$};

        \draw[lightgray] (0, 0) -- (4, 4);


        \draw[] (2.15, -0.1) -- (2.15, 0);
        \node at (2.15, -0.2) {$a$};

        \draw[white] (0, -0.48) -- (0.1, -0.48);

        \draw[thick, decorate, decoration={snake, amplitude=0.5, segment length=25}] (2.15, 2.15) -- (2.71, 1.28);
        \draw[thick, decorate, decoration={snake, amplitude=0.5, segment length=20}]  (2.15, 2.15) -- (1.5, 1.1);

        \draw[thick, decorate, decoration={snake, amplitude=0.6, segment length=30}] (0, 1) -- (1, 3);
        \draw[thick, decorate, decoration={snake, amplitude=0.5mm, segment length=15mm}] (1, 1.7) -- (1.5, 0);
        \draw[thick, decorate, decoration={snake, amplitude=0.6mm, segment length=15mm}] (2.71, 4) -- (3.5, 2.5);
        \draw[thick, decorate, decoration={snake, amplitude=0.4mm, segment length=20mm}] (3.5, 1.9) -- (4, 2.9);

    \end{scope}
\end{tikzpicture}
&
\begin{tikzpicture}
    \begin{scope}[scale=1.5]
        \draw[black] (0, 0) -- (0, 4) -- (4, 4) -- (4, 0) -- (0, 0);

        \node at (0.3, 3.7) {$\T_\varepsilon$};

        \draw[lightgray] (0, 0) -- (4, 4);


        \draw[] (-0.05, 2.5) -- (0.05, 2.5);
        \node at (-0.2, 2.5) {$a_2$};

        \draw[] (-0.05, 1.8) -- (0.05, 1.8);
        \node at (-0.2, 1.8) {$a_1$};

        \draw[] (2.5, -0.05) -- (2.5, 0.05);
        \node at (2.5, -0.2) {$a_2$};

        \draw[] (1.8, -0.05) -- (1.8, 0.05);
        \node at (1.8, -0.2) {$a_1$};

        \node[blue] at (0.59, -0.26) {$\hat{H}_0^{(1)}$};

        \node[teal] at (2.15, -0.26) {$\hat{H}_1$};

        \node[blue] at (3.65, -0.26) {$\hat{H}_0^{(2)}$};

        \draw[red] (1.8, 1.8) -- (2.5, 1.8) -- (2.5, 2.5) -- (1.8, 2.5) -- (1.8, 1.8);

        \draw[thick, decorate, decoration={snake, amplitude=0.5, segment length=25}] (2.15, 2.65) -- (2.71, 1.28);
        \draw[thick, decorate, decoration={snake, amplitude=0.5, segment length=20}]  (2.15, 2.65) -- (1.5, 1.1);

        \draw[thick, decorate, decoration={snake, amplitude=0.6, segment length=30}] (0, 1) -- (1, 3);
        \draw[thick, decorate, decoration={snake, amplitude=0.5mm, segment length=15mm}] (1, 1.7) -- (1.5, 0);
        \draw[thick, decorate, decoration={snake, amplitude=0.6mm, segment length=15mm}] (2.71, 4) -- (3.5, 2.6);
        \draw[thick, decorate, decoration={snake, amplitude=0.4mm, segment length=20mm}] (3.5, 1.7) -- (4, 2.9);


        \draw[blue, line width=0.5mm] (0.41, 0) -- (0.74, 0);
        \draw[blue, line width=0.5mm] (3.5, 0) -- (3.8, 0);

        \draw[teal, line width=0.5mm] (2.075, 0) -- (2.22, 0);

    \end{scope}
\end{tikzpicture}
\end{tabular}
\caption{The map $\T$ on the left and the perturbed map $\T_\varepsilon$ on the right. \textit{The box} is in red and the holes are colored.}\label{fig: holes}
\end{figure}
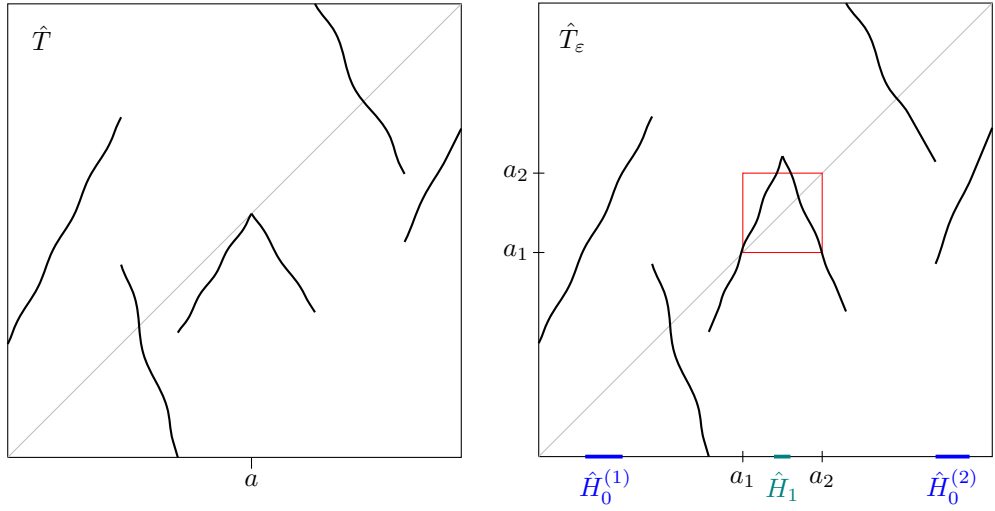

Notice that our framework has some similarities to that of metastable systems in \cite{metastable, difcoef10}. Fist, notice that the maps $\T_\varepsilon$ have two almost invariant sets, $\B$ and $I \setminus \B$.
See Fig~\ref{fig: holes} for an illustration of the holes $\hat H_0 = (I\setminus\B) \cap \Te^{-1}\B$ and $\hat H_1 = \B \cap \Te^{-1}(I\setminus \B)$ -- specifically for this example, $\hat H_0=\hat H_0^{(1)} \cup \hat H_0^{(2)}$ has two connected components.

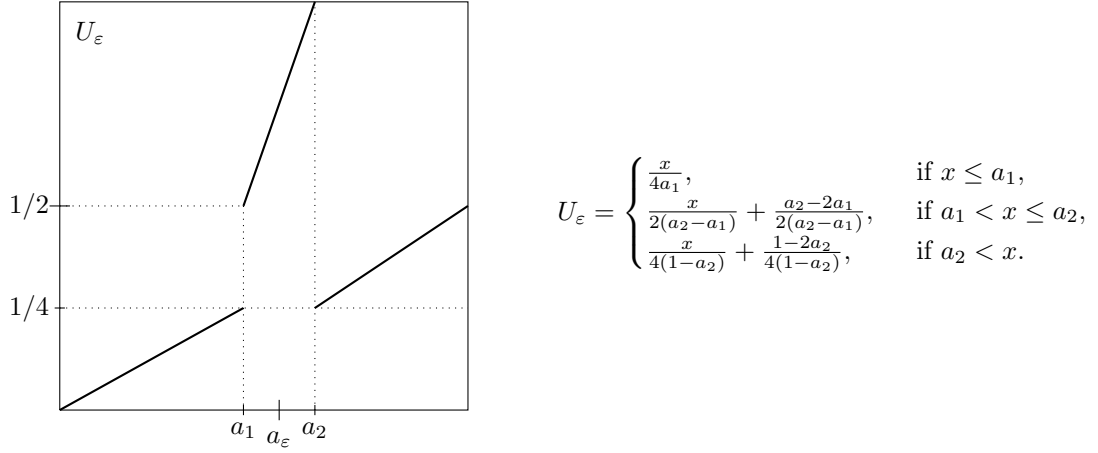
\begin{figure}[H]
\centering
\begin{subfigure}{0.4\textwidth}

    \begin{tikzpicture}
        \begin{scope}[scale=1.35]
        \draw[black] (0, 0) -- (0, 4) -- (4, 4) -- (4, 0) -- (0, 0);

        \draw[black, thick] (2-0.2, 2) -- (2+0.5, 4);
        \draw[black, thick] (0, 0) -- (2-0.2, 1);
        \draw[black, thick] (2+0.5, 1) -- (4, 2);

        \draw[dotted] (0, 2) -- (2-0.2, 2);
        \draw[dotted] (2-0.2, 2) -- (2-0.2, 0);
        \draw[dotted] (2+0.5, 0) -- (2+0.5, 4);
        \draw[dotted] (0, 1) -- (4, 1);

        \draw[] (-0.1, 2) -- (0.1, 2);
        \node at (-0.3, 2) {1/2};
        \draw[] (2-0.2, -0.05) -- (2-0.2, 0.05);
        \draw[] (2.15, -0.1) -- (2.15, 0.1);
        \draw[] (2+0.5, -0.05) -- (2+0.5, 0.05);
        \draw[] (-0.05, 1) -- (0.05, 1);
        \node at (-0.3, 1) {1/4};

        \node at (2.15, -0.29) {$a_\varepsilon$};
        \node at (1.8, -0.2) {$a_1$};
        \node at (2.5, -0.2) {$a_2$};

        \node at (0.3, 3.7) {$U_\varepsilon$};
        \end{scope}
    \end{tikzpicture}
    \end{subfigure}
    \begin{subfigure}{0.4\textwidth}
        \begin{tikzpicture}
            \node at (0,0) {$\;$};
            \node at (4, 3.1)
            {
            $
            U_\varepsilon
            =
            \begin{cases}
                \frac{x}{4a_1}, \; &\text{ if } x\leq a_1,
                \\
                \frac{x}{2(a_2 - a_1)} +  \frac{a_2 -2 a_1}{2(a_2-a_1)}, \; &\text{ if } a_1 < x \leq a_2,
                \\
                \frac{x}{4(1-a_2)} +  \frac{1 - 2 a_2}{4(1-a_2)}, \; &\text{ if } a_2 < x .
            \end{cases}
            $
            };
        \end{tikzpicture}
    \end{subfigure}
    \caption{ The graph of the conjugation map $U_\varepsilon$.}\label{fig: conjugation}
\end{figure}

To further motivate the connection between our framework and that of metastable systems, we introduce a conjugation $U_{\varepsilon}$ in Fig \ref{fig: conjugation}, which is a piecewise linear bijection of $I$ onto itself, such that the box is mapped onto $I_R = [1/2, 1]$ and the outside of the box is mapped onto $I_L = [0, 1/2[$.

Using $U_\varepsilon$ we define the map
\begin{align*}
    T_{\varepsilon} = U_\varepsilon \circ \T_\varepsilon \circ U_\varepsilon^{-1}.
\end{align*}

\begin{figure}[H]
\centering
\begin{subfigure}{0.4\textwidth}
\begin{tikzpicture}
    \begin{scope}[scale=1.5]
        \draw[black] (0, 0) -- (0, 4) -- (4, 4) -- (4, 0) -- (0, 0);

        \node at (0.3, 3.7) {$\T_\varepsilon$};

        \draw[lightgray] (0, 0) -- (4, 4);

        \draw[black] (2.15, -0.1) -- (2.15, 0);
        \node at (2.15, -0.2) {$a_\varepsilon$};

        \draw[white] (0, -0.5) -- (0.1, -0.5);

        \draw[red] (1.8, 1.8) -- (2.5, 1.8) -- (2.5, 2.5) -- (1.8, 2.5) -- (1.8, 1.8);

        \draw[thick, decorate, decoration={snake, amplitude=0.5, segment length=25}] (2.15, 2.65) -- (2.71, 1.28);
        \draw[thick, decorate, decoration={snake, amplitude=0.5, segment length=20}]  (2.15, 2.65) -- (1.5, 1.1);

        \draw[thick, decorate, decoration={snake, amplitude=0.6, segment length=30}] (0, 1) -- (1, 3);
        \draw[thick, decorate, decoration={snake, amplitude=0.5mm, segment length=15mm}] (1, 1.7) -- (1.5, 0);
        \draw[thick, decorate, decoration={snake, amplitude=0.6mm, segment length=15mm}] (2.71, 4) -- (3.5, 2.6);
        \draw[thick, decorate, decoration={snake, amplitude=0.4mm, segment length=20mm}] (3.5, 1.7) -- (4, 2.9);

    \end{scope}
\end{tikzpicture}
\end{subfigure}
\begin{subfigure}{0.4\textwidth}
\begin{tikzpicture}
    \begin{scope}[scale=1.5]
        \draw[black] (0, 0) -- (0, 4) -- (4, 4) -- (4, 0) -- (0, 0);

        \node at (0.6, 3.7) {$T_\varepsilon$};

        \draw[dotted] (0, 2) -- (4, 2);
        \draw[dotted] (2, 0) -- (2, 4);
        \draw[dotted] (0, 1) -- (4, 1);

        \draw[] (-0.1, 2) -- (0, 2);
        \node at (-0.3, 2) {\small{$1/2$}};

        \draw[] (-0.1, 1) -- (0, 1);
        \node at (-0.3, 1) {\small{$1/4$}};

        \draw[] (2, -0.1) -- (2, 0);
        \node at (2, -0.2) {\small{$1/2$}};

        \draw[thick, decorate, decoration={snake, amplitude=0.3mm, segment length=28mm}] (2, 2) -- (2.8, 4);
        \draw[thick, decorate, decoration={snake, amplitude=0.3mm, segment length=18mm}] (3.3, 4) -- (4, 2);

        \draw[thick, decorate, decoration={snake, amplitude=0.3mm, segment length=8mm}]  (0, 0.48) -- (0.25, 1);

        \draw[blue, thick, decorate, decoration={snake, amplitude=0.3mm, segment length=22mm}]  (0.25, 2) -- (0.35, 4);

        \draw[thick, decorate, decoration={snake, amplitude=0.3mm, segment length=12mm}] (0.35, 1) -- (0.5, 1.4);

        \draw[thick, decorate, decoration={snake, amplitude=0.3mm, segment length=8mm}]  (0.5, 0.9) -- (0.8, 0);

        \draw[thick, decorate, decoration={snake, amplitude=0.3mm, segment length=8mm}]  (0.8, 0.6) -- (1, 1);

        \draw[thick, decorate, decoration={snake, amplitude=0.3mm, segment length=8mm}]  (1, 1) -- (1.15, 0.75);

        \draw[thick, decorate, decoration={snake, amplitude=0.3mm, segment length=8mm}]  (1.15, 2) -- (1.6, 1.05);

        \draw[thick, decorate, decoration={snake, amplitude=0.3mm, segment length=8mm}]  (1.6, 0.9) -- (1.65, 1);

        \draw[blue, thick, decorate, decoration={snake, amplitude=0.3mm, segment length=20mm}]  (1.65, 2) -- (1.75, 4);

        \draw[thick, decorate, decoration={snake, amplitude=0.3mm, segment length=8mm}]  (1.75, 1) -- (2, 1.4);

        \draw[teal, thick, decorate, decoration={snake, amplitude=0.3mm, segment length=8mm}] (2.8, 1) -- (3.05, 1.1);
        \draw[teal, thick, decorate, decoration={snake, amplitude=0.3mm, segment length=8mm}] (3.05, 1.1) -- (3.3, 1);

        \draw[blue, line width=0.5mm] (0.25, 0) -- (0.35, 0);
        \draw[blue, line width=0.5mm] (1.65, 0) -- (1.75, 0);
        \node at (1.05, -0.3) {\textcolor{blue}{$H_0$}};

        \draw[teal, line width=0.5mm] (2.8, 0) -- (3.3, 0);
        \node at (3.05, -0.3) {\textcolor{teal}{$H_1$}};

    \end{scope}
\end{tikzpicture}
\end{subfigure}
\caption{The perturbed map $\T_\varepsilon$ and the conjugated perturbed map $T_\varepsilon$ with the holes colored.}\label{fig: conjugated holes}
\end{figure}
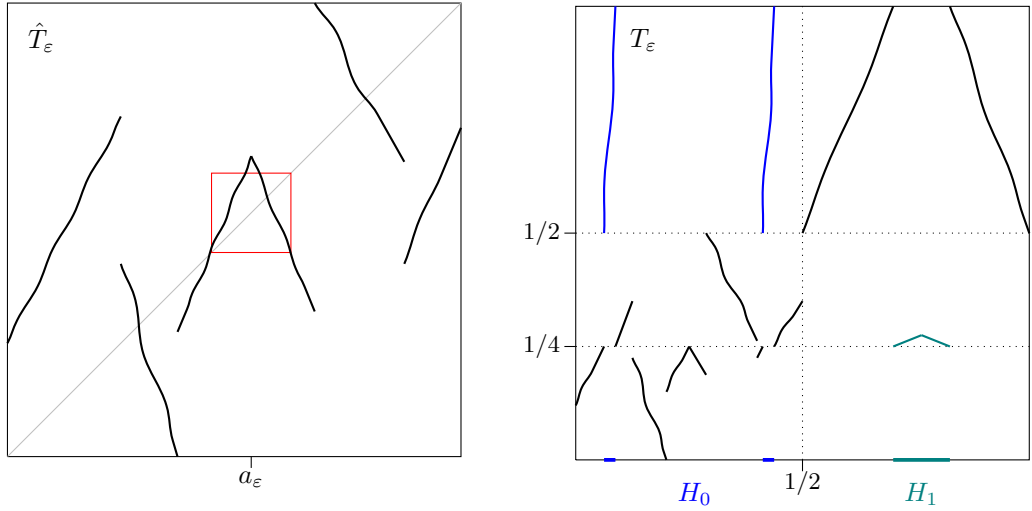

The map $T_\varepsilon$ has a unique ACIM $\mu_{\varepsilon} = U_*\hat{\mu}_{\varepsilon}$ with density $\phi_{\varepsilon}$.
We define the holes $H_0= (T_\varepsilon^{-1}I_R) \cap I_L$ and $H_1= (T_\varepsilon^{-1}I_L) \cap I_R$, see Fig \ref{fig: conjugated holes}.

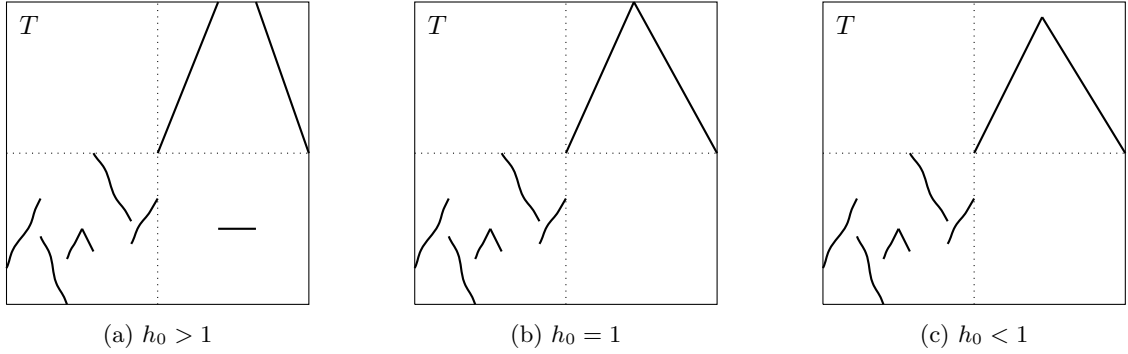
\begin{figure}[H]
    \begin{subfigure}{0.32\textwidth}
    \centering
        \begin{tikzpicture}
    \begin{scope}[scale=1]
        \draw[black] (0, 0) -- (0, 4) -- (4, 4) -- (4, 0) -- (0, 0);

        \node at (0.3, 3.7) {$T$};

        \draw[dotted] (0, 2) -- (4, 2);
        \draw[dotted] (2, 0) -- (2, 4);

        \draw[thick] (2, 2) -- (2.8, 4);
        \draw[thick] (3.3, 4) -- (4, 2);

        \draw[thick] (2.8, 1) -- (3.3, 1);

        \draw[thick, decorate, decoration={snake, amplitude=0.3mm, segment length=8mm}]  (0, 0.48) -- (0.45, 1.4);
        \draw[thick, decorate, decoration={snake, amplitude=0.3mm, segment length=8mm}]  (0.45, 0.9) -- (0.8, 0);
        \draw[thick, decorate, decoration={snake, amplitude=0.1mm, segment length=4mm}]  (0.8, 0.6) -- (1, 1);
        \draw[thick, decorate, decoration={snake, amplitude=0.3mm, segment length=8mm}]  (1, 1) -- (1.15, 0.7);
        \draw[thick, decorate, decoration={snake, amplitude=0.3mm, segment length=8mm}]  (1.15, 2) -- (1.65, 1.1);
        \draw[thick, decorate, decoration={snake, amplitude=0.3mm, segment length=8mm}]  (1.65, 0.8) --  (2, 1.4);

    \end{scope}
\end{tikzpicture}

    \subcaption{$h_0>1$}

    \end{subfigure}
    \begin{subfigure}{0.32\textwidth}
    \centering
        \begin{tikzpicture}
    \begin{scope}[scale=1]
        \draw[black] (0, 0) -- (0, 4) -- (4, 4) -- (4, 0) -- (0, 0);

        \node at (0.3, 3.7) {$T$};

        \draw[dotted] (0, 2) -- (4, 2);
        \draw[dotted] (2, 0) -- (2, 4);

        \draw[thick] (2, 2) -- (2.9, 4);
        \draw[thick] (2.9, 4) -- (4, 2);

        \draw[thick, decorate, decoration={snake, amplitude=0.3mm, segment length=8mm}]  (0, 0.48) -- (0.45, 1.4);
        \draw[thick, decorate, decoration={snake, amplitude=0.3mm, segment length=8mm}]  (0.45, 0.9) -- (0.8, 0);
        \draw[thick, decorate, decoration={snake, amplitude=0.1mm, segment length=4mm}]  (0.8, 0.6) -- (1, 1);
        \draw[thick, decorate, decoration={snake, amplitude=0.3mm, segment length=8mm}]  (1, 1) -- (1.15, 0.7);
        \draw[thick, decorate, decoration={snake, amplitude=0.3mm, segment length=8mm}]  (1.15, 2) -- (1.65, 1.1);
        \draw[thick, decorate, decoration={snake, amplitude=0.3mm, segment length=8mm}]  (1.65, 0.8) --  (2, 1.4);

    \end{scope}
\end{tikzpicture}

        \subcaption{$h_0=1$}

    \end{subfigure}
    \begin{subfigure}{0.32\textwidth}
    \centering
        \begin{tikzpicture}
    \begin{scope}[scale=1]
        \draw[black] (0, 0) -- (0, 4) -- (4, 4) -- (4, 0) -- (0, 0);

        \node at (0.3, 3.7) {$T$};

        \draw[dotted] (0, 2) -- (4, 2);
        \draw[dotted] (2, 0) -- (2, 4);

        \draw[thick] (2, 2) -- (2.9, 3.8);
        \draw[thick] (2.9, 3.8) -- (4, 2);

        \draw[thick, decorate, decoration={snake, amplitude=0.3mm, segment length=8mm}]  (0, 0.48) -- (0.45, 1.4);
        \draw[thick, decorate, decoration={snake, amplitude=0.3mm, segment length=8mm}]  (0.45, 0.9) -- (0.8, 0);
        \draw[thick, decorate, decoration={snake, amplitude=0.1mm, segment length=4mm}]  (0.8, 0.6) -- (1, 1);
        \draw[thick, decorate, decoration={snake, amplitude=0.3mm, segment length=8mm}]  (1, 1) -- (1.15, 0.7);
        \draw[thick, decorate, decoration={snake, amplitude=0.3mm, segment length=8mm}]  (1.15, 2) -- (1.65, 1.1);
        \draw[thick, decorate, decoration={snake, amplitude=0.3mm, segment length=8mm}]  (1.65, 0.8) --  (2, 1.4);

    \end{scope}
\end{tikzpicture}

    \subcaption{$h_0<1$}

    \end{subfigure}

\caption{The dependence of $T_0$ on $h_0$.}\label{fig: h0}
\end{figure}

We define $T=T_0$ as the pointwise limit of $T_\varepsilon$ when $\varepsilon\to 0$. We can see a qualitative dependence of $T$ on $h_0$ (defined in \eqref{eq: h}) in Fig \ref{fig: h0}

From now on, we will work with the conjugated map.
Notice that the map $T_\varepsilon$ only experiences anomalous behavior in the holes. It can be easily computed that the derivative of $T_\varepsilon$ is asymptotic to $\varepsilon^{-1}$ on $H_0$ and to $\varepsilon$ on $H_1$. Furthermore, its distorsion\footnote{Given an interval $J$ and a smooth monotone map $f:J\to\mathbb{R}$, by distorsion of $f$ on $J$ we mean the supremum of $\frac{|f'(x_1)|}{|f'(x_2)|}$ for $x_1,x_2\in J$.} is asymptotic to $\varepsilon$ on $H_0$, on $H_1$ and on the connected components of  $I_R\setminus H_1$.

\subsection{Lasota-Yorke inequalities and spectral consequences}\label{section: lasota-yorke intro}

\subsubsection{For a single map}\label{section: lasota-yorke for a single map}

Here we recall briefly some important facts on transfer operators of piecewise expanding maps of the interval and on their spectral properties when acting on functions of bounded variation, see for example \cite{laws of chaos}.

For a function $f:I\to \R$ its \emph{total variation} is defined by
\begin{equation*}
    \V{f} = \sup \left\{\sum_{i=0}^{n-1}|f(t_{i+1}) - f(t_{i})|\; : n\in \N^+, \; \{t_0, \dots, t_n\} \text{ is a partition of } I \right\}.
\end{equation*}
The $f:I\to \R$ has \emph{bounded variation}, if
\begin{equation*}
    \BV{f} = \V{f} + \Lnorm{f} = \V{f} + \int_I |f|
\end{equation*}
is finite. The space of functions with bounded variation endowed with the $\BV{\cdot}$ norm is denoted by $BV(I)\subset L^1(I)$, the unit ball of which is compact in $L^1(I)$.
Let $T:I \to I$ be a piecewise expanding map of the interval $I$. Then the associated Perron-Frobenius operator $\perron$ is defined as

\begin{align*}
\perron:& BV(I) \to BV(I),
\quad
&&\perron f(x) = \sum_{y \in T^{-1} \{x\}} \frac{f(y)}{|\D T(y)|}.
\end{align*}
This is the unique operator that satisfies
\begin{align*}
    \int_{T^{-1}A} f = \int_{A}\perron f
\end{align*}
for all $A\subset I$ and $f\in BV(I)$. We are interested in the fixed points of $\perron$, which are invariant densities for $T$.

If there is an $n_0 \in \N^+$, an $r<1$ and a $K>0$ such that

\begin{align}\label{eq: one-step lasota}
    \V{\perron^{n_0}f} \leq r\V{f} + K\Lnorm{f}
\end{align}
holds for all $f\in BV$, then there are finite number of points $\lambda_i$, $i\leq n$ of modulus one in the spectrum of $\perron$, all with finite eigen-spaces, and $1$ is an eigenvalue. This implies that there are finite rank projections $P_i$ and an operator $Q$ that are mutually orthogonal such that
\begin{align*}
    \perron = \sum_{i\leq n}\lambda_{i} P_i + Q.
\end{align*}
It also follows that $Q$ has norm less, than $1$.
Inequality \eqref{eq: one-step lasota} is equivalent with the statement that there is an $\alpha<1$ and a $C>0$ such that

\begin{align}
    \V{\perron^n f} \leq C(\alpha^n \V{f} + \Lnorm{f})
\end{align}
holds for all $n\in \N^+$ and $f\in BV(I)$. This inequality implies that the norm of $R$ is at most $\alpha$.

The fixed point $f\in BV$ of $\perron$ corresponds to the density of an invariant measure. The mixing of the map is equivalent to $1$ being a simple eigenvalue.

In \cite{Ionescu-Marinescu} and \cite{Hennion} the authors developed a general framework for studying such operators, but the above framework is sufficient for our purposes.

\subsubsection{The construction of the estimate}\label{section: lasota-yorke methodology}

From \cite[Proposition 3.1]{Stronger Lasota-Yorke}, we can see that a map $T$ satisfying (C1), the Lasota-Yorke type estimate \eqref{eq: one-step lasota} with $n_0=1$ depends on the following attributes of the map. A simple estimate for the contraction coefficient $r$ is $2\max_{I\setminus C_0}|\D T|^{-1}$, while $C$ depends on the ratio $|\D^2 T|/|\D T|^2$ and on the size of the smallest continuity interval of $\D T$. However, if $T$ only has surjective branches, then $r = \max|\D T|^{-1}$ and $C=\max|\D^2 T|/|\D T|^2$.

\subsubsection{Family of maps}\label{section: lasota-yorke for a family}

In the same framework of the previous section, we assume that there is a family of mixing maps $T_\varepsilon$, $\varepsilon\geq0$ with Perron-Frobenius operators $\perron_\varepsilon$ satisfying Lasota-Yorke type inequalities,

\begin{align}\label{eq: Lasota-Yorke}
    \V{\perron_\varepsilon^n f} \leq C_\varepsilon(\alpha_\varepsilon^n \V{f} + \Lnorm{f}), \quad \varepsilon \geq0, \, f\in BV, \, n\in \N^+.
\end{align}

When there is an $\alpha<1$ and $C>0$ such that $\alpha_\varepsilon<\alpha$ and $C_\varepsilon<C$ for all $\varepsilon>0$, we say that the operators satisfy a \emph{uniform Lasota-Yorke inequality},
\begin{align*}
    \V{\perron_\varepsilon^n f} \leq C(\alpha^n \V{f} + \Lnorm{f}), \quad \varepsilon \geq0, \, f\in BV, \, n\in \N^+.
\end{align*}
 This implies the ACIM stability of the family $T_\varepsilon$, $\varepsilon \geq 0$ by \cite{KL99}.

However, there are two inconvenient possibilities. It could happen, that $\alpha_\varepsilon\to 1$, or $C_\varepsilon \to \infty$.

The former means that the spectral gaps $1-\|R_\varepsilon\|$ are not bounded from below. In this case, even if every $\perron_\varepsilon$ has $1$ as a simple eigenvalue, it could happen that there is a family of eigenfunctions $f_\varepsilon \in BV$ with $\Lnorm{f_\varepsilon}=1$ corresponding to eigenvalues $\beta_\varepsilon<1$ such that $\beta_\varepsilon\to 1$. In this case, since $\perron_0$ also has $1$ as a simple eigenvalue, $f_\varepsilon$ might not have a limit in $L^1$.

If $C_\varepsilon \to \infty$, it can happen that the invariant densities $f_\varepsilon$ have variational norms converging to infinity, hence they might not have a limit in $L^1$.

In both cases however, the family of measures $\nu_\varepsilon$ absolutely continuous with respect to the Lebesgue measure (with densities $f_\varepsilon$) may have a weak limit $\nu$ which is, however, not absolutely continuous with respect the Lebesgue measure. That is, the family is not ACIM stable.

It follows from the form of the estimates mentioned in the previous section and detailed in \cite{Stronger Lasota-Yorke}, that if the family $T_\varepsilon$ has sufficient lower bound on the derivatives $|\D T_\varepsilon|$, it admits a uniform Lasota-Yorke inequality. Otherwise, one has to study the family $T_\varepsilon^{n_0(\varepsilon)}$, which has adequate derivative bounds. When the family $\T_\varepsilon$, $\varepsilon>0$ is such that $\inf_\varepsilon\inf \D\T>\lambda>1$, then there is a uniform $n_0\in \N^+$ such that $n_0(\varepsilon)\leq n_0$.  However, when constructing the estimate for $\T_\varepsilon^n$, one has to examine the length of the shortest continuity interval of $T_\varepsilon^{n_0}$ to estimate $C_\varepsilon$. If $T_0$ has a periodic  critical point other than possibly $0$ or $1$, this length might converge to zero, resulting in $C_\varepsilon\to \infty$.

Notice that if a $C^2$ family has sufficiently steep branches,  $\inf_{I}|\D\T_\varepsilon|>\lambda>2$ for all $\varepsilon\geq 0$, then one can directly construct the estimate \eqref{eq: one-step lasota} for $T_\varepsilon$ (meaning $n_0=1)$. This condition can be significantly relaxed, see \cite[Theorem 4.5]{Stronger Lasota-Yorke}.
Otherwise, if the limit map $T$ has no critical periodic points other than possibly $0$ or $1$, there is an $n_0\in \N^+$ such that one can construct \eqref{eq: one-step lasota} for the family.
In both of these cases, the family is ACIM stable.

Notice, however, that the existence of the turning fixed point $a$ defined in \ref{C1} directly violates the nonexistence of critical periodic points. By \ref{C6}, for any $n_0\geq2$, there are continuity intervals of $\T_\varepsilon^{n_0}$ asymptotic to $\varepsilon$, so there is no uniform Lasota-Yorke estimate for $\hat \perron_{\varepsilon}$.
Nevertheless, when considering the dynamics on appropriate subintervals, and studying this way families of maps with holes, it may still be possible to obtain uniform Lasota-Yorke inequalities for the corresponding transfer operators. For example, in the next section we prove the existence of a uniform Lasota-Yorke inequality for the map $T_{\varepsilon, 0}: I_L\cap T_\varepsilon^{-1}I_L \to I_L$, where $T_\varepsilon= U\circ \T_\varepsilon \circ U^{-1}$ is the conjugated map.
\begin{remark}\label{remark: C6 opposite}
    If instead of \ref{C6} $\T_\varepsilon$ satisfies
    \begin{align}\label{eq: opposite of c6}
        sign(\partial_-\T (a))\cdot \T_\varepsilon(a_\varepsilon)
        \leq
        sign(\partial_-\T (a))\cdot a_\varepsilon,
    \end{align}
    then $a_\varepsilon$ does not introduce shrinking continuity intervals for $T_\varepsilon^{n_0}$. Combining this with \ref{C2} and \ref{C4}, no other critical points give rise to such complications. As a result, if \eqref{eq: opposite of c6} holds instead of \ref{C6}, the family admits a uniform Lasota-Yorke inequality, hence it is ACIM stable.
\end{remark}

\subsection{Open systems}\label{section: open systems}

We will rely heavily on the following proposition throughout the proofs. The proposition is found in \cite[Proposition 3]{difcoef10}, where the authors consider more conventional maps. However, the requirements for the system can be relaxed as follows.

Let $T_\varepsilon$, $\varepsilon\geq0$ be a family of piecewise expanding maps such that there is a partition $I=I_0\cup I_1$ such that both subintervals $I_j$ are invariant under $T_0$, and support a unique absolutely continuous $T_0$-invariant measure, with density $\phi_j\in BV(I_j)$. Furthermore, there are holes $H_{j,\varepsilon}=I_j\cap T_\varepsilon^{-1}I_{1-j}$ with  $\lim_{\varepsilon\to 0}Leb(H_{j, \varepsilon})/\varepsilon = \beta_j\in \R^+$.
For $j \in \{0, 1\}$, let $\perron_{j,\varepsilon}$ be the transfer operator acting on $BV(I_j)$, where we consider $H_{j,\varepsilon}$ as a hole, i.e.,
\begin{align}\label{eq: open perron-frobenius definition}
\mathcal{L}_{j,\varepsilon} A(x) = \sum_{y \in (I_j \cap T_{\varepsilon}^{-1} \{x\})} \frac{A(y)}{|\D T_{\varepsilon}(y)|}.
\end{align}
Assume that the operators $\perron_{j,\varepsilon}$ admit a uniform Lasota-Yorke inequality \eqref{eq: Lasota-Yorke} and that $\perron_{j,\varepsilon} : [0, 1[ \to \mathcal{B}(BV, L^1) $ is continuous (as a function of $\varepsilon$) at $0$ for $j=0,1$.

 Then, as described in \cite[Proposition 3]{difcoef10}, the setting of \cite{KL99} applies: for small $\varepsilon > 0$, $\perron_{j,\varepsilon}$ has an isolated simple eigenvalue $\lambda_{j,\varepsilon} < 1$ of multiplicity one, with $\lambda_{j,\varepsilon} \to 1$ as $\varepsilon \to 0$, and otherwise, the spectrum has a uniform spectral gap. More precisely, if $Q_{j,\varepsilon}$ is the spectral projection corresponding to $\lambda_{j,\varepsilon}$, the following statement holds.
\begin{prop}\label{prop: open}
    There exists a measure $\nu_{j,\varepsilon}$ and a function $\phi_{j,\varepsilon} \in BV(I_j)$ such that
    \[
    \nu_{j,\varepsilon}(\phi_{j,\varepsilon}) = 1, \quad Q_{j,\varepsilon} = \nu_{j,\varepsilon}(\cdot)\phi_{j,\varepsilon}
    \]
    and
    \begin{itemize}
        \item[(a)] $\lambda_{j,\varepsilon} = 1 - \beta_j \varepsilon + o(\varepsilon).$
        \item[(b)] As $\varepsilon \to 0$, $\nu_{j,\varepsilon}(A) \to \int A(y) dy$ where the convergence holds in strong topology in $BV^*$, and $\phi_{j,\varepsilon} \to \phi_j$ in $L^1$. Moreover,
        \[
        \lim_{\varepsilon \to 0} \sup_{H_{j,\varepsilon}} |\phi_{j,\varepsilon} - \phi_j| = 0.
        \]
        \item[(c)] There exists $C > 0$ such that, for all small $\varepsilon$, $A \in BV(I_j)$, and all $n \geq 0$,
        \[
        \BV{ \lambda_{j,\varepsilon}^{-n} \mathcal{L}_{j,\varepsilon}^n A - Q_{j,\varepsilon} A } \leq C\theta^n \BV{A}.
        \]
    \end{itemize}
\end{prop}
Recall that our maps only experience anomalous behaviour in the holes, thus $T_\varepsilon$ converges to $T$ in the Skohorod metric (\cite[Section 11.2]{laws of chaos}), hence by \cite[Lemma 11.2.1.]{laws of chaos} $\| \perron_{j,\varepsilon}-\perron_j \|_{BV\to L^1}\to 0$ for $j=0,1$. Since $\| \perron_\varepsilon-(\perron_{0, \varepsilon} + \perron_{1, \varepsilon}) \|_{BV\to L^1}\to 0$, $\| \perron_\varepsilon-\perron \|_{BV\to L^1}\to 0$ also holds.

Also, $\|\perron_\varepsilon\|_{L^1}$ and $\|\perron \|_{BV}$ are both uniformly bounded, thus one can prove the following by induction
\begin{lemma}\label{lemma: continuity}
    For any $n\in\N$, $\perron^n :[0, 1[ \to \mathcal{B}(BV, L^1) $ is continuous at $0$.
\end{lemma}
This means that for any $f\in BV$, for any $n\in \N^+$,
$\int_I \perron_\varepsilon^n f = \int_I \perron^n f + O(\varepsilon)\BV{f}$.

By construction $T_{\varepsilon, 0}: I_0\cap T_\varepsilon^{-1}I_0 \to I_0$ always has $1/4$ as a fixed point and $T_{\varepsilon, 0}(x)=1/4$ for all $x\in \partial H_0$. Notice that the critical point $1/4$ and all critical points of $\partial H_0$ are our own creations through the conjugation. However, since these critical points map directly onto the critical fixed point $1/4$ for all $\varepsilon>0$, they do not give rise to shrinking continuity intervals to any fixed power $(T_{\varepsilon, 0})^{n_0}$. By \ref{C2} and \ref{C4} there are no other critical points that would prohibit us from taking an appropriate power of $T_{\varepsilon, 0}$ when constructing the uniform Lasota-Yorke inequality.

Hence the family $T_{\varepsilon, 0}: I_0\cap T_\varepsilon^{-1}I_0 \to I_0$ satisfies a uniform Lasota-Yorke inequality \eqref{eq: Lasota-Yorke}.

As a consequence, Proposition \ref{prop: open} holds for the conjugated family $T_{\varepsilon,0}$, $\varepsilon\geq 0$.

\section{The proof of Theorem \ref{theo: convex}}\label{section: proof of theo 1}

\subsection{The value of the $l.h.r.$}

To calculate the numerical value of $l.h.r.$, first recall that by \ref{C2}, $\hat\phi$ is continuous at every $x\in H$. Hence,
\begin{equation*}
    \hat{\mu}(\hat{H}_0)
    =
    \sum_{x \in H}|\D \T x|^{-1}\hat\phi(x)Leb(\B) + o(\varepsilon).
\end{equation*}

We can calculate $a_1$ based on the fact that it is a fixed point in the vicinity of $a_\varepsilon$, hence

\begin{align*}
     o(\varepsilon)
     =&
     \partial_-\T_\varepsilon(a_\varepsilon)(a_1-a_\varepsilon) - 1\cdot (a_1 - \T_\varepsilon(a_\varepsilon)).
    \\
    a_1
    =&
    \frac{\partial_-\T_\varepsilon(a_\varepsilon)a_\varepsilon - \T_\varepsilon(a_\varepsilon) } {\partial_-\T_\varepsilon(a_\varepsilon) - 1}
    +
    o(\varepsilon)
\end{align*}

Also,

\begin{align*}
     o(\varepsilon)
     =&
     \partial_+\T_\varepsilon(a_\varepsilon)(a_2-a_\varepsilon) - 1\cdot (a_1 - \T_\varepsilon(a_\varepsilon))
    \\
    a_2
    =&
    \frac{\partial_+\T_\varepsilon(a_\varepsilon)a_\varepsilon - \T_\varepsilon(a_\varepsilon) } {\partial_+\T_\varepsilon(a_\varepsilon)}
    +
    \frac{\partial_-\T_\varepsilon(a_\varepsilon)a_\varepsilon - \T_\varepsilon(a_\varepsilon) } {\partial_+\T_\varepsilon(a_\varepsilon)(\partial_-\T_\varepsilon(a_\varepsilon) - 1)}
    +
    o(\varepsilon)
\end{align*}

\begin{align*}
    Leb(\B)
    =&
    a_2 - a_1
    =
    \frac{\partial_+\T_\varepsilon(a_\varepsilon)a_\varepsilon - \T_\varepsilon(a_\varepsilon) } {\partial_+\T_\varepsilon(a_\varepsilon)}
    +
    (\partial_+\T_\varepsilon(a_\varepsilon)^{-1} - 1)
    \frac{\partial_-\T_\varepsilon(a_\varepsilon)a_\varepsilon - \T_\varepsilon(a_\varepsilon) } {\partial_-\T_\varepsilon(a_\varepsilon) - 1} + o(\varepsilon)
    \\
    =&
    \left(
    \partial_+\T_\varepsilon(a_\varepsilon)(\partial_-\T_\varepsilon(a_\varepsilon) - 1)
    \right)^{-1}
    \left(
    \partial_+\T_\varepsilon(a_\varepsilon)
    -
    \partial_-\T_\varepsilon(a_\varepsilon)
    \right)
    \left(
    \T_\varepsilon(a_\varepsilon)-a_\varepsilon
    \right)
    +
    o(\varepsilon)
    \\
    \sim&
    \left(
    \partial_+\T(a)(\partial_-\T(a) - 1)
    \right)^{-1}
    \left(
    \partial_+\T(a)
    -
    \partial_-\T(a)
    \right)
    c_a\varepsilon
    +
    o(\varepsilon)
    \\
    =&
    -c_a\left(
    (\partial_+\T(a) - \mathbb{1}_{\partial_+\T(a)  > 0})(\partial_-\T(a) - \mathbb{1}_{\partial_-\T(a)  > 0})
    \right)^{-1}
    \left(
    |\partial_+\T(a)|
    +
    |\partial_-\T(a)|
    \right)
    \varepsilon
    +
    o(\varepsilon)
\end{align*}
This expression is appropriate independent of the "orientation" of the turning point.

Also,
\begin{align*}
    m_\B(\hat{H}_1)=\frac{Leb(\hat H_1)}{Leb(\B)} = \frac{h_\varepsilon-1}{h_\varepsilon} + o(\varepsilon).
\end{align*}

We can see from the definition of $c_\pm$ \eqref{eq: cpm} that
\begin{align*}
    \partial_{\pm}\T_\varepsilon a_\varepsilon
    =&
    \partial_{\pm}\T a + c_\pm\varepsilon + o(\varepsilon).
\end{align*}

Using that for a $y>0$, $g(\varepsilon)=(y+c\varepsilon)^{-1} = y^{-1} - c y^{-2}\varepsilon + o(\varepsilon)$, we can see that

\begin{align*}
    \frac{h_\varepsilon-1}{h_\varepsilon}
    =&
    1-h_\varepsilon^{-1}
    =
    |c_-| \partial_{-}\T (a)^{-2}\varepsilon + | c_+ |\partial_{+}\T (a)^{-2}\varepsilon + o(\varepsilon).
\end{align*}

\begin{align*}
    l.h.r. = \frac{|c_-| \partial_{-}\T (a)^{-2} +  |c_+| \partial_{+}\T (a)^{-2}}
    {-c_a \left(
    (\partial_+\T(a) - \mathbb{1}_{\partial_+\T(a)  > 0})(\partial_-\T(a) - \mathbb{1}_{\partial_-\T(a)  > 0})
    \right)^{-1}
    \left(
    |\partial_+\T(a)|
    +
    |\partial_-\T(a)|
    \right)
     \sum_{x \in H}|\D \T x|^{-1}\hat\phi(x)}
\end{align*}

\subsection{Case (I)}

If $h_\varepsilon\leq 1$ for all $\varepsilon > 0$, then the box is an invariant interval of $\T_\varepsilon$.
Hence, we can decompose the Perron-Frobenius operator $\hat\perron_\varepsilon$ into two parts,

\begin{align*}
    \hat\perron_\varepsilon = \hat\perron_{0, \varepsilon} + \hat\perron_{1, \varepsilon},
\end{align*}
where
\begin{align*}
    \hat\perron_{0, \varepsilon}:& BV(I\setminus \B) \to BV(I)
    \\
    \hat\perron_{1, \varepsilon}:& BV(\B) \to BV(\B).
\end{align*}

First we show that for any fixed $\varepsilon>0$, and for all $f\in BV$, $\hat\perron_\varepsilon^{n} f$
converges to $\hat\phi_{1, \varepsilon}$, the ACIM of $\hat{T}|_{\mathcal{B}}$.
Fix a $\varepsilon>0$. We can see that by Proposition \ref{prop: open},
\begin{align*}
    \BV{(\hat\perron_\varepsilon^n f)|_{I\setminus\B}}
    =&
    \BV{\hat\perron_{0, \varepsilon}^n f}
    \leq
    \lambda_{0, \varepsilon}^n \BV{Q_{0, \varepsilon} f} + \BV{f}O(\theta^n)
\end{align*}

For the fixed $\varepsilon>0$, we can choose an $n_0$ such that
\begin{align*}
    \BV{(\hat\perron_\varepsilon^{n_0} f)|_{I\setminus\B}} < \varepsilon\BV{f}
\end{align*}
for all $f\in BV$.

$\hat\perron_{1, \varepsilon}$ is the transfer operator of a piecewise  expanding map on $\B$, hence it has ACIM $\hat\phi_{1, \varepsilon}\in BV$ Thus,

\begin{align*}
    \BV{\hat\perron_\varepsilon^{n+n_0} f - \hat\phi_{1, \varepsilon}}
    \leq&
    \BV{\hat\perron_\varepsilon^n (\hat\perron_\varepsilon^{n_0} f)|_\B
    -
    \hat\phi_{1, \varepsilon}}
    +
    \BV{\hat\perron_\varepsilon^n (\hat\perron_\varepsilon^{n_0}   f)|_{I\setminus\B}}
    \\
    \leq&
    \BV{\hat\perron_{1, \varepsilon}^n (\hat\perron_\varepsilon^{n_0} f)|_\B
    -
    \hat\phi_{1, \varepsilon}}
    +
    \varepsilon\BV{f}\BV{\hat\perron_\varepsilon^n}.
\end{align*}

Considering that for a fixed $\varepsilon>0$, $\perron_\varepsilon$ satisfies a Lasota-Yorke inequality, we can see that $\hat\perron_\varepsilon^{n} f$ converges to $\hat\phi_{1, \varepsilon}$ for for all $f\in BV$. Since $\hat\phi_{1, \varepsilon}$ is supported in $\B$,

\begin{align*}
    \hat\mu_\varepsilon(g)
    =&
    \int_{\B}g\hat\phi_{1, \varepsilon} \to g(a).
\end{align*}
holds for any $g: I \to \R$ observable, which is continuous at $a$.

Notice that the proof did not utilize the conjugated family $T_\varepsilon$, $\varepsilon > 0$.

We will continue with the parallel analysis of cases (II) and (III).

\subsection{Uniform Lasota-Yorke inequalities}\label{section: lasota yorke}

\subsubsection{Case (II) and (III)}\label{section: lasota yorke (III)}

We will show that the Perron-Frobenius operators $\perron_\varepsilon$ satisfy a uniform Lasota-Yorke inequality with respect to a specific norm. In order to reach this goal, we have to understand the behavior of the following operators:
\begin{align*}
    \perron_{0, \varepsilon}:& BV(I_L \setminus H_0) \to BV(I_L), \quad
    &&\perron_{0, \varepsilon}f(x) = \sum_{y\in T_{\varepsilon}^{-1}x \cap I_L} \frac{1}{\D T_{\varepsilon}y}f(y)
    \\
    \perron_{1, \varepsilon}:& BV(I_R \setminus H_1) \to BV(I_R), \quad
    &&\perron_{1, \varepsilon}f(x) = \sum_{y\in T_{\varepsilon}^{-1}x \cap I_R} \frac{1}{\D T_{\varepsilon}y}f(y)
    \\
    \perron_{0\to 1, \varepsilon}:& BV(H_0) \to BV(T_\varepsilon(H_0)), \quad
    &&\perron_{0\to 1, \varepsilon}f(x) = \sum_{y\in T_{\varepsilon}^{-1}x \cap I_L} \frac{1}{\D T_{\varepsilon}y}f(y)
    \\
    \perron_{1\to 0, \varepsilon}:& BV(H_1) \to BV(T_\varepsilon(H_1)), \quad
    &&\perron_{1\to 0 \varepsilon}f(x) = \sum_{y\in T_{\varepsilon}^{-1}x \cap I_R} \frac{1}{\D T_{\varepsilon}y}f(y)
\end{align*}

As discussed in Secion \ref{section: open systems}, the open transfer operators $\perron_{0, \varepsilon}$ satisfy a uniform Lasota-Yorke inequality

\begin{align*}
    \V{\perron_{0, \varepsilon}^nf} \leq C_1(\alpha_0^n\V{f} + \Lnorm{f})
\end{align*}
for some $C_1>0$ and $\alpha_0<1$.

Regarding the other operators, recall the discussion in Section \ref{section: lasota-yorke methodology}.

The map of the dynamics of $T_\varepsilon$ from $I_R\cap T_\varepsilon ^{-1} I_R$ to $I_R$ converges in $C^2$ (as described in Condition \ref{C5}) to the dynamics of $T$ from $I_R \cap T^{-1} I_R$ to $I_R$,
which is a tent map. Since both branches of these maps are surjective, we can see from \ref{section: lasota-yorke methodology} that for $\varepsilon<\varepsilon_0$, all $\perron_{1, \varepsilon}$ satisfy a uniform Lasota-Yorke inequality,

\begin{align*}
    \V{\perron_{1, \varepsilon}^nf} \leq C_2(\alpha_1^n\V{f} + \varepsilon\Lnorm{f})
\end{align*}
for some $C_2>0$ and $\alpha_1<1$. In particular, $\alpha_1^{-1}=\inf_{\varepsilon<\varepsilon_0}\{\partial_-T_\varepsilon(a), \partial_+T_\varepsilon(a)\}$.

Regarding $\perron_{0 \to 1, \varepsilon}$, we can see that the the derivatives are asymptotic to $\varepsilon^{-1}$, so by the definition of the operator \eqref{eq: open perron-frobenius definition},

\begin{align*}
    \V{\perron_{0 \to 1, \varepsilon}f} \leq C_3\varepsilon\BV{f}
\end{align*}

for some $C_3>0$.

The last family of operators is $\perron_{1 \to 0, \varepsilon}$.
Notice that the map $T_\varepsilon|_{H_1}: H_1 \to T_\varepsilon(H_1)$ only has surjective branches.
The derivatives are asymptotic to $\varepsilon$, and
\begin{align*}
    \frac{\|\D^2T_\varepsilon|_{H_1}\|}{\min|\D T_\varepsilon|^2}
    =
    \frac{\|\D^2\T_\varepsilon|_{\hat H_1}\|}{\min|\D \T_\varepsilon|^2},
\end{align*}
hence
\begin{align*}
    \V{\perron_{1 \to 0, \varepsilon} f} \leq C_4(\varepsilon^{-1}\V{f} + \Lnorm{f}).
\end{align*}

We can assume that the above inequalities hold for a common $C>1$ and $\alpha<1$.
The Lasota-Yorke inequality for $\perron_\varepsilon$ in this case will be derived in the next chapter.

\subsection{The bounded $\varepsilon$-variance}\label{section: e-variance}

Notice that if $\hat f\in BV(I)$ an observable for the original system $\hat{T}: I\to I$, then the conjugate observable $f=\hat f \circ U^{-1}$ is such that $\V{f|_{I_R}}\leq \varepsilon \V{\hat{f}}$.

We introduce a new notion, the \textit{$\varepsilon$-variance} of $f\in BV(I)$,
\begin{align*}
    \Va{f} = \V{f|_{I_L}} + \varepsilon^{-1}\V{f|_{I_R}} + \BV{f|_{I_L}} + \BV{f|_{I_R}} \geq \V{f}.
\end{align*}
With this, we define the $BV_\varepsilon$ norm,
\begin{align*}
    \BVa{f} = \Va{f} + \Lnorm{f}.
\end{align*}

We will show that in case (III), for a large enough $\hat K>0$, the set
\begin{align*}
    \mathcal{F}_\varepsilon = \{f\in BV(I) \; :\; \Va{f}\leq \hat K, \; \Lnorm{f}\leq 1 \}
\end{align*}
is eventually invariant.

Note that

\begin{align*}
        \left\{ (\hat f + \varepsilon^{-1}t 1_\B) \circ U^{-1} \; \bigg| \;  t\in \R,\, \hat f\in BV(I), \, \V{\hat f} \leq \hat K, \, \Lnorm{\hat f + \varepsilon^{-1}t 1_\B}\leq 1
    \right\}
    \subset
    \mathcal{F}_\varepsilon,
\end{align*}
(where $t\in \R$ is independent from $\varepsilon$.)
Hence, in particular, the conjugates of all functions in a ball of $BV$ are included.

We need a short analysis of the behavior of the map $\T_\varepsilon$ near $a_\varepsilon$, see Figure~\ref{fig: holes}, right panel. A key observation is that for small $\varepsilon$, the images $\T^k\hat{H}_1$ remain in $I\setminus \B$ for long intervals of time $k\ge 1$. Indeed, this hole (from the box to the exterior) maps to the right side of the box $a_2$, then to the left side, where the right endpoint of $\T^2 \hat{H}_1$, $a_1$, is a fixed point of $\hat T$. Hence, the following subsequent images of the hole have the same right endpoint, and are becoming longer.
As a consequence, $T_\varepsilon$ admits similar properties, see Figure~\ref{fig: conjugated holes}, right panel for reference.
The second image $T^2H_1$ has a fixed endpoint, and thus its subsequent images remain connected and on the same branch of $T$ for a prolonged time.
In particular, let
\begin{align*}
    I_0 =& T H_1
    \\
    I_k =& T^{k}I_0 \setminus T^{k-1}I_0, k\geq 1
\end{align*}
in fact, we will only consider $I_k$ for $k\leq N_\varepsilon$, which is defined as follows.

We will use the notation $H_\varepsilon = H_0 \cup H_1$.
Since $Leb(H_1)\asymp \varepsilon$, and $|\D T_\varepsilon|_{H_1}|\asymp \varepsilon$, $Leb(I_0)\asymp \varepsilon^2$. Furthermore, $L_1 = \lim_{\varepsilon\to 0} Leb(I_0)/\varepsilon^2$ exists.
Set a $\gamma<1$ and define
\begin{align}\label{eq: Na}
    N_\varepsilon
    =
    \min\{k\geq 1\; : \; Leb(T^{k+1} I_0)\geq\varepsilon^{\gamma}L_1\}
    \in
    \left[ \frac{(\gamma-2)\log\varepsilon}{\log \sup_{I \setminus H_\varepsilon}|DT|} - 1,
    \frac{(\gamma-2)\log\varepsilon}{\log \inf_{I\setminus H_\varepsilon}|DT|}\right].
\end{align}
Then
\begin{align}\label{eq: expansion on irregular part}
    |D T^{N_\varepsilon}(x)|\asymp \varepsilon^{\gamma-2}
\end{align}
uniformly on $I_0$ and
\begin{align}
    |\cup_{k=0}^{N_\varepsilon} I_k|
    \ll
    \varepsilon^{\gamma}.
\end{align}
Notice that $I_k$ is on the same branch $B_{I_1}$ of $T_\varepsilon$ as $I_1$ for $1<k\leq N_\varepsilon$.
The exact value of $\gamma$ will only be important in the proof of Lemma \ref{prop: growth lemma}. Notice that $I_{k+1}=T_\varepsilon I_k$ for $k\neq1$ but $T_\varepsilon I_1 = I_1 \cup I_2$. For this reason, we need one more open Perron-Frobenius operator,
\begin{align*}
    \perron_{I_1, \varepsilon}:& BV(I_1 \cap T_\varepsilon^{-1} I_1) \to BV(I_1), \quad
    &&\perron_{I_1, \varepsilon}f(x) = \sum_{y\in T_{\varepsilon}^{-1}x \cap I_1} \frac{1}{\D T_{\varepsilon}y}f(y).
\end{align*}
For an $f\in BV(I)$, let us define the following decomposition.
\begin{align*}
    f_L &= f|_{I_L}
    \\
    f_R &= f|_{I_R}
    \\
    0 = f_k &: I_k \to \R , \quad 0\leq k \leq N_\varepsilon
\end{align*}
Using this, we can see the decomposition of the evolved density,
\begin{align*}
    (\perron_\varepsilon f)_L
    &=
    \perron_{0, \varepsilon} f_L + \perron_\varepsilon f_{N_\varepsilon}
    \\
    (\perron_\varepsilon f)_R &= (\perron_\varepsilon f)|_{ I_R} = \perron_{1, \varepsilon} f_R + \perron_{0 \to 1, \varepsilon}f_L
    \\
    (\perron_\varepsilon f)_0 &= \perron_{1 \to 0,  \varepsilon}f_R
    \\
    (\perron_\varepsilon f)_1 &= \perron_{0,  \varepsilon}f_0 + \perron_{I_1,  \varepsilon}f_1
    \\
    (\perron_\varepsilon f)_k &= \perron_{0,  \varepsilon}(f_{k-1}), \quad 2 \leq k\leq N_\varepsilon
\end{align*}

This partition is useful because the evolved density $\perron_\varepsilon^n f$ shows anomalous behavior on $I_k$, since the derivative of $T$ in $H_1$ converges to zero. By separating the $f_k$-s from the two regular parts, $f_L$ and $f_R$, we "freeze" the evolution process, and only let them back when they have already regularized.

In Proposition\ref{prop: Lasota-Yorke (III)}, we will use the above defined $\varepsilon$-variance and the set $\mathcal{F}_\varepsilon$.

\begin{prop}\label{prop: Lasota-Yorke (III)}
     In cases (II) and (III), the maps $T_\varepsilon: I \to I$ satisfy a uniform Lasota-Yorke inequality, with respect to the $\varepsilon$-variance, for any $f\in BV_\varepsilon(I)$. That is, there are constants, $\beta<1$, $K\geq 1$ such that
    \begin{align*}
        \Va{\perron_\varepsilon^{n} f}
        \leq&
        K(\beta^n\Va{f} + \Lnorm{f})
        \\
        \Lnorm{\perron_\varepsilon^{n} f}
        \leq&
        K \Lnorm{f}
    \end{align*}
    holds for all $f\in BV_\varepsilon(I) $ and $n\in \N^+$.
\end{prop}

This will imply that $\mathcal{F}_\varepsilon$ is eventually invariant for $\perron_\varepsilon$.

Note that in Section \ref{section: lasota-yorke for a family} $BV(I)$ is the same space for all $\varepsilon>0$.
However, we have a family of spaces $(BV_\varepsilon(I), \BVa{\cdot}) \subset (L^1(I), \Lnorm{\cdot})$.
By Proposition \ref{prop: Lasota-Yorke (III)}, there is a Lasota-Yorke inequality for all $\varepsilon$.
This implies that for every $\varepsilon>0$, the operator $\perron_\varepsilon$ has an invariant density $\f_\varepsilon\in BV_\varepsilon$.
Moreover, since the constants in Proposition \ref{prop: Lasota-Yorke (III)} are uniform in $\varepsilon$, there is a $B>0$ (independent of $\varepsilon$) such that $\BVa{\f_\varepsilon} \leq B$.

\begin{customproof}{Proof of Proposition \ref{prop: Lasota-Yorke (III)}}
Recall from subsection \ref{section: lasota yorke (III)} that the various transfer operators satisfy the following inequalities.
\begin{align*}
    \V{\perron_{0, \varepsilon}^n f} \leq C(\alpha^n\V{f} + \Lnorm{f})
\end{align*}

\begin{align*}
    \V{\perron_{1, \varepsilon}^n f} \leq C(\alpha^n\V{f} + \varepsilon\Lnorm{f})
\end{align*}

\begin{align*}
    \V{\perron_{0\to 1, \varepsilon} f} \leq \varepsilon C(\V{f} + \Lnorm{f})
\end{align*}

\begin{align*}
    \V{\perron_{1 \to 0, \varepsilon} f} \leq C(\varepsilon^{-1}\V{f} + \Lnorm{f})
\end{align*}

\begin{align*}
    \V{\perron_{I_1, \varepsilon}^n f} \leq C(\alpha^n\V{f} + \Lnorm{f}).
\end{align*}

Since all of these estimates are adequate for both cases (II) and (III), we can prove the estimates for them simultaneously.
We will utilize the following estimate in the calculations. If we have a subinterval $J\subset I$, $f\in BV(J)$ such that $\hat f \in BV(I)$ is an extension of $f$ that is zero outside $J$, then we have

\begin{align}\label{eq: variance of an extension}
    \V{\hat{f}} \leq \V{f} + 2 \BV{f} \ll \BV{f}.
\end{align}

We will show that there exists some $n_0\ge 1$, and constants $r<1$ and $K>0$  such that, for any  $\varepsilon$ small enough, and any $f\in BV_\varepsilon(I)$, an estimate of the type \eqref{eq: one-step lasota} holds:
\begin{align*}
\Va{\perron_\varepsilon^{n_0} f}\le r \cdot \Va{f} + K\cdot \Lnorm{f}.
\end{align*}
This automatically implies $\perron^{n_0}f\in BV_\varepsilon(I)$, and thus Proposition~\ref{prop: Lasota-Yorke (III)} follows by standard iteration.

The left regular part of $f$ will evolve by the open transfer operator acting on $f_L$ and by $f_R$ evolving through the intervals $I_k$ and flowing into $(\perron^{n_0} f)_L$. However, since $n_0$ is fixed, we can assume (by taking small enough $\varepsilon$) that  $n_0\leq N_\varepsilon$. Thus within this $n_0$ steps $f_R$ can not evolve through the irregular intervals, hence there is no contribution form the right side. This means, that only the left open transfer operator evolves the left part of $f$. In the computation below, we assume $n\le n_0\le N_{\varepsilon}$

\begin{align*}
    \V{(\perron_\varepsilon^{n} f)_L}
    =&
    \V{\perron_{0, \varepsilon}^{n} f_L}
    \leq
    C(\alpha^n\V{f_L} + \Lnorm{f})
\end{align*}

\begin{align*}
    \V{(\perron_\varepsilon^n f)_R}
    =&
    \V{\perron_{1, \varepsilon}^{n} f_R + \sum_{k=1}^{n-1} \perron_{1, \varepsilon}^{k} \perron_{0 \to 1, \varepsilon}\perron^{n-k-1}_{0, \varepsilon}f_L}
    \ll
    \BV{\perron_{1, \varepsilon}^{n} f_R} + \sum_{k=1}^{n-1} \BV{\perron_{1, \varepsilon}^{k} \perron_{0 \to 1, \varepsilon}\perron^{n-k-1}_{0, \varepsilon}f_L}
    \\
    \leq&
    C\left(\alpha^n\V{f_R} + \varepsilon \Lnorm{f_R}\right)
    +
    \sum_{k=1}^{n-1}
    C\left(\alpha^k\V{\perron_{0 \to 1, \varepsilon}\perron^{n-k-1}_{1, \varepsilon}f_L} + \varepsilon \Lnorm{\perron_{0 \to 1, \varepsilon}\perron^{n-k-1}_{0, \varepsilon}f_L}\right)
    + \Lnorm{f}
    \\
    \leq&
    C\alpha^n\V{f_R} + \varepsilon C\Lnorm{f}
    +
    \sum_{k=1}^{n-1}
    C^2\alpha^k\varepsilon(\V{\perron^{n-k-1}_{0, \varepsilon}f_L} + \Lnorm{\perron^{n-k-1}_{0, \varepsilon}f_L}) + \varepsilon C\Lnorm{\perron_{0 \to 1, \varepsilon}\perron^{n-k-1}_{0, \varepsilon}f_L}
    + \Lnorm{f}
    \\
    \leq&
    C\alpha^n\V{f_R}
    +
    \sum_{k=1}^{n-1}
    C^2\alpha^k\varepsilon (\V{\perron^{n-k-1}_{0, \varepsilon }f_L} + \Lnorm{f}) + \varepsilon n3C\Lnorm{f}
    \\
    \leq&
    C\alpha^n\V{f_R}
    +
    \sum_{k=1}^{n-1}
    C^3\alpha^k\varepsilon \alpha^{n-k-1}\V{f_L} + \varepsilon n4C^2\Lnorm{f}
    \\
    \leq&
    C\alpha^n\V{f_R}
    +
    \varepsilon n\alpha^{n-1}C^3\V{f_L}
    +
    \varepsilon n4C^2\Lnorm{f}
    \\
    \ll&
    \alpha^n\V{f_R}
    +
    \varepsilon n\alpha^{n-1}\V{f_L}
    +
    \varepsilon n\Lnorm{f}
\end{align*}

\begin{align*}
    \V{(\perron_\varepsilon^n f)_0}
    \leq&
    \V{\perron_{1 \to 0, \varepsilon }(\perron_\varepsilon^{n-1} f)_R}
    \\
    \leq&
    C\varepsilon^{-1}\V{(\perron_\varepsilon^{n-1} f)_R}
    +
    C\Lnorm{(\perron_\varepsilon^{n-1}f)_R}
    \\
    \ll&
    \varepsilon^{-1}
    \alpha^{n-1}\V{f_R} +  (n-1)\alpha^{n-2} \V{f_L}
    +
    (n-1)\Lnorm{f}
\end{align*}

\begin{align*}
    \V{(\perron_\varepsilon^n f)_1}
    =&
    \V{
    \perron_{0, \varepsilon }(\perron_\varepsilon^{n-1} f)_0
    +
    \perron_{I_1, \varepsilon }(\perron_\varepsilon^{n-1} f)_1
    }
    \\
    =&
    \V{\sum_{k=1}^{n-1} \perron^{k-1}_{I_1, \varepsilon }\perron_{0, \varepsilon}(\perron_\varepsilon^{n-k}f)_0}
    \\
    \ll&
    \varepsilon^{-1}
    (n-1)\alpha^{n-1}\V{f_R} +  (n-1)^2\alpha^{n-2} \V{f_L}
    +
    (n-1)^2\Lnorm{f}
\end{align*}

For $k>1$,

\begin{align*}
    \V{(\perron_\varepsilon^{n} f)_k}
    \leq&
    \V{\perron^{k-1}_{0, \varepsilon}(\perron_\varepsilon^{n-k+1} f)_1}
     \\
     \ll&
     \varepsilon^{-1}
    (n-1)\alpha^{n-1}\V{f_R}
    +
    (n-1)^2\alpha^{n-2} \V{f_L}
    +
    (n-1)^2\Lnorm{f}
\end{align*}

Note that this inequality holds for $k=0, 1$ as well. Hence,

\begin{align*}
    \sum_{k=0}^{n-1}\V{(\perron_\varepsilon^n f)_k}
    \ll&
    n
    \left(
    \varepsilon^{-1}
    (n-1)\alpha^{n-1}\V{f_R}
    +
    (n-1)^2\alpha^{n-2} \V{f_L}
    +
    (n-1)^2\Lnorm{f}
    \right)
    \\
    \leq&
    \varepsilon^{-1}
    n^2\alpha^{n-1}\V{f_R}
    +
    n^3\alpha^{n-2} \V{f_L}
    +
    n^3\Lnorm{f}
\end{align*}

Combining all of this, we get

\begin{align*}
    \Va{\perron_\varepsilon^n f}
    =&
    \V{(\perron_\varepsilon^n f)_L + \sum_{k=0}^{n-1}(\perron_\varepsilon^n f)_k}
    +
    \varepsilon^{-1}\V{(\perron_\varepsilon^n f)_R}
    +
    \BV{(\perron_\varepsilon^n f)_L + \sum_{k=0}^{n-1}(\perron_\varepsilon^n f)_k}
    +
    \BV{(\perron_\varepsilon^n f)_R}
    \\
    \ll&
    \BV{(\perron_\varepsilon^n f)_L}
    +
    \BVa{(\perron_\varepsilon^n f)_R}
    +
    \sum_{k=0}^{n-1}\BV{(\perron_\varepsilon^n f)_k}
    \\
    \ll&
    \V{(\perron_\varepsilon^n f)_L}
    +
    \varepsilon^{-1}\V{(\perron_\varepsilon^n f)_R}
    +
    \sum_{k=0}^{n-1}\V{(\perron_\varepsilon^n f)_k}
    +
    \Lnorm{f}
    \\
    \ll&
    n^3\alpha^{n-2}\V{f_L}
    +
    n^2\alpha^{n-1}\varepsilon^{-1}\V{f_R}
    +
    n^3
    \Lnorm{f}
    \\
    \leq&
    n^3\alpha^{n-2}\Va{f}
    +
    n^3
    \Lnorm{f}
\end{align*}

Thus, for a large enough $n_0$, we have

\begin{align*}
    \Va{\perron_\varepsilon ^{n_0}f}\leq \alpha\Va{f} + \tilde{C}\Lnorm{f}
\end{align*}

for some $\tilde C>0$.
This implies that for an $f\in BV_\varepsilon$, $\perron_\varepsilon^{n_0}f\in BV_\varepsilon$, hence we can assume $(\perron_\varepsilon^{n_0} f)_L = (\perron_\varepsilon^{n_0} f)|_{I_L}$ and $(\perron_\varepsilon^{n_0} f)_k = 0$ for all $k\leq N_\varepsilon$.
Thus,

\begin{align*}
    \Va{\perron_\varepsilon^{n}f}
    \leq&
    K(\beta^n\Va{f} + \Lnorm{f})
\end{align*}
holds for some $\beta<1$ and $K\geq 0$.
\end{customproof}

\subsection{The invariant density and the $l.h.r.$}\label{section: invariant density}

We established in the previous sections that $\perron_\varepsilon$ has a fixed density $\f_\varepsilon \in BV_\varepsilon(I)$ in cases (II) and (III). For ease of notation, we will drop the $\varepsilon$ index, and write $\f$ going forward. We want to show that
\begin{align}\label{eq: convergence to the l.h.r.}
    \frac{\int_{I_L}\f}{\int_{I_R}\f}\to l.h.r.
\end{align}
We can decompose $\f$ as follows:

\begin{align*}
    \f_R =& \f|_{I_R},
    \\
    \f_0 =& \perron_{1 \to 0, \varepsilon}\f|_{I_R},
    \\
    \f_1 =& \sum_{n=1}^{\infty}\perron_{I_1, \varepsilon}^{n-1}\perron_{0,  \varepsilon}\f_0,
    \\
    \f_{k+1} =& \perron_{0, \varepsilon}\f_k, \quad 1\leq k < N_\varepsilon,
    \\
    \f_L =& \f|_{I_L}-\sum_{k=0}^{N_\varepsilon}\f_k.
\end{align*}

It has to satisfy
\begin{align*}
    \f_R
    =
    (\perron_\varepsilon \f)_R
    =
    \perron_{1, \varepsilon}\f_R + \perron_{0\to 1, \varepsilon}\f_L
\end{align*}

We will first prove \eqref{eq: convergence to the l.h.r.} in case (II).
Recall that
\begin{align*}
    \BV{\perron_{0\to 1, \varepsilon}\f_L}
    \ll
    \varepsilon\BV{\f_L} + \Lnorm{\f_L|_{H_1}}
    \ll
    \varepsilon\BV{\f_L}.
\end{align*}
Let $\lambda_\varepsilon=\BV{\perron_{1, \varepsilon}}$ and $\lambda=\BV{\perron_{1}}$.
Notice that in case (II), $\lambda_\varepsilon\leq\lambda<1$, thus
\begin{align*}
    \BV{\perron_{0\to 1, \varepsilon}\f_L}
    =
    \BV{\f_R-(\perron_{1, \varepsilon}\f_R)}
    \geq
    \BV{(1-\perron_{1, \varepsilon})\f_R}
    \geq
    (1-\lambda)\BV{\f_R}.
\end{align*}
Combining the last two estimates, we can see that $\BV{\f_R}=O(\varepsilon) $ in case (II). Since $\int_{I}\f=1$, this proves \eqref{eq: convergence to the l.h.r.} in this case.

Now we turn to case (III). We will show that $\f_L$ is close in $BV$ to the a.c.c.i.m. $\phi_{0,\varepsilon}$  of $\perron_{0, \varepsilon}$. Since $\phi_{0, \varepsilon}$ converges in $BV$ to the ACIM $\phi_0$ of $\perron_0$, this will imply $\int_{H_0}\f \sim \Lnorm{\f_L}\int_{H_0}\phi_{0, \varepsilon}\sim \Lnorm{\f_L}\mu_0(H_0)$. Since $\f\in BV_\varepsilon$, $\int_{H_1}\f \sim \Lnorm{\f|_{I_R}}m_{I_R}(H_1)$. We will complete the proof by showing that $\int_{H_0}\f = \int_{H_1}\f$.

Regarding $\f_L$,
\begin{align*}
    \f_L
    =
    (\perron_\varepsilon \f)_L
    =
    \perron_{0, \varepsilon}\f_L + \perron_{0, \varepsilon}\f_{I_{N_\varepsilon}}
\end{align*}

Thus,
\begin{align*}
    \BV{\f_L - \perron_{0, \varepsilon}\f_L }
    =&
    \BV{\perron_{0, \varepsilon}\f_{N_\varepsilon}}
    \leq
    \BV{(\perron_{0, \varepsilon}^{N_\varepsilon+1}\perron_{1 \to 0, \varepsilon}\f_R)}
    \\
    \ll&
    \alpha^{N_\varepsilon+1}\frac{\V{\f_R}}{\varepsilon} + \Lnorm{\f|_{H_1}}
    \\
    \leq&
    \alpha^{N_\varepsilon+1}\Va{\f_R} + \BV{\f_R}Leb(H_1) \to 0.
\end{align*}
Note that the last expression goes to zero as $\varepsilon\to 0$ because $Leb(H_1)\to 0$.

Now, let  $g_\varepsilon\in BV(I_L)$ be belong to the eigenspace of $\perron_{0, \varepsilon}$ with the largest eigenvalue $\lambda_\varepsilon<1$ (that is, a multiple of the a.c.c.i.m) such that $\Lnorm{g_\varepsilon}=\Lnorm{\f_L}$. Then there is some $h_\varepsilon\in BV(I_L)$ such that
\begin{align*}
    \f_L =& g_\varepsilon + h_\varepsilon
    \\
    \perron_{0,\varepsilon}\f_L =& \lambda_\varepsilon g_\varepsilon + \perron_{0, \varepsilon}h_\varepsilon
\end{align*}
holds. Since $h_\varepsilon$ is in the complement of the eigenspace spanned by $g_\varepsilon$,
\begin{align*}
    \BV{\f_L - g_\varepsilon }
    =
    \BV{h_\varepsilon}
    \ll&
    \BV{(\lambda_\varepsilon -\perron_{0, \varepsilon})h_\varepsilon}
    =
    \BV{ (\lambda_\varepsilon -\perron_{0, \varepsilon})\f_L}
    \\
    \leq&
    \BV{ (1 -\perron_{0, \varepsilon})\f_L}
    +
    (1-\lambda_\varepsilon)\BV{\f_L}
    \\
    =&
    \BV{ \f_L - \perron_{0, \varepsilon} \f_L}
    +
    (1-\lambda_\varepsilon)\BV{\f_L}=V_L.
\end{align*}
Since $1-\lambda_\varepsilon \asymp Leb(H_0)$ by Propositon \ref{prop: open}, this expression tends to zero.

We will utilize the following identity several times. For positive functions $f,g\in BV(I)$ and subset $A \subset I$ with $\int_A g >0$,
\begin{align*}
    \left|1-\frac{\int_A f}{\int_A g} \right|
    =
    \left|\frac{\int_A g-f}{\int_A g}\right|
    \leq
    \frac{\int_A |g-f|}{\int_A g}
    \leq
    \frac{Leb(A) \BV{g-f}}{Leb(A)\inf_A g}
    =
    \frac{\BV{g-f}}{\inf_A g}
\end{align*}
which implies
\begin{align}
    1 - \frac{\BV{g-f}}{\inf_A g}
    \leq
    \frac{\int_A f}{\int_A g}
    \leq
    1 + \frac{\BV{g-f}}{\inf_A g}.
\end{align}

Recall that $\mu_i$ is the ACIM with density $\phi_i=\perron_i\phi_i$.
Denote by $H_{i, +}$ the components of $H_i$ on which $\phi_i$ is positive for small enough $\varepsilon$.
By the uniform Lasota-Yorke inequality satisfied by $\perron_{0, \varepsilon}$, $\varepsilon \geq 0$, $\BV{\phi_{i, \varepsilon} - \phi_i} =L_i \to 0$. We denote $C_i= (\inf_{H_{i, +}}\phi_{i, \varepsilon})^{-1} \to (\inf_{\lim H_{i, +}}\phi_i)^{-1}$. Then,
\begin{align*}
    \int_{H_0}\f
    =&
    \int_{H_{0, +}}\f + o(\varepsilon)
    =
    \int_{H_{0, +}}\f_L + o(\varepsilon)
    \leq
    \int_{H_{0, +}}g_\varepsilon(1 + C_0 V_L)  + o(\varepsilon)
    \\
    =&
    \Lnorm{\f_L}\mu_0(H_{0, +})(1+C_0 V_L)(1+ C_0 L_0) + o(\varepsilon)
    =
    \Lnorm{\f_L}\mu_0(H_{0})(1+C_0 V_L)(1+ C_0 L_0) + o(\varepsilon)
\end{align*}
Doing the analogous lower estimates we can see that
\begin{align*}
    \Lnorm{\f_L}\mu_0(H_0)(1- C_0 V_L)(1- C_0 L_0 )
    \leq
    \int_{H_0}\f
    \leq
    \Lnorm{\f_L}\mu_0(H_0)(1+ C_0 V_L)(1+ C_0 L_0 )
\end{align*}

Since $H_1 = H_{1, +}$,
\begin{align*}
    \int_{H_1}\f
    =
    \int_{H_1}\f_R
    \leq
    \Lnorm{\f_R}\int_{H_1}2(1 + C_1 \V{\f_R})
    \leq
    \Lnorm{\f_R}\mu_1(H_1)(1 + C_1\varepsilon\Va{\f}).
\end{align*}

From this, it is clear that

\begin{align*}
    \Lnorm{\f_R}\mu_1(H_1)(1 -  C_1 \varepsilon\Va{\f})
    \leq
    \int_{H_1}\f
    \leq
    \Lnorm{\f_R}\mu_1(H_1)(1 +  C_1 \varepsilon\Va{\f})
\end{align*}

Notice that these integrals are equal, because
\begin{align*}
    0=&\Lnorm{\f_R}
    -
    \Lnorm{(\perron_\varepsilon \f)_R}
    =
    \Lnorm{\perron_{0\to 1, \varepsilon} \f_L} - \Lnorm{\perron_{1 \to 0, \varepsilon}\f_R}
    =
    \Lnorm{\f_L|_{H_0}} - \Lnorm{f_R|_{H_1}}.
\end{align*}

Combining these, we can see that

\begin{align*}
    \frac{\Lnorm{\f_L}}{\Lnorm{\f_R}}\to \lim_{\varepsilon\to 0}\frac{\mu_1(H_1)}{\mu_0(H_0)}
\end{align*}
If we further consider
\begin{align*}
        \Lnorm{\f|_{I_L}}
        =&
        \Lnorm{\f_L} + \sum_{i=0}^{N_\varepsilon} \Lnorm{\f_i}
        \leq
        \Lnorm{\f_L} + C N_\varepsilon \Lnorm{\f_R|_{H_1}}
        \\
        \leq&
        \Lnorm{\f_L} + C N_\varepsilon \Lnorm{\f_R}\mu_1(H_1)(1 + C_1 \varepsilon \Va{\f_R})u
        \leq
        \Lnorm{\f_L} +  \Lnorm{\f_R}O(\varepsilon\log\varepsilon),
    \end{align*}
then we can see that
\begin{align*}
    \frac{\Lnorm{\f|_{I_L}}}{\Lnorm{\f|_{I_R}}}\to \lim_{\varepsilon\to 0}\frac{\mu_1(H_1)}{\mu_0(H_0)}
    =
    \lim_{\varepsilon\to 0}\frac{m_\B(\hat H_1)}{\hat\mu(\hat H_0)}
    =
    l.h.r.
\end{align*}

\section{Proof of Theorem \ref{theo: markov}}\label{section: proof of theo 2}

\subsection{A Growth lemma}\label{section: growth lemma}

To prove Theorem \ref{theo: markov}, we need a growth Lemma, that is, give some curve $J\subset I$, need to control the lengths of the components of $T\varepsilon^n J$.

\begin{remark}
Since the derivative of $T_{\varepsilon}$ on $H_1$ tends to $0$ as $\varepsilon\to 0$ -- as reflected by the fact that for any density $f\in BV(H_1)$, $\perron_\varepsilon f$ becomes asymptotic to $\varepsilon^{-1}$ -- there is no hope that a general curve will evolve regularly. What we know, however, is that after passing through $H_1$, every curve stays connected and gets expanded by a factor at least $\inf_{I\setminus H_\varepsilon}|\D T_\varepsilon|$ for $N_\varepsilon$ steps.
\end{remark}

This motivates the following definition.
\begin{definition} Let
\begin{align*}
    B_\varepsilon = \cup_{k=1}^{N_\varepsilon}T^{k} H_1 = \cup_{k=0}^{N_\varepsilon-1}I_k.
\end{align*}
Given a curve $J\subset I$ and a time $n\in \N^+$ we only want to look at parts of $J_n$ that are not in $B_\varepsilon$, thus we define
\begin{align*}
    \overline J = J \setminus B_\varepsilon.
\end{align*}
\end{definition}

Recall that for a curve $J$ we will use $Leb_J$ for the Lebesgue measure restricted to $J$. This arclength measure induces a metric $d_J$ on $J$, and a conditional measure $m_J$ on $J$, meaning
\begin{align*}
    m_J(A)=\frac{Leb_J(A\cap J)}{Leb(J)}
\end{align*}
for any set $A$.
We will use the notation $J_n=T_\varepsilon^n J$, which is regarded as a collection of curves (which may overlap). For an evolved curve $J_n = \{J_n^{(i)}\}_{i\in I}$ $\overline{J_n}= \left\{\overline{J_n^{(i)}}\right\}_{i\in I}$.
For an element $x\in J$, let
\begin{align*}
    r^{(J)}_n(x)=d_{J_n}(T_\varepsilon^n x, \partial J_n)
\end{align*}
denote the distance of $T_\varepsilon^n x$ from the closest boundary point of $J_n$.
In the following lemma, one will encounter expressions akin to
\begin{align*}
    m_J(x\in \overline J: r_n(x)<\varepsilon') = m_J(x\in \overline J: r_n^{(\overline J)}(x)<\varepsilon').
\end{align*}
That is, we may measure probabilities w.r.t. one curve and calculate $r_n$ with respect to another. The distance $r_n$ is always measured according to the second expression $\overline J$, unless it is indicated otherwise.

A natural question is how $m_J(i\in \overline J: r_0(x)<\alpha)$ is related to $m_J(i\in J: r_0(x)<\alpha)$ as a function of $\alpha>0$. For any curve $J\in I$, $\overline J$ has at most two connected components, thus
\begin{align}\label{eq: 2-reasonable}
    m_J(i\in \overline J: r_0(x)<\alpha) \leq 2 m_J(i\in J: r_0(x)<\alpha)
\end{align}
holds for any $\alpha>0$.

\begin{lemma}[Growth lemma]\label{prop: growth lemma} There exists a $\overline\Lambda>1$ and a $\overline{C}, \overline{c} >0$ such that for all $\varepsilon, \varepsilon'>0$ small enough, all curves $J \subset I$,   and all $n\in\N^+$ we have the following estimates.

\begin{itemize}
\item[(a)]
For $J^{(L)} = J \cap I_L \cap T_\varepsilon^{-n}I_L$ and $J^{(R)} =  J \cap I_R \cap T_\varepsilon^{-n}I_R \, (=T_\varepsilon^{-n} \overline{ J^{(R)}_n})$,
\begin{align*}
    m_{J}\left(
    x\in T_\varepsilon^{-n}\overline{ J^{(L)}_n}\cup T_\varepsilon^{-n}\overline{ J^{(R)}_n}: r_n(x)<\varepsilon'
    \right)
    \leq&
    \overline C m_{J}(x\in J: r_0(x)<\varepsilon'/\overline\Lambda^n)
    +
    \overline c \varepsilon'.
\end{align*}

\item[(b)]
For $ J^{(L\to R)} =  J \cap I_L \cap T_\varepsilon^{-n}I_R \, (=T_\varepsilon^{-n} \overline{ J^{(L\to R)}_n})$,
\begin{align*}
    m_{J}\left(
    x\in T_\varepsilon^{-n}\overline{ J^{(L\to R)}_n}: r_n(x)<\varepsilon'
    \right)
    \leq&
    \overline C m_{J}(x\in J: r_0(x)<\varepsilon\varepsilon'/\overline\Lambda^n)
    +
    \overline c \varepsilon'
\end{align*}

\item[(c)]
For $ J^{(R\to L)} = J \cap I_R \cap T_\varepsilon^{-n}I_L$ we have
\begin{align*}
    m_{J}\left(
    x\in T_\varepsilon^{-n}\overline{ J^{(R\to L)}_n}: r_n(x)<\varepsilon'
    \right)
    \leq&
    \overline C m_{J}(x\in J: r_0(x)<\varepsilon'/\overline\Lambda^{n-N_\varepsilon})
    +
    \overline c \varepsilon'.
\end{align*}

\end{itemize}
\end{lemma}

Notice that in $(a)$ for example, if $T_\varepsilon^{-n} \overline {J^{(\xi)}_n}= J_\xi$ then the statement is the one we would expect from a regular system, see \cite[Chapter 5.9]{chaotic billiards}.
It is clear that we have to be careful about the set $B_\varepsilon$. Since $|DT|_{H_1}|\ll \varepsilon$, we have to wait $N_\varepsilon\asymp \log\varepsilon$ time for a density to regain its regularity.

Also notice that since $B_\varepsilon$ may cut $J_\xi$ into many connected components until time $n$, $r^{\left(T_\varepsilon^{-n} \overline {J^{(\xi)}_n}\right)}_n(x) \leq r^{(J_\xi)}_n(x)$ for $x \in T_\varepsilon^{-n} \overline {J^{(\xi)}_n}$

\begin{customproof}{Proof}
    The proof is structured as follows. First, we fix an $n_0\in\N^+$ (unrelated to the $n_0$ in the proof of Proposition \ref{prop: Lasota-Yorke (III)}) and construct a one-step estimate for $T_\varepsilon^{n_0}$. After considering what happens to a curve $J\subset H_1$ in the next $N_\varepsilon$ steps, we will modify our previous $n_0$-step estimate to an asynchronous backwards estimate for an evolved curve $\overline{J_m}$. This will mean that in one estimate we will have components pulled back to different points in time. Iterating this to get the general estimate directly would be hard to follow, since there would be numerous expressions but there would be only a few distinct types. For this reason, we will analyse the possible behaviours of points in $n$ steps, construct a backwards estimate for an arbitrary evolved curve $W_n$ having the same evolution, and finally counting how many times the corresponding expression would occur in the estimate for a general curve $J$ when using the asynchronous evolution. Utilizing this, we will construct an estimate for $T_\varepsilon^{-n}\overline {J_n}$ from which we can extract the three separate ones in the statement of the lemma.

    Let $\Lambda_* = \inf_{\varepsilon \geq 0}\inf_{I\setminus (H_0\cup H_1) |\D T_\varepsilon|}>1$, and set any $1<\Lambda<\Lambda_*$ and $n_0\in \N^+$ such that $\Lambda^{n_0}/(5\tau)> 1$.
    \\
    For an $x\in I$ and $i=0,1$ we define
    \begin{align*}
        K_i^{(0)}(x)&=\min\{k\in\N: T_\varepsilon^k x \in H_i\},
    \end{align*}
    as the first time $x$ enters $H_i$. By \ref{C2} and \ref{C4}, there is a $\delta>0$ such that for any $J_\xi\subset I_\xi$ with $Leb(J)\leq\delta$, if $k_J=\min_{J}\min\{K_0^{(0)}, K_1^{(0)}\}\leq n_0$, then $T^{k_J+1}J$ gets cut into at most $3$ parts, two of which stays in $I_\xi$ while the remaining one is evolved through the hole. If a curve $J\subset I$ satisfies $Leb(J)\leq \delta$, we will say that $J$ is a \emph{$\delta$-curve}.

    When evolving the $\delta$-curve $J$ by $n_0$ steps, we can partition it into
        \begin{align*}
        J=[J\cap (n_0\leq k_J)]\cup [J\cap (K_0^{(0)}<n_0\leq K_1^{(0)})]\cup [J\cap (K_0^{(0)}< K_1^{(0)} < n_0)]\cup [J\cap (K_1^{(0)}<n_0<K_0^{(0)})]
        \end{align*}
    We can look at the evolution of each partition element.
    \\
    First, $J\cap (n_0\leq k_J)$ has at most four components, two in $I_L$ and two in $I_R$. All of these $J^{(i)}$ components obey the following estimate.
    \begin{align*}
        m_{J^{(i)}}(x\in J^{(i)}: r_{n_0}(x)< \varepsilon')
        \leq&
        m_{J^{(i)}}(x\in J^{(i)}: r_{0}(x)< \varepsilon'/\Lambda^{n_0})
    \end{align*}
    All of these components might evolve to be longer than $\delta$, hence we cut them up into $\delta$-curves, which results in a $c\varepsilon'$ error (see e.g.~\cite[section 5.9]{chaotic billiards}).
        \begin{align*}
            m_{J\cap (n_0\leq k_J)}(x\in J\cap (n_0\leq k_J): r_{n_0}(x)< \varepsilon')
            \leq&
            4m_{J\cap (n_0\leq k_J)}(x\in J\cap (n_0\leq k_J): r_{0}(x)< \varepsilon'/\Lambda^{n_0})
            +
            c\varepsilon'.
        \end{align*}
    As for $J\cap (K_0^{(0)}<n_0\leq K_1^{(0)})$, the only component evolves in $I_L$ for $k_0 = k_{J\cap (K_0^{(0)}<n_0\leq K_1^{(0)})}$ steps and $T_\varepsilon^{k_0} [J\cap (K_0^{(0)}<n_0\leq K_1^{(0)})] \subset H_0$. After this, it gets enlargened by a factor $C^{-1}\varepsilon^{-1}\Lambda$ after which we cut it into $\delta^{-1}$ components to make them $\delta$-curves, all of which evolve in $I_R$ splitting into at most $2$ components. Finally, we cut the components into $\delta$-curves that introduces a $c\varepsilon'$ error.
    \begin{align*}
        &m_{J\cap (K_0^{(0)}<n_0\leq K_1^{(0)})}(x\in J\cap (K_0^{(0)}<n_0\leq K_1^{(0)}): r_{n_0}(x)<\varepsilon')
        \\
        \leq&
        \tau m_{[J\cap (K_0^{(0)}<n_0\leq K_1^{(0)})]_{k_0+1}}(x\in [J\cap (K_0^{(0)}<n_0\leq K_1^{(0)})]_{k_0+1}: r_{n_0-1-k_0}(x)<\varepsilon')
        \\
        \leq&
        2\tau m_{[J\cap (K_0^{(0)}<n_0\leq K_1^{(0)})]_{k_0+1}}(x\in [J\cap (K_0^{(0)}<n_0\leq K_1^{(0)})]_{k_0+1}: r_{0}(x)<\varepsilon'/\Lambda^{n_0-1-k_0})
        +
        c\varepsilon'
        \\
        \leq&
        2\tau^2 m_{[J\cap (K_0^{(0)}<n_0\leq K_1^{(0)})]_{k_0}}(x\in [J\cap (K_0^{(0)}<n_0\leq K_1^{(0)})]_{k_0}: r_{1}(x)<\varepsilon'/\Lambda^{n_0-1-k_0})
        +
        c\varepsilon'
        \\
        \leq&
        2\tau^2\delta^{-1} m_{[J\cap (K_0^{(0)}<n_0\leq K_1^{(0)})]_{k_0}}(x\in [J\cap (K_0^{(0)}<n_0\leq K_1^{(0)})]_{k_0}: r_{0}(x)<C\varepsilon\varepsilon'/\Lambda^{n_0-k_0})
        +
        c\varepsilon'
        \\
        \leq&
        2\tau^3\delta^{-1} m_{[J\cap (K_0^{(0)}<n_0\leq K_1^{(0)})]}(x\in [J\cap (K_0^{(0)}<n_0\leq K_1^{(0)})]: r_{k_0}(x)<C\varepsilon\varepsilon'/\Lambda^{n_0-k_0})
        +
        c\varepsilon'
        \\
        \leq&
        2\tau^3\delta^{-1} m_{[J\cap (K_0^{(0)}<n_0\leq K_1^{(0)})]}(x\in [J\cap (K_0^{(0)}<n_0\leq K_1^{(0)})]: r_{k_0}(x)<C\varepsilon\varepsilon'/\Lambda^{n_0})
        +
        c\varepsilon'
    \end{align*}
    The component $J\cap (K_0^{(0)}< K_1^{(0)} < n_0)$ evolves similarly, but after exiting $H_0$ and getting cut into $\delta^{-1}$ components, all of these components enter $H_1$ at some (possibly different) $k_1$ time. After this, they get contracted by a factor $C^{-1}\varepsilon\Lambda$ when evolving through $H_1$, and finally get mapped into $I_{n_0-k_1-1}$.
    \begin{align*}
        &m_{[J\cap (K_0^{(0)}< K_1^{(0)} < n_0)]}(x\in [J\cap (K_0^{(0)}< K_1^{(0)} < n_0)]: r_{n_0}(x)<\varepsilon')
        \\
        \leq&
        \tau m_{[J\cap (K_0^{(0)}< K_1^{(0)} < n_0)]_{k_1+1}}(x\in [J\cap (K_0^{(0)}< K_1^{(0)} < n_0)]_{k_1+1}: r_{0}(x)<\varepsilon'/\Lambda^{n_0-k_1-1})
        \\
        \leq&
        2\tau^2 m_{[J\cap (K_0^{(0)}< K_1^{(0)} < n_0)]_{k_1}}(x\in [J\cap (K_0^{(0)}< K_1^{(0)} < n_0)]_{k_1}: r_{0}(x)<\frac{C\varepsilon'}{\Lambda^{n_0-k_1}\varepsilon})
        \\
        \leq&
        2\tau^3 m_{[J\cap (K_0^{(0)}< K_1^{(0)} < n_0)]_{k_0+1}}(x\in [J\cap (K_0^{(0)}< K_1^{(0)} < n_0)]_{k_0+1}: r_{0}(x)<\frac{C\varepsilon'}{\Lambda^{n_0-k_0-1}\varepsilon})
        \\
        \leq&
        2\tau^4 \delta^{-1} m_{[J\cap (K_0^{(0)}< K_1^{(0)} < n_0)]_{k_0}}(x\in [J\cap (K_0^{(0)}< K_1^{(0)} < n_0)]_{k_0}: r_{0}(x)<\frac{C^2\varepsilon'}{\Lambda^{n_0-k_0}})
        \\
        \leq&
        2\tau^5 \delta^{-1} m_{[J\cap (K_0^{(0)}< K_1^{(0)} < n_0)]}(x\in [J\cap (K_0^{(0)}< K_1^{(0)} < n_0)]: r_{0}(x)<\frac{C^2\varepsilon'}{\Lambda^{n_0}})
    \end{align*}
    The last component, $J\cap (K_1^{(0)}<n_0<K_0^{(0)})$ evolves in $I_R$ for $k_1$ steps, evolves through $H_1$ then gets mapped into $I_{n_0-k_1-1}$.
    \begin{align*}
        &m_{[J\cap (K_1^{(0)}<n_0<K_0^{(0)})]}(x\in [J\cap (K_1^{(0)}<n_0<K_0^{(0)})]: r_{n_0}(x)<\varepsilon')
        \\
        \leq&
        \tau m_{[J\cap (K_1^{(0)}<n_0<K_0^{(0)})]_{k_1+1}}(x\in [J\cap (K_1^{(0)}<n_0<K_0^{(0)})]_{k_1+1}: r_{0}(x)<\varepsilon'/\Lambda^{n_0-k_1-1})
        \\
        \leq&
        \tau^2 m_{[J\cap (K_1^{(0)}<n_0<K_0^{(0)})]_{k_1}}(x\in [J\cap (K_1^{(0)}<n_0<K_0^{(0)})]_{k_1}: r_{0}(x)<\frac{C\varepsilon'}{\Lambda^{n_0-k_1}\varepsilon})
        \\
        \leq&
        \tau^3 m_{[J\cap (K_1^{(0)}<n_0<K_0^{(0)})]}(x\in [J\cap (K_1^{(0)}<n_0<K_0^{(0)})]: r_{0}(x)<\frac{C\varepsilon'}{\Lambda^{n_0}\varepsilon})
    \end{align*}

    Combining these four estimates we have
    \begin{align}\label{eq: growth lemma for n0 steps}
        m_{J}(x\in J: r_{n_0}(x)<\varepsilon')
        \leq&
        m_{J}(x\in J\cap (k_J\geq n_0): r_{n_0}(x)<\varepsilon') \nonumber
        \\&+
        m_{J}(x\in J\cap (K_0^{(0)}<n_0\leq K_1^{(0)}): r_{n_0}(x)<\varepsilon') \nonumber
        \\&+
        m_{J}(x\in J\cap (K_0^{(0)}< K_1^{(0)}< n_0): r_{n_0}(x)<\varepsilon') \nonumber
        \\&+
        m_{J}(x\in J\cap (K_1^{(0)}<n_0): r_{n_0}(x)<\varepsilon') \nonumber
        \\
        \leq&
        4m_{J}(x\in J: r_{0}(x)<\varepsilon'/\Lambda^{n_0}) + c\varepsilon' \nonumber
        \\&+
        2\tau^3\delta^{-1}m_{J}(x\in J: r_{0}(x)<C\varepsilon\varepsilon'/\Lambda^{n_0})
        \\&+
        2\tau^3\delta^{-1}m_{J}(x\in J: r_{0}(x)<C^2\varepsilon'/\Lambda^{n_0}) \nonumber
        \\&+
        2\tau^3 m_{J}(x\in J: r_{0}(x)<\frac{C\varepsilon'}{\Lambda^{n_0}\varepsilon}). \nonumber
    \end{align}

    The last two expressions represent the at most $2+2\delta^{-1}$ components that can enter $\cup_{k_0}^{n_0-1}I_k\subset B_\varepsilon$. In the further arguments of the proof we will use induction to prove an estimate for $k n_0$ steps for arbitrary $k\in\N^+$. In order to accomplish this, we need to establish how the components of an already evolved curve $J_m$ evolve backwards when we pull them back during the induction.
    Recall that $\overline{ J_m} = J_m\setminus B_\varepsilon$.
    We will include an additional expression representing the components exiting $B_\varepsilon$. These components develop according to $(T_\varepsilon|_{I_L})^{n_0}$, but we can assume that they do not include any critical point of it, and stay $\delta$-curves. In every step of the induction, a new instance of this expression will appear. These expressions have a strict evolution path. Throughout the induction, these expressions go backwards through $B_\varepsilon$, where they stay connected. Using \eqref{eq: expansion on irregular part}, we will evolve these components $\lfloor N_\varepsilon/n_0\rfloor n_0$ steps in one induction step. We will denote $\lfloor N_\varepsilon/n_0\rfloor+1$ by $M_\varepsilon$

    \begin{align}\label{eq: growth lemma asynchronous n_0 steps}
        m_{J_{m+n_0}}(x\in \overline{ J_{m+n_0}}: r_{0}(x)<\varepsilon')
        \leq&
        \tau m_{J_m}(x\in \overline {J_m}\cap (k_J\geq n_0): r_{n_0}(x)<\varepsilon') \nonumber
        \\&+
        \tau m_{J_m}(x\in \overline {J_m}\cap (K_0^{(0)}<n_0\leq K_1^{(0)}): r_{n_0}(x)<\varepsilon') \nonumber
        \\&+
        \tau m_{J_{m+n_0 - (M_\varepsilon-1) n_0}}(x\in  J_{m+n_0 -(M_\varepsilon - 1)n_0}\cap \cup_{k_0}^{n_0-1}I_k: r_{(M_\varepsilon - 1 ) n_0}(x)<\varepsilon') \nonumber
        \\
        \leq&
        4\tau m_{J_m}(x\in \overline {J_m}: r_{0}(x)<\varepsilon'/\Lambda^{n_0}) +c\varepsilon' \nonumber
        \\&+
        2\tau^4\delta^{-1}m_{J_m}(x\in \overline {J_m}: r_{0}(x)<C\varepsilon\varepsilon'/\Lambda^{n_0}) +c\varepsilon' \nonumber
        \\&+
        \tau m_{J_{m+n_0 - (M_\varepsilon-1) n_0}}(x\in J_{m+n_0 - (M_\varepsilon-1) n_0}\cap \cup_{k_0}^{n_0-1}I_k: r_{0}(x)<\varepsilon'/\Lambda^{(M_\varepsilon-1) n_0}) \nonumber
        \\
        \leq&
        4\tau m_{J_m}(x\in \overline {J_m}: r_{0}(x)<\varepsilon'/\Lambda^{n_0}) +c\varepsilon' \nonumber
        \\&+
        2\tau^4\delta^{-1}m_{J_m}(x\in \overline {J_m}: r_{0}(x)<C\varepsilon\varepsilon'/\Lambda^{n_0}) +c\varepsilon' \nonumber
        \\&+
        \tau^2 m_{J_{m + n_0 - M_\varepsilon n_0}}(x\in J_{m + n_0 - M_\varepsilon n_0}\cap T_\varepsilon^{-n_0}\cup_{k_0}^{n_0-1}I_k: r_{n_0}(x)<\varepsilon'/\Lambda^{(M_\varepsilon-1) n_0}) \nonumber
        \\
        \leq&
        4\tau m_{J_m}(x\in \overline {J_m}: r_{0}(x)<\varepsilon'/\Lambda^{n_0}) +c\varepsilon' \nonumber
        \\&+
        2\tau^4\delta^{-1}m_{J_m}(x\in \overline {J_m}: r_{0}(x)<C\varepsilon\varepsilon'/\Lambda^{n_0}) +c\varepsilon'
        \\&+
        2\tau^5\delta^{-1}m_{J_{m + n_0 - M_\varepsilon n_0}}(x\in \overline {J_{m + n_0 - M_\varepsilon n_0}}: r_{0}(x)<C^2\varepsilon'/\Lambda^{M_\varepsilon n_0}) \nonumber
        \\&+
        2\tau^5 m_{J_{m + n_0 - M_\varepsilon n_0}}(x\in \overline {J_{m + n_0 - M_\varepsilon n_0}}: r_{0}(x)<\frac{C\varepsilon'}{\Lambda^{M_\varepsilon n_0}\varepsilon}) \nonumber
    \end{align}
    Notice that in the last two expressions we evolved $J_{m+n_0}$ backwards by $\lfloor N_\varepsilon/n_0\rfloor n_0+n_0$ steps. If $m<\lfloor N_\varepsilon/n_0\rfloor n_0$, these expressions are interpreted as zero. Also notice that the second expression represents components that got mapped through $H_0$ in the last $n_0$ steps, hence they will not be evolved by the same rule until they get pulled through $B_\varepsilon$.

    Notice that when iterating \eqref{eq: growth lemma asynchronous n_0 steps} we see more and more expressions arise, which are probabilities (and error terms) measured at times $n-kn_0$ for some $k\leq n/n_0$, and time zero. We can visualize the non-error expressions during the iteration as nodes of a tree seen in Figure \ref{fig: evolution tree}.

\begin{figure}[H]
    \centering
    \begin{tikzpicture}

        \node at (0, 3) {$m_{J_{n}}(x\in \overline{ J_{n}}: r_{0}(x)<\varepsilon')$};

        \draw[] (2, 3) -- (4, 6);
        \draw[] (2, 3) -- (4, 4);
        \draw[] (2, 3) -- (4, 2);
        \draw[] (2, 3) -- (4, 0);

        \node at (5.8, 6) {$m_{J_{n-n_0}}(x\in \overline{ J_{n-n_0}} ...)$};
        \node at (5.8, 4) {$m_{J_{n-n_0}}(x\in \overline{ J_{n-n_0}} ...)$};
        \node at (6.2, 2) {$m_{J_{n-M_\varepsilon n_0}}(x\in \overline{ J_{n-M_\varepsilon n_0}} ...)$};
        \node at (6.2, 0) {$m_{J_{n-M_\varepsilon n_0}}(x\in \overline{ J_{n-M_\varepsilon n_0}} ...)$};

        \draw[] (8.3, 0) -- (9, 0.5);
        \draw[] (8.3, 0) -- (9, 0.16);

        \draw[] (8.3, 2) -- (9, 2.5);
        \draw[] (8.3, 2) -- (9, 1.83);
        \draw[] (8.3, 2) -- (9, 1.5);

        \draw[] (7.5, 4) -- (9, 5);
        \draw[] (7.5, 4) -- (9, 3.66);
        \draw[] (7.5, 4) -- (9, 3);

        \draw[] (7.5, 6) -- (9, 6.5);
        \draw[] (7.5, 6) -- (9, 6.16);
        \draw[] (7.5, 6) -- (9, 5.83);
        \draw[] (7.5, 6) -- (9, 5.5);

        \node at (10, 3) {$\dots \quad \dots$};

        \draw[] (10.9, 8) -- (10.5, 7.75);
        \node at (12, 8) {$m_{J}(x\in \overline{ J} ...)$};

        \draw[] (10.9, 7.25) -- (10.5, 7.2);
        \node at (12, 7.25) {$m_{J}(x\in \overline{ J} ...)$};

        \node at (12, 6) {$\vdots$};
        \node at (12, 0) {$\vdots$};

        \draw[] (10.9, -1) -- (10.5, -0.75);
        \node at (12, -1) {$m_{J}(x\in \overline{ J} ...)$};

    \end{tikzpicture}
    \caption{The visualization of the tree of estimates throughout the iteration of \eqref{eq: growth lemma asynchronous n_0 steps}}\label{fig: evolution tree}
    \end{figure}
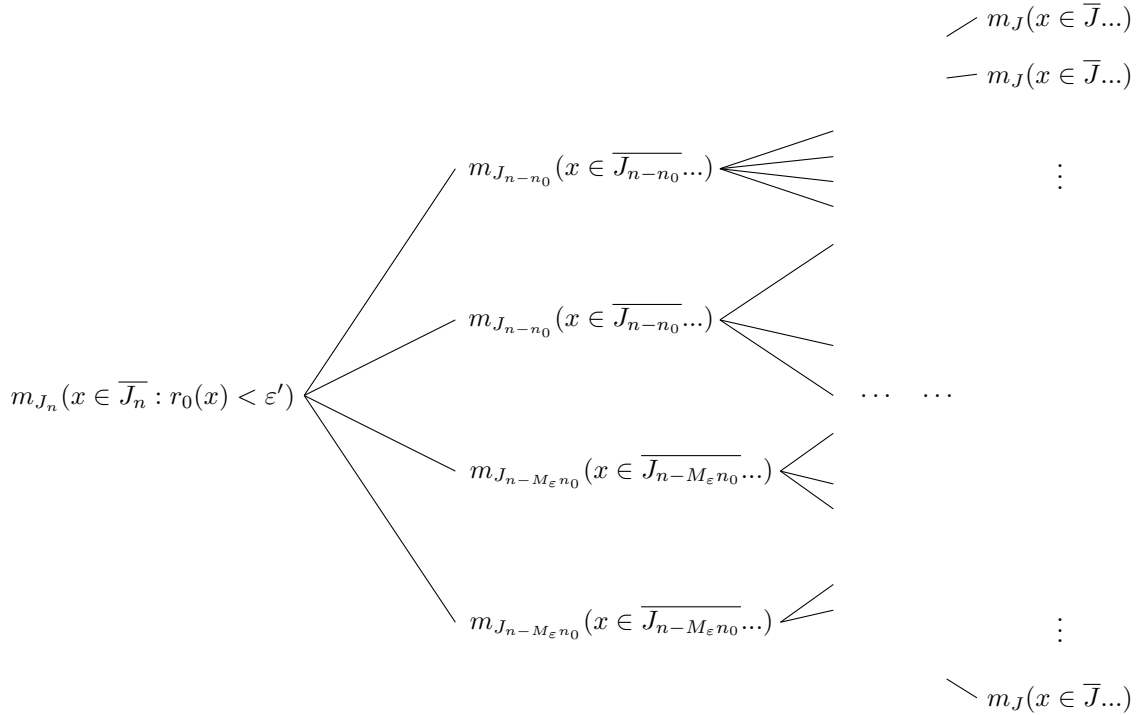
    We aim to avoid presenting the iteration of \eqref{eq: growth lemma asynchronous n_0 steps} directly due to the complexity of the final estimate, hence we propose a rigorous workaround. During the iteration the different expressions represent points of $\overline{J_n}$ that underwent different (backwards) evolution. We will define the notion and the set of evolution futures $\mathfrak{F}$ below.  By deconstructing $J$ into components with a given evolution future $f\in \mathfrak{F}$, represented by different paths in the evolution tree in Figure \ref{fig: evolution tree}, we arrive at subsets $W^{f}\subset J$. Then, when iterating \eqref{eq: growth lemma asynchronous n_0 steps} on $\overline{W^{f}_n}$, we only get one expression in each step, thus we can easily construct a growth estimate for it. This will be one expression present in the last layer of the evolution tree in Figure \ref{fig: evolution tree}. By summing the estimates for all possible evolution futures $f \in \mathfrak{F}$ we get an estimate for $\overline{J_n}$.

    For this reason we will categorize subsets $W$ of $J$ that are the collections of points represented by a path in the evolution tree. This is equivalent to saying that for all \emph{appropriate} $k\leq n/n_0$, one of the following four statements hold:
    \begin{align*}
        W_{n-(k-1)n_0}
        \subset&
        (I_L \cap T_\varepsilon^{n_0}I_L)\cup (I_R \cap T_\varepsilon^{n_0}I_R),
        \\
    W_{n-(k-1)n_0}
    \subset&
    I_R \cap T_\varepsilon^{n_0}I_L,
    \\
    W_{n-(k-1)n_0}
    \subset&
    I_L \cap T_\varepsilon^{M_\varepsilon n_0}I_R,
    \\
    W_{n-(k-1)n_0}
    \subset&
    \{x\in I_L: \exists l < n_0: T_\varepsilon^l x \in I_R, T_\varepsilon^{n_0}x\in B_\varepsilon\}
    \subset I_L\cap T_\varepsilon^{M_\varepsilon n_0}I_L.
    \end{align*}
    (We say \emph{appropriate}, since if $W_{n-(k-1)n_0}\subset I_L \cap T_\varepsilon^{M_\varepsilon n_0}I_R$, the next timestamp for $W$ is $n-(k-1-M_\varepsilon)n_0$.)
    If $W$ satisfies this, we will say that it has a \emph{joint future} up to time $n$.

    In order to precisely describe all possible joint futures, on top of $K_i^{(0)}$, we introduce the functions
    \begin{align*}
        K_i^{(j)}(x)&=\min\{k\geq K_i^{(j-1)}: T_\varepsilon^k x \in H_i\}
    \end{align*}
    for $i=0, 1$ and $j\in \N^+$.

    If $\min_{W}K^{(0)}<n$, then the evolution of $W$ can be encoded by two elements of information. The first hole $W$ entered during its evolution $i=\arg\min_{j=0,1}K_j^{(0)}$, and the vector $l=(l_k, l_{k-1}, \dots, l_0)$ of times it entered a hole, which means $n-l_{2t}n_0\leq K_i^{(t)}(x) < n-(l_{2t}-1)n_0$ and $n-l_{2t+1}n_0\leq K_{1-i}^{(t)}(x) < n-(l_{2t+1}-1)n_0$ for all $x\in W$. If $W$ is a curve with joint future  can be encoded by the pair $(i, l)$, we say that $W$ has a \emph{joint future of type $(i,l)$}. If $\min_{W}K^{(0)}\geq n$, we say that $W$ has a joint future of type $0$. The set of all possible joint future types $\mathfrak{F}$ is a subset of
    $\overline{\mathfrak{F}} = \left\{ (i, l)\in \{0, 1\}\times\cup_{j=1}^{\infty} (\N^+)^j: 1\leq l_{t}\leq l_{t-1}\leq n/n_0 \; \forall t \right\}\cup \{0\}$\footnote{Any visit to $H_0$ must be preceded by at least $N_\varepsilon$ steps, but a curve can hit $H_0$ and then $H_1$ in $n_0$ steps.}.  When iterating \eqref{eq: growth lemma asynchronous n_0 steps}, each $ f \in \mathfrak{F}$ has exactly one expression representing the estimate for  $W^{f}= \cup\{W\subset J: W \,\text{ has a joint future of type } \, f\}$. Since  $J=\cup_{f\in\mathfrak{F}}W^{f}$, we will derive an estimate for all future types $f\in \mathfrak{F}$ then sum these estimates to construct one for $J$.

    It is important to note that even though $W^{(f)}$ may have several connected components, at the end of iteration we will see the probability $m_J$ and distance $r^{\overline J}$ calculated with respect to $J$ and $\overline J$, as is indicated in Figure \ref{fig: evolution tree}. This is because the disconnected sets $W^{(f)}$ are "adapted" to the estimate of the evolution under $T_\varepsilon$. This means that every separation between components of $W^{(f)}$ would be accounted for at some time in the estimate either in the exponential bound or the error term.

    First, suppose that $W = W^{(0)}$. Then,

    \begin{align}\label{eq: growth lemma stays in place}
        m_{W_n}(x\in  \overline{ W_n}: r_0(x)<\varepsilon')
        \leq
        (4\tau)^{1+n/n_0}m_{J}(x\in  \overline J: r_0(x)<\varepsilon'/\Lambda^{n})
        +
        c\sum_{l=0}^{n/n_0}\left(\frac{4\tau}{\Lambda^{n_0}}\right)^l\varepsilon'.
    \end{align}

    Let $W=W^{(1, l_0)}$, which has the joint future $0\leq K_1^{(0)}< n \leq K_0^{(0)}$.
    \begin{align*}
        m_{W_n}(x\in  \overline {W_n}: r_0(x)<\varepsilon')
        \leq&
        \frac{2\tau^3}{(4\tau)^{N_\varepsilon/n_0}}(4\tau)^{(n-K_1^{(0)})/n_0} m_{W_{K_1^{(0)}}}(x\in \overline {W_{K_1^{(0)}}}: r_0(x)<\frac{C\varepsilon'}{\Lambda^{n-K_1^{(0)}}\varepsilon})
        \\&+
        c\sum_{l=0}^{(n-K_1^{(0)}-N_\varepsilon)/n_0-1}\left(\frac{4\tau}{\Lambda^{n_0}}\right)^l\varepsilon'
        \\
        \leq&
         \frac{2\tau^3}{(4\tau)^{N_\varepsilon/n_0}}(4\tau)^{1+n/n_0} m_{J}(x\in \overline J: r_0(x)<\frac{C\varepsilon'}{\Lambda^{n}\varepsilon})
        \\&+
        c \frac{2C\tau^3}{(4\tau)^{N_\varepsilon/n_0}\varepsilon}\sum_{l=(n-K_1^{(0)})/n_0}^{n/n_0}\left(\frac{4\tau}{\Lambda^{n_0}}\right)^l\varepsilon'
        \\&+
        c\sum_{l=0}^{(n-K_1^{(0)}-N_\varepsilon)/n_0-1}\left(\frac{4\tau}{\Lambda^{n_0}}\right)^l\varepsilon'
        \end{align*}
        We will show that for an appropriate $\gamma<1$,  $\varepsilon^{-1}\leq\Lambda^{N_\varepsilon}$. Using \eqref{eq: Na}, we can see that
    \begin{align*}
        N_\varepsilon \log \Lambda
        \leq&
        |\log \varepsilon |
        \\
        \frac{N_\varepsilon \log \Lambda}{|\log \varepsilon |}
        \leq&
        \frac{(2-\gamma)\log\Lambda}{\log\inf_{I\setminus H_\varepsilon}|\D T_\varepsilon|}
        < 1
    \end{align*}
    holds for $\gamma<1$ large enough.

    Hence, continuing the previous estimate,

    \begin{align}\label{eq: growth lemma right to left}
        m_{W_n}(x\in  \overline {W_n}: r_0(x)<\varepsilon')
        \leq&
         \frac{2\tau^3}{(4\tau)^{N_\varepsilon/n_0}}(4\tau)^{1+n/n_0} m_{J}(x\in \overline J: r_0(x)<\frac{C\varepsilon'}{\Lambda^{n-N_\varepsilon}}) \nonumber
        \\&+
        c 2C\tau^3\sum_{l=(n-K_1^{(0)}-N_\varepsilon)/n_0}^{(n-N_\varepsilon)/n_0 - 1}\left(\frac{4\tau}{\Lambda^{n_0}}\right)^l\varepsilon'  \nonumber
        \\&+
        c\sum_{l=0}^{(n-K_1^{(0)}-N_\varepsilon)/n_0}\left(\frac{4\tau}{\Lambda^{n_0}}\right)^l\varepsilon'  \nonumber
        \\
        \leq&
        \frac{2\tau^3}{(4\tau)^{N_\varepsilon/n_0}}(4\tau)^{1+n/n_0} m_{J}(x\in \overline J: r_0(x)<\frac{C\varepsilon'}{\Lambda^{n-N_\varepsilon}})
        \\&+
        c2C\tau^3\sum_{l=0}^{(n-N_\varepsilon)/n_0}\left(\frac{4\tau}{\Lambda^{n_0}}\right)^l\varepsilon'  \nonumber
    \end{align}

    Similarly, if $W=W^{(0, l_0)}$ with the joint future $0\leq K_0^{(0)}< n \leq K_1^{(0)}$, then
    \begin{align}\label{eq: growth lemma left to right}
        m_{W_n}(x\in  \overline {W_n}: r_0(x)<\varepsilon')
        \leq&
        \tau^3\delta^{-1}(4\tau)^{(n-K_0^{(0)})/n_0} m_{W_{K_0^{(0)}}}(x\in \overline {W_{K_0^{(0)}}}: r_0(x)<\frac{C\varepsilon\varepsilon'}{\Lambda^{n-K_0^{(0)}}})  \nonumber
        \\&+
        c\sum_{l=0}^{(n-K_0^{(0)})/n_0-1}\left(\frac{4\tau}{\Lambda^{n_0}}\right)^l\varepsilon' \nonumber
        \\
        \leq&
        \tau^3\delta^{-1}(4\tau)^{1+n/n_0} m_{J}(x\in \overline J: r_0(x)<\frac{C\varepsilon\varepsilon'}{\Lambda^{n}})
        \\&+
        c\sum_{l=0}^{n/n_0}\left(\frac{4\tau}{\Lambda^{n_0}}\right)^l\varepsilon' \nonumber
    \end{align}

    The set $W = W^{(0, l)}$, which has the joint future $0\leq K_0^{(0)}<K_1^{(0)}<\dots < K_1^{(j)}< n \leq K_0^{(j+1)}$ has an upper bound

    \begin{align}\label{eq: growth lemma left to left}
        m_{W_n}(x\in \overline {W_n}: r_0(x)<\varepsilon')
        \leq&
        \frac{2\tau^3}{(4\tau)^{N_\varepsilon/n_0}}(4\tau)^{(n-K_1^{(j)})/n_0} m_{W_{K_1^{(j)}}}(x\in \overline {W_{K_1^{(j)}}}: r_0(x)<\frac{C\varepsilon'}{\Lambda^{n-K_1^{(j)}}\varepsilon})  \nonumber
        \\&+
        c\sum_{l=0}^{(n-K_1^{(j)}-N_\varepsilon)/n_0-1}\left(\frac{4\tau}{\Lambda^{n_0}}\right)^l\varepsilon' \nonumber
        \\
        \leq&
        \frac{2\delta^{-1}\tau^6}{(4\tau)^{N_\varepsilon/n_0}}(4\tau)^{(n-K_0^{(j)})/n_0} m_{W_{K_0^{(j)}}}\left(x\in \overline {W_{K_0^{(j)}}}: r_0(x)<\frac{C^2\varepsilon'}{\Lambda^{n-K_0^{(j)}}}\right)  \nonumber
        \\&+
        c 2C\tau^3\sum_{l=(n-K_1^{(j)}-N_\varepsilon)/n_0}^{(n-K_0^{(j-1)}-N_\varepsilon)/n_0-1}\left(\frac{4\tau}{\Lambda^{n_0}}\right)^l\varepsilon' \nonumber
        \\&+
        c\sum_{l=0}^{(n-K_1^{(j)}-N_\varepsilon)/n_0-1}\left(\frac{4\tau}{\Lambda^{n_0}}\right)^l\varepsilon' \nonumber
        \\
        \leq&
        \left(\frac{2\delta^{-1}\tau^6}{(4\tau)^{N_\varepsilon/n_0}}\right)^{j}(4\tau)^{(n-K_0^{(0)})/n_0} m_{W_{K_0^{(0)}}}\left(x\in \overline {W_{K_0^{(0)}}}: r_0(x)<\frac{C^{2j}\varepsilon'}{\Lambda^{n-K_0^{(0)}}}\right)  \nonumber
        \\&+
        c2C\tau^3\sum_{l=0}^{(n-K_0^{(0)}-(j+1)N_\varepsilon)/n_0-1}\left(\frac{4\tau}{\Lambda^{n_0}}\right)^l\varepsilon' \nonumber
        \\
        \leq&
        \left(\frac{2C^2\delta^{-1}\tau^6}{(4\tau)^{N_\varepsilon/n_0}}\right)^{j}(4\tau)^{1+n/n_0} m_{J}\left(x\in \overline J: r_0(x)<\frac{\varepsilon'}{\Lambda^{n}}\right)
        \\&+
        c2C\tau^3\sum_{l=0}^{n/n_0}\left(\frac{4\tau}{\Lambda^{n_0}}\right)^l\varepsilon' \nonumber
    \end{align}

    Similarly, if $W = W^{(0, l)}$ with joint future $0\leq K_1^{(0)}<K_0^{(0)}<\dots < K_0^{(j)}< n \leq K_1^{(j+1)}$, then

    \begin{align}\label{eq: growth lemma right to right}
        m_{W_n}(x\in \overline {W_n}: r_0(x)<\varepsilon')
        \leq&
        \left(\frac{2C^2\delta^{-1}\tau^6}{(4\tau)^{N_\varepsilon/n_0}}\right)^{j}(4\tau)^{1+n/n_0} m_{J}\left(x\in \overline J: r_0(x)<\frac{\varepsilon'}{\Lambda^{n}}\right)
        \\&+
        c2C\tau^3\sum_{l=0}^{n/n_0}\left(\frac{4\tau}{\Lambda^{n_0}}\right)^l\varepsilon'. \nonumber
    \end{align}

    If $W = W^{(0, l)}$ has the joint future $0\leq K_1^{(0)}<K_0^{(0)}<\dots < K_1^{(j)}< n \leq K_0^{(j)}$, then

    \begin{align}\label{eq: growth lemma right to left rounds}
        &m_{W_n}(x\in \overline {W_n}: r_0(x)<\varepsilon')
        \leq
        \\
        \leq&
        \left(\frac{2C^2\delta^{-1}\tau^6}{(4\tau)^{N_\varepsilon/n_0}}\right)^{j}(4\tau)^{(n-K_0^{(0)})/n_0} m_{W_{n-K_0^{(0)}}}\left(x\in \overline {W_{n-K_0^{(0)}}}: r_0(x)<\frac{\varepsilon'}{\Lambda^{n-K_0^{(0)}}}\right)  \nonumber
        \\&+
        c2C\tau^3\sum_{l=0}^{(n-K_0^{(0)}-jN_\varepsilon)/n_0-1}\left(\frac{4\tau}{\Lambda^{n_0}}\right)^l\varepsilon' \nonumber
        \\
        \leq&
        \frac{2\tau^3}{(4\tau)^{N_\varepsilon/n_0}}
        \left(\frac{2\delta^{-1}\tau^6}{(4\tau)^{N_\varepsilon/n_0}}\right)^{j}
        (4\tau)^{(n-K_1^{(0)})/n_0} m_{W_{n-K_1^{(0)}}}
        \left(x\in \overline {W_{n-K_1^{(0)}}}: r_0(x)<\frac{C\varepsilon'}{\Lambda^{n-K_1^{(0)}-N_\varepsilon}}\right)  \nonumber
        \\&+
        c2C\tau^3\sum_{l=0}^{(n-K_1^{(0)}-(j+1)N_\varepsilon)/n_0-1}\left(\frac{4\tau}{\Lambda^{n_0}}\right)^l\varepsilon' \nonumber
        \\
        \leq&
        \frac{2\tau^3}{(4\tau)^{N_\varepsilon/n_0}}
        \left(\frac{2\delta^{-1}\tau^6}{(4\tau)^{N_\varepsilon/n_0}}\right)^{j}
        (4\tau)^{1+n/n_0}
        m_{J}\left(x\in \overline J: r_0(x)<\frac{C\varepsilon'}{\Lambda^{n-N_\varepsilon}}\right)
        \\&+
        c2C\tau^3\sum_{l=0}^{n/n_0}\left(\frac{4\tau}{\Lambda^{n_0}}\right)^l\varepsilon' \nonumber
    \end{align}

    If $W = W^{(0, l)}$ has the joint future $0\leq K_0^{(0)}<K_1^{(0)}<\dots < K_0^{(j)}< n \leq K_1^{(j)}$, then

    \begin{align}\label{eq: growth lemma left to right rounds}
        m_{W_n}(x\in \overline {W_n}: r_0(x)<\varepsilon')
        \leq&
        \tau^3\delta^{-1}
        \left(\frac{2C^2\delta^{-1}\tau^6}{(4\tau)^{N_\varepsilon/n_0}}\right)^{j}(4\tau)^{1+n/n_0} m_{J}\left(x\in \overline J: r_0(x)<\frac{C\varepsilon\varepsilon'}{\Lambda^{n}}\right)
        \\&+
        c2C\tau^3\sum_{l=0}^{n}\left(\frac{4\tau}{\Lambda^{n_0}}\right)^l\varepsilon'. \nonumber
    \end{align}

    We showed what types of expressions would arise when iterating the asynchronous estimate \eqref{eq: growth lemma asynchronous n_0 steps}. Notice that the estimates only depend on $i=0,1$ and the dimension of $l$. Hence we define

    \begin{align*}
        J^{(0)}
        =&
        J\cap (n\leq K^{(0)})
        =
        W^{(0)},
    \end{align*}
    Now fix a $j\in \N^+$.
    \begin{align*}
    J^{(j)}_{L}
        =&
        J\cap(0\leq K_0^{(0)}<K_1^{(0)}<\dots < K_1^{(j)}< n \leq K_0^{(j+1)})
        =
        \bigcup_{\substack{(0, l)\in \mathfrak{F} \\ \dim l = 2j}}W^{(0, l)},
    \end{align*}
    \begin{align*}
        J^{(j)}_{R}
        =&
        J\cap(0\leq K_1^{(0)}<K_0^{(0)}<\dots < K_0^{(j)}< n \leq K_1^{(j+1)})
        =
        \bigcup_{\substack{(1, l)\in \mathfrak{F} \\ \dim l = 2j}}W^{(1, l)},
    \end{align*}
    For any $j\in \N$.
    \begin{align*}
        J^{(j)}_{L\to R}
        =&
        J\cap(0\leq K_0^{(0)}<K_1^{(0)}<\dots < K_1^{(j)}< n \leq K_1^{(j)})
        =
        \bigcup_{\substack{(0, l)\in \mathfrak{F} \\ \dim l = 2j+1}}W^{(0, l)},
    \end{align*}
    \begin{align*}
        J^{(j)}_{R\to L}
        =&
        J\cap(0\leq K_1^{(0)}<K_0^{(0)}<\dots < K_1^{(j)}< n \leq K_0^{(j)})
        =
        \bigcup_{\substack{(1, l)\in \mathfrak{F} \\ \dim l = 2j+1}}W^{(1, l)},
    \end{align*}

    By subadditivity of the probability measure $m_J$ the only further information we need is the size of the index sets of the unions in the above partitions.

    We will start with $\{(0, l)\in \mathfrak{F}: \dim l = 2j\}$ and $\{(1, l)\in \mathfrak{F}: \dim l = 2j\}$
    We can neglect the fact that after a visit to $H_1$, at least $N_\varepsilon$ steps have to precede $n$ and the next visit to $H_0$. Recall that a point can visit both $H_0$ and $H_1$ in $n_0$ steps (only in this order).
    Hence, we are trying to place $2j$ visits in $n/n_0+1$ units of time, but every second visit can be in the same unit time as the previous. We will account for this by imagining that we are placing visits in $n/n_0+1+j$ places. If a visit for $H_1$ is in the unit time after the visit for $H_0$, it means that they occured in the same $n_0$ steps.
    From this, we can see that  $\#\{(i, l)\in \mathfrak{F}: \dim l = 2j\} \leq \binom{n/n_0+1+j}{2j}$ for both $i=0,1$.

    By the same logic $\#\{(i, l)\in \mathfrak{F}: \dim l = 2j+1\}\leq \binom{n/n_0 + 1 + j}{2j+1}$ for $i=0,1$.

    Combining everything we get

    \begin{align*}
        &m_{J}\left(x\in T_\varepsilon^{-n}\overline {J_n}: r_n(x)<\varepsilon'\right)
        \leq
        \tau
        m_{J_n}\left(x\in \overline {J_n}: r_0(x)<\varepsilon'\right)\leq
        \\
        \leq&
        \tau m_J\left(x\in
        \overline {J^{(0)}_n}
        \cup
        \bigcup_{i=1}^{n/n_0+1}
        \overline {(J^{(i)}_L)_n}
        \cup \overline {(J^{(i)}_R)_n}
        \cup
        \bigcup_{j=0}^{n/n_0+1}
        \overline {(J^{(j)}_{L \to R})_n}
        \cup \overline {(J^{(j)}_{R \to L})_n}:
        r_0(x)<\varepsilon'
        \right)
        \\=&
        \tau m_J\left(x\in \overline {J^{(0)}_n}: r_0(x)<\varepsilon'\right)
        +
        \tau \sum_{i=1}^{n/n_0+1}
        m_J\left(x\in \overline {(J^{(i)}_L)_n}: r_0(x)<\varepsilon'\right)
        +
        m_J\left(x\in \overline {(J^{(i)}_R)_n}: r_0(x)<\varepsilon'\right)
        \\&+
        \tau \sum_{j=0}^{n/n_0+1}
        m_J\left(x\in \overline {(J^{(i)}_{L \to R})_n}: r_0(x)<\varepsilon'\right)
        +
        m_J\left(x\in \overline {(J^{(i)}_{R \to L})_n}: r_0(x)<\varepsilon'\right)
        \\
        \leq&
        \tau(4\tau)^{1+n/n_0}m_{J}\left(x\in  \overline J: r_0(x)<\varepsilon'/\Lambda^{n}\right)
        \\&+
        2\tau\sum_{i=1}^{n/n_0+1}
        \binom{n/n_0 +1 +i}{2i}\left(\frac{2C^2\delta^{-1}\tau^6}{(4\tau)^{N_\varepsilon/n_0}}\right)^{i}(4\tau)^{1+n/n_0} m_{J}\left(x\in \overline J: r_0(x)<\frac{\varepsilon'}{\Lambda^{n}}\right)
        \\&+
        \sum_{j=0}^{n/n_0+1}
        \binom{n/n_0 +1 +j}{2j+1}\tau^4\delta^{-1}
        \left(\frac{2C^2\delta^{-1}\tau^6}{(4\tau)^{N_\varepsilon/n_0}}\right)^{j}(4\tau)^{1+n/n_0} m_{J}\left(x\in \overline J: r_0(x)<\frac{C\varepsilon\varepsilon'}{\Lambda^{n}}\right)
        \\&+
        \binom{n/n_0 +1 +j}{2j+1}\frac{2\tau^4}{(4\tau)^{N_\varepsilon/n_0}}
        \left(\frac{2\delta^{-1}\tau^6}{(4\tau)^{N_\varepsilon/n_0}}\right)^{j}
        (4\tau)^{1+n/n_0}
        m_{J}\left(x\in \overline J: r_0(x)<\frac{C\varepsilon'}{\Lambda^{n-N_\varepsilon}}\right)
        \\&+
        c2C\tau^4(1-4\tau/\Lambda^{n_0})^{-1}\varepsilon'.
    \end{align*}
     We need to prove that the sums behave well. For this, let $p_\varepsilon^2= \frac{2C^2\delta^{-1}\tau^6}{(4\tau)^{N_\varepsilon/n_0}}$.
     Then for any $j\leq n/n_0+1$,
     \begin{align*}
         \max\left\{\binom{n/n_0+1+j}{2j}, \binom{n/n_0+1+j}{2j+1}\right\}
         \leq&
         \frac{(n/n_0 + 1 + j)^{2j+1}}{(2j)!}
         \leq
         \frac{(2(n/n_0+1))^{2j+1}}{(2j)!}
     \end{align*}
     For this reason, it is enough to show appropriate estimates for the right hand side.

     \begin{align*}
         \sum_{j=0}^{n/n_0+1}
         \frac{(2(n/n_0+1))^{2j+1}}{(2j)!}p_\varepsilon^{2j}
         \leq&
         2(n/n_0+1)\exp\left( 2(n/n_0+1) p_\varepsilon \right)
     \end{align*}

     \begin{align*}
         2(n/n_0+1)\exp \left( 2(n/n_0+1) p_\varepsilon \right)
         (4\tau)^{n/n_0}/\Lambda^n
         \leq&
         2(n/n_0+1)e^{2p_\varepsilon}\left(\frac{4\tau e^{2p_\varepsilon}}{\Lambda^{n_0}}\right)^{n/n_0}
         \leq
         \tilde C /\tilde\Lambda^{n}
     \end{align*}
    for an appropriate $\tilde C >0$ and $\tilde \Lambda >1$ for all $\varepsilon>0$ small enough.

    Thus we can see that

    \begin{align*}
        &m_J(x\in \overline J: r_n(x)<\varepsilon')
        \\
        \leq&
        \tilde C \tau(4\tau)m_{J}(x\in  \overline J: r_0(x)<\varepsilon'/\tilde \Lambda^{n})
        \\&+
        2\tilde C \tau(4\tau) m_{J}\left(x\in \overline J: r_0(x)<\frac{\varepsilon'}{\tilde \Lambda^{n}}\right)
        \\&+
        \tilde C \tau^4\delta^{-1}(4\tau) m_{J}(x\in \overline J: r_0(x)<\frac{C\varepsilon\varepsilon'}{\tilde \Lambda^{n}})
        \\&+
        \tilde C \frac{2\tau^4}{(4\tau)^{N_\varepsilon/n_0}}(4\tau) m_{J}(x\in \overline J: r_0(x)<\frac{C\varepsilon'}{\tilde \Lambda^{n-N\varepsilon}})
        \\&+
        c2C\tau^3(1-4\tau/\Lambda^{n_0})^{-1}\varepsilon'
        \\
        \leq&
        12\tilde C \tau^2 m_{J}(x\in  \overline J: r_0(x)<\varepsilon'/\tilde \Lambda^{n})
        \\&+
        4\tilde C \tau^5\delta^{-1} m_{J}(x\in \overline J: r_0(x)<\frac{C\varepsilon\varepsilon'}{\tilde \Lambda^{n}})
        \\&+
        \tilde C \frac{8\tau^5}{(4\tau)^{N_\varepsilon/n_0}}
        m_{J}\left(x\in \overline J: r_0(x)<\frac{C\varepsilon'}{\tilde\Lambda^{n-N_\varepsilon}}\right)
        \\&+
        c2C\tau^4(1-4\tau/\Lambda^{n_0})^{-1}\varepsilon'.
    \end{align*}
    This holds for any $\delta$-curve $J$. To extend the estimate from $\delta$-curves and using \eqref{eq: 2-reasonable}, we set $\overline C= 2\cdot 12 C \tilde C \tau^5\delta^{-2}$ and $\overline c = c2C\tau^4(1-4\tau/\Lambda^{n_0})^{-1}\delta^{-1}$: This way, for any curve $J$,

    \begin{align*}
        &m_J(x\in T_\varepsilon^{-n}\overline {J_n}: r_n(x)<\varepsilon')
        \\
        \leq&
        \overline C m_{J}(x\in   J: r_0(x)<\varepsilon'/\overline \Lambda^{n})
        \\&+
        \overline C m_{J}(x\in J: r_0(x)<\frac{\varepsilon\varepsilon'}{\overline \Lambda^{n}})
        \\&+
        \overline C m_{J}\left(x\in  J: r_0(x)<\frac{\varepsilon'}{\overline\Lambda^{n-N_\varepsilon}}\right)
        \\&+
        \overline c \varepsilon'.
    \end{align*}

    for any curve $J$. Notice that the first expression comes from estimates regarding $\overline {J^{(L)}_n}  = \cup_{j=0}^{\infty} \overline {(J_L^{(j)})_n}$ and $\overline {J^{(R)}_n}  = \cup_{j=0}^{\infty} \overline {(J_L^{(j)})_n}$, the second is an estimate for $\overline {J^{(L\to R)}_n}  = \cup_{j=0}^{\infty} \overline {(J_{L\to R}^{(j)})_n}$, and the third one is corresponding to $\overline {J^{(R\to L)}_n}  = \cup_{j=0}^{\infty} \overline {(J_{R\to L}^{(j)})_n}$. Hence the statement of the lemma is proven.

\end{customproof}

\subsection{Tightness lemma}
Recall that $\mu_0$ and $\mu_1$ is the ACIM of $T_0|_{I_L}$ and $T_0|_{I_R}$, respectively. Recall also from section~\ref{section: process} that $t_k^{\varepsilon}$ is the time when the process associated to $\hat{T}_{\varepsilon}$ -- or after conjugacy, to $T_{\varepsilon}$ -- switches between intervals $I_L$ and $I_R$, while $\Time^\varepsilon_{k+1}=t_{k+1}^{\varepsilon}-t_k^{\varepsilon}$.
The following is a useful proposition from \cite{difcoef10} on the tightness of the process of switching intervals, it states that sudden consecutive switches are unlikely (note this process runs on a $\frac{1}{\varepsilon}$ time scale.)

\begin{lemma}\label{prop: tightness}
    For $j=0$ or $1$, given $S>0$ and $\delta>0$ arbitrary small, there exists $\sigma>0$ such that, for all $\varepsilon$ sufficiently small,
    \[
    \mu_j(\exists k \text{ with } t^\varepsilon_k \leq S/\varepsilon \text{ and } \Time^\varepsilon_{k+1} \leq \sigma/\varepsilon) \leq \delta.
    \]
\end{lemma}

\begin{customproof}{Proof}
    The proof in our framework differs from the one presented for \cite[Proposition 2]{difcoef10} only because we have a quite non-standard growth lemma. Hence, we will briefly explain the ideas of the proof, and elaborate on the part that is specific for our framework.

    For $k=0$, we can utilize Proposition \ref{prop: open}.

    Set $\xi(j)\in\{L, R \}$. The corresponding limit measure $\mu_j$ is absolutely continuous, so it is sufficient to prove the proposition for $Leb$. We will denote probabilities with respect to Lebesgue measure by $\mathbb{P}$. Taking a small $r$, as in \cite{difcoef10}, one can define an ``$r$-inessential visit'' to the hole as one where the length of the smoothness component of $T^n_\varepsilon I_\xi \cap H_\varepsilon$ is less than $r \varepsilon$. Let $x\in E_{n, \varepsilon}$ denote the event that $x$ has an $r$-essential visit to $H_\varepsilon = H_0\cup H_1$ at time $n$. They also introduce
    \begin{align*}
        S_{n,m} = \sum_{i=n+1}^{n+m} 1_{T_\varepsilon^i x \in H_\varepsilon}.
    \end{align*}

    We need that the probability of an $r$-inessential visit to $H_\varepsilon$ up to time $S/\varepsilon$ is small enough. Recall that the probability of a point $x\in I_R$ hitting $H_0$ before time $N_\varepsilon$ is $0$. Hence, by Lemma \ref{prop: growth lemma} (a) and (c), the probability of an $r$-inessential visit to $H_0$ up to time $k\in\N^+$ $Cr\varepsilon k$. We have the same estimate for $H_1$ by Lemma \ref{prop: growth lemma} (a) and (b).  Thus, the probability of an $r$-inessential visit to $H_\varepsilon$ up to time $S/\varepsilon$ is at most $CSr$. Hence, it suffices to show that
    \begin{align*}
        \max_{n \leq S/\varepsilon} \mathbb{P}(S_{n,\sigma/\varepsilon} > 0 | E_{n,\varepsilon}(x)) \to 0, \; \text{as } \varepsilon \to 0, \sigma \to 0.
    \end{align*}
    By \ref{C2} for any fixed $M_0$, $S_{n, M_0}(x)=0$
    for all $x \in E_{n,\varepsilon}$ provided that $\varepsilon$ is small enough.
    This is when we have to use Lemma \ref{prop: growth lemma} and branch off from the original proof.

    As we do not have a uniform growth lemma, but one that depends on the starting and ending hole, we need to handle the cases separately.
    For this purpose, we will use the notation
    \begin{align*}
        L_{i,l}^k
        =
        Leb(x\in I_{\xi(j)}\cap E_{n,\varepsilon}: T_\varepsilon^{n+k}x\in H_l,\; T^nx\in H_i).
    \end{align*}
    for $i,l\in \{0,1\}$.
    We use the notation $J_{n, i} = T^n(I_{\xi(j)}\cap (T_\varepsilon^{-n} H_i) \cap E_{n,\varepsilon})$
    and $\beta_i=\lim_{\varepsilon\to 0}\mu_i(H_i)/\varepsilon$. Then
    \begin{align*}
            L_{i,l}^k
            \leq&
            Leb(I_{\xi(j)}) m_{I_{\xi(j)}}(x\in I_{\xi(j)}: T_\varepsilon^{n+k}x\in H_l,\; T^nx\in H_i, \; E_{n,\varepsilon}(x))
            \leq
            \tau Leb(I_{\xi(j)})
            m_{J_{n, i}}(T_\varepsilon^{k}x\in H_l)
            \\
            \leq&
            \frac{\tau}{2}
            m_{J_{n, i}}(x \in T_\varepsilon^{-k} \overline{(J_{n, i})_k}: r_k(x)< \beta_l \varepsilon).
    \end{align*}
    By Lemma \ref{prop: growth lemma} (a), and using that all components in $J_{n,i}$ are longer that $\varepsilon$:
    \begin{align*}
        L_{i,i}^k
        \leq&
        \frac{\overline C \tau}{2}
        m_{J_{n, i}}(r_0(x)< \beta_i \varepsilon/\overline \Lambda^k) + \frac{\tau \overline c \beta_i \varepsilon}{2}
        \\
        \leq&
        \frac{\overline C \tau\beta_i}{\overline \Lambda^k r} + c \varepsilon
        .
    \end{align*}

    Notice that $Leb(H_i\cap T_\varepsilon^{-k}H_i)=0$ for all $1\leq k \leq N_\varepsilon$. Hence,
    \begin{align*}
        \sum_{k=1}^{\sigma/\varepsilon}
        L_{i, i}^k
        \leq&
        \sum_{k=N_\varepsilon}^{\sigma/\varepsilon}
        L_{i,i}^k
        \leq
        \frac{\overline C\tau\beta_i}{\overline \Lambda^{N_\varepsilon}(1-\overline\Lambda)r} + c\sigma,
    \end{align*}
    which can be made arbitrarily small.
    By the same argument, using Lemma \ref{prop: growth lemma} (c), we can see that
    \begin{align*}
        L_{0, 1}^k
        \leq&
        \frac{\overline C \tau\beta_1  \varepsilon}{\overline \Lambda^k r} + c \varepsilon,
    \end{align*}
    thus
    \begin{align*}
         \sum_{k=1}^{\sigma/\varepsilon}
        L_{0, 1}^k
        \leq&
        \frac{\overline C \tau\beta_1 \varepsilon}{(1-\overline\Lambda)r} + c\sigma.
    \end{align*}
    Finally, by Lemma \ref{prop: growth lemma} (b),
    \begin{align*}
        L_{1, 0}^k
        \leq&
        \frac{\overline C \tau\beta_0 }{\overline \Lambda^{k-N_\varepsilon} r} + c \varepsilon,
    \end{align*}
    Notice that $Leb(H_0\cap T_\varepsilon^{-k}H_1)=0$ for all $k \leq N_\varepsilon$, but this is not sufficient.
    We will show, that for each $N_0\in \N^+$
    \begin{align*}
        L_{1, 0}^k
        =
        0
    \end{align*}
    holds for all $k<N_0 + N_\varepsilon$  for $\varepsilon>0$ small enough.
    A point $x\in H_1$ hits all $I_j$ $j\leq N_\varepsilon$, before reaching $H_0$. Furthermore, $T^{N_\varepsilon+1}x\in I_{N_\varepsilon}\subset \cup_{j=1}^{N_\varepsilon}I_j$, which is a shrinking neighborhood of $1/2$, which is a fixed point of continuity of $T_\varepsilon^{m}$ for all $m$. This means that for any $N_0\in \N^+$ we can find a neighborhood of $1/2$ such that any point of that neighborhood will not intersect $H_0$ until time $N_0$. If $\varepsilon$ is small enough, $T^{N_\varepsilon+1}H_0$ is within that neighborhood.

    \begin{align*}
        \sum_{k=1}^{\sigma/\varepsilon}
        L_{1, 0}^k
        \leq&
        \sum_{k=N_\varepsilon + N_0}^{\sigma/\varepsilon}
        L_{1, 0}^k
        \leq
        \frac{\overline C \tau\beta_0 }{\overline \Lambda^{N_0}(1-\overline\Lambda)r} + c\sigma.
    \end{align*}
\end{customproof}

\subsection{Proof of Theorem \ref{theo: markov}}
Notice that we can calculate the probabilities of events regarding the system $\T_\varepsilon$ w.r.t.~the measures $\hat{\mu}_{\varepsilon}$ by calculating the probabilities of the conjugate events regarding the system $T_\varepsilon$ w.r.t. the conjugate  measures $\mu_{\varepsilon}$.

The beginning of the proof presented for \cite[Theorem 1]{difcoef10} is suitable for our case, but it introduces the notations used in the part specific to our framework, hence we recall it here to make our exposition self contained.

The proof is by induction on $p$. Let $\xi(j)\in\{L, R \}$ denote the side corresponding to the density $\phi_j$. First, for $p = 1$, if the initial distribution on $I_{\xi(j)}$ is according to some $\rho\in BV(I_{\xi(j)})$, then by Proposition \ref{prop: open} (c),
\begin{align*}
    \mu_j(\tau^\varepsilon_1 =n, z(t^\varepsilon_1)=r)
    =&
    \int_{H_j}\perron_{j, \varepsilon}^{n}(\rho) \, dx
    =
    \lambda_{j, \varepsilon}^{n}\int_{H_j}Q_{j, \varepsilon}\rho \, dx + O(\theta^n).
\end{align*}
If $n\varepsilon \approx t$ then one can further utilize Proposition \ref{prop: open} to ensure that
\begin{align*}
    \sum_{n \in \Delta/\varepsilon} \lambda_{j, \varepsilon}^{n}\int_{H_j}Q_{j, \varepsilon}\rho \, dx
    \to
    \beta_j\int_{\Delta}e^{-\beta_j t} \, dt \int_{I_{\xi(j)}} \rho \, dx
\end{align*}
as a Riemann integral.

Now, suppose the statement is known for some $p$. Define
\[
\Omega = \{ \varepsilon \Time^\varepsilon_k \in \Delta_k, z(t^\varepsilon_k) = r_k, \text{ for } k = 1, \dots, p \}.
\]
To carry out the induction step, it suffices to show that
\begin{equation}\label{eq: finite dimensional convergence induction}
\mu_j( \varepsilon \Time^\varepsilon_{p+1} \in \Delta_{p+1}, z(T_\varepsilon^{p+1}) = r_{p+1}, \Omega) \to \mu_j(\Omega) \beta_{r_p} \int_{\Delta_{p+1}} e^{-\beta_{r_p} t} dt.
\end{equation}
Set
\[
(\perron_\Omega A)(x) = \sum_{T_\varepsilon^{t^\varepsilon_k} y = x} \frac{A(y)}{(T_\varepsilon^{t^\varepsilon_k})'(y)},
\]
where the sum is taken over $y \in \Omega$. Then
\begin{equation}
\label{eq:generalp}
    \mu_j(\Time^\varepsilon_{p+1} = n, z(T_\varepsilon^{p+1}) = r_{p+1}, \Omega) = \int_{H_{r_p}} \perron_{j,\varepsilon}^n (\perron_\Omega(1_{I_j} \phi_j)) dx.
\end{equation}
Taking a small $\sigma$, this may be rewritten as
\[
\perron_{j,\varepsilon}^n (\perron_\Omega(1_{I_j} \phi_j))
=
\perron_{j,\varepsilon}^{n - \sigma / \varepsilon} \left( \perron_{j,\varepsilon}^{\sigma / \varepsilon} (\perron_\Omega(1_{I_j} \phi_j)) \right).
\]

In order to use induction, we need to prove that for a fixed $p$,
$\BV{\perron_{j,\varepsilon}^{\sigma / \varepsilon} (\perron_\Omega(1_{I_j} \phi_j)) }$
is uniformly bounded in $\varepsilon>0$. In the setting of \cite{difcoef10}, this follows from the regular behavior of $\perron_\varepsilon$ with respect to the norm $\BV{\cdot}$. For our systems, we will show that, for any $p\ge 1$, there exist $K_p>0$ such that

\begin{align}\label{eq: induction}
    \BV{\perron_{\Omega_p}}\leq \frac{K_p}{\varepsilon^{2p}}
\end{align}
holds for all $p\geq 1$, where $\Omega_p = \Omega$ is is used to annotate the dependency of the set $\Omega$ on $p$.

Notice that for $p=1$,
\begin{align*}
    \perron_{\Omega_1} f = \sum_{n=a_1/\varepsilon}^{b_1/\varepsilon} \perron_\varepsilon ( \perron_{r_1, \varepsilon}^{n-1} f 1_{H_{r_1}})
\end{align*}

\begin{align*}
    \BV{\perron_{\Omega_1} f}
    \leq&
    \sum_{n=a_1/\varepsilon}^{b_1/\varepsilon} \BV{\perron_\varepsilon (( \perron_{r_1, \varepsilon}^{n-1} f) 1_{H_{r_1}})}
    \leq
    \sum_{n=a_1/\varepsilon}^{b_1/\varepsilon} \BV{\perron_\varepsilon}
    \BV{ (\perron_{r_1, \varepsilon}^{n-1} f) 1_{H_{r_1}}}
    \\
    \leq&
    \frac{K}{\varepsilon}\sum_{n=a_1/\varepsilon}^{b_1/\varepsilon}
    \BV{ (\perron_{r_1, \varepsilon}^{n-1} f) 1_{H_{r_1}}}
    \leq
    \frac{K'}{\varepsilon}\sum_{n=a_1/\varepsilon}^{b_1/\varepsilon}
    \BV{ \perron_{r_1, \varepsilon}^{n-1} f}
    \\
    \leq&
    \frac{K'}{\varepsilon}\sum_{n=a_1/\varepsilon}^{b_1/\varepsilon}
     \lambda_{r_1, \varepsilon}^{n-1} \nu_{r_1}(f) \BV{\phi_{r_1, \varepsilon}}
     +
     \theta^{n-1}\BV{f}
     \\
    \leq&
    \frac{K''}{\varepsilon}
     \frac{
     e^{-\beta_{r_1}a_1}-e^{-\beta_{r_1}b_1}}{1-\lambda_{r_1, \varepsilon}} \nu_{r_1}(f) \BV{\phi_{r_1, \varepsilon}}
     +
     \frac{K'}{\varepsilon}\frac{\theta^{a_1/\varepsilon}}{1-\theta}\BV{f}
     \leq
     \frac{K_1}{\varepsilon^2}\BV{f},
\end{align*}
where we have used Proposition~\ref{prop: open}. Now assume that we have \eqref{eq: induction} for $p-1$. Then,
\begin{align*}
    \BV{\perron_{\Omega_{p}} f}
    \leq&
    \sum_{n=a_p/\varepsilon}^{b_p/\varepsilon} \BV{\perron_\varepsilon ( (\perron_{r_p, \varepsilon}^{n-1} \perron_{\Omega_{p-1}} f) 1_{H_{r_p}})}
    \\
    \leq&
    \frac{K''}{\varepsilon}
     \frac{
     e^{-\beta_{r_p}a_p}-e^{-\beta_{r_p}b_p}}{1-\lambda_{r_p, \varepsilon}} \nu_{r_p}(\perron_{\Omega_{p-1}} f) \BV{\phi_{r_p, \varepsilon}}
     +
     \frac{K'}{\varepsilon}\frac{\theta^{a_p/\varepsilon}}{1-\theta}\BV{\perron_{\Omega_{p-1}} f}
     \leq
     \frac{K_p}{\varepsilon^{2p}}\BV{f}.
\end{align*}

With this in mind we can see that
\begin{align*}
    \BV{\perron_{j,\varepsilon}^{\sigma / \varepsilon} (\perron_\Omega(f)) }
    \leq&
    C(\alpha^{\sigma/\varepsilon}\V{\perron_\Omega f}
    +
    \Lnorm{\perron_\Omega f})
    +
    \Lnorm{\perron_\Omega f}
    \\
    \leq&
    C\alpha^{\sigma/\varepsilon}\BV{\perron_\Omega f}
    +
    (C+1)\Lnorm{f}
    \\
    \leq&
    C\alpha^{\sigma/\varepsilon}\varepsilon^{-2p}\BV{f}
    +
    (C+1)\Lnorm{f},
\end{align*}
which is uniformly bounded.

Hence, as in \cite{difcoef10}, utilizing the $p = 1$ case, we find that the expression \eqref{eq:generalp} is asymptotic to
\begin{align*}
    &\beta_{r_p} \int_{\Delta_{p+1}} e^{-\beta_{r_p} t} dt \int_I \perron_{j,\varepsilon}^{\sigma / \varepsilon} (\perron_\Omega(1_{I_j} \phi_j)) dx (1 + o_{\sigma \to 0}(1)).
\end{align*}

Then, Using Lemma \ref{prop: tightness},

\begin{align*}
    \mu_j(\Omega)\geq\int_{I} \perron_{j,\varepsilon}^{\sigma / \varepsilon} (\perron_\Omega(1_{I_j} \phi_j)) dx
    &
    \\
    =\mu_j((\Time_{p+1}^{\varepsilon}>\sigma/\varepsilon) \cap \Omega)
    &\geq
    \mu_j\left((\forall k \quad t^\varepsilon_k \leq pS/\varepsilon \Rightarrow \Time^\varepsilon_{k+1} > \sigma/\varepsilon)\cap \Omega \right)
    \\
    &=
    \mu_j\left(\neg(\exists k \text{ with } t^\varepsilon_k \leq pS/\varepsilon \text{ and } \Time^\varepsilon_{k+1} \leq \sigma/\varepsilon)\cap \Omega \right)
    \\
    &=
    \mu_j(\Omega) + o_{\sigma \to 0}(1)
\end{align*}

Since $\sigma$ is arbitrary, this proves \eqref{eq: finite dimensional convergence induction}, completing the induction and the proof of Theorem \ref{theo: markov}.

\section{The proof of Theorem \ref{theo: diffusion}}\label{section: proof of theo 3}

Recall that Lemma \ref{lemma: continuity} implies that for any $f\in BV$, for any $n\in \N^+$,
$\int_I \perron_\varepsilon^n f = \int_I \perron^n f + O(\varepsilon)\BV{f}$. One of the advantages of this is that $\perron$ can be decomposed into the direct sum of operators acting on the invariant intervals of $T$. We will need one more tool to prove Theorem \ref{theo: diffusion}.

\subsection{The contraction lemma}

The proof of the following lemma which is stated for the conjugated system relies on Proposition \ref{prop: Lasota-Yorke (III)} and Theorem \ref{theo: markov}, see \cite[Proposition 2]{difcoef10}.

\begin{prop}\label{prop: contraction}
    There exists a $\kappa>0$  and $\eta<1$ such that for all $f\in BV_\varepsilon$ with $\int f = 0$,
    \begin{align*}
        \BVa{\perron_\varepsilon^{ n\kappa / \varepsilon}f} \ll \eta^n\BVa{f}
    \end{align*}
\end{prop}
The main idea of the proof is that if for every $\delta<1$, there is a $\kappa>0$ such that
\begin{align*}
    \|\mathcal{L}^{\kappa/2\varepsilon}_{\varepsilon}A\|_{L^{1}} \leq \delta\BV{A}
\end{align*}
holds, then Proposition \ref{prop: Lasota-Yorke (III)} would yield the statement of the lemma.
By the mixing of $T$ on both $I_L$ and $I_R$, and Lemma \ref{lemma: continuity}, for any $c>0$, if the expected value of $A\mathbb{1}_{I_\xi}$ is small enough for each $\xi\in\{L, R\}$, then for an any $n_0$ large enough,
\begin{align*}
\|\mathcal{L}^{n_{0}}_{\varepsilon}A\|_{L^{1}} \leq c\BV{A}.
\end{align*}
In the last step, it is shown that by the convergence of the process described in Theorem \ref{theo: markov}, for large enough $\kappa>0$, the expected value of $(\perron_\varepsilon^{\frac{\kappa}{2\varepsilon}-n_0}A)\mathbb{1}_{\xi}$ is small enough for each $\xi\in\{L, R\}$ for any $A\in BV_0$.

\subsection{Proof of Theorem \ref{theo: diffusion}}

Notice that the diffusion coefficient of an observable $\hat X\in BV(I)$ for the system $\T_\varepsilon$ is the same as the diffusion coefficient of the conjugate observable $X\in BV_\varepsilon(I)$ with respect to the system $T_\varepsilon$. Let $\bar{X}=X-\mu_{\varepsilon}(X)$.

\begin{align*}
    \varepsilon\sum_{n\geq S/\varepsilon}\mu_{\varepsilon}(\bar{X}\bar{X}\circ T_{\varepsilon}^{n})
    \leq&
    \sum_{k\geq S}
    \mu_{\varepsilon}(\bar{X}\bar{X}\circ T_{\varepsilon}^{k/\varepsilon})
    \leq
    \sum_{k\geq S}
    \int \bar{X}(\bar{X}\circ T_{\varepsilon}^{k/\varepsilon}) \phi_\varepsilon
    \\
    \leq&
    \sum_{k\geq S}
    \int (\perron_{\varepsilon}^{k/\varepsilon}\bar{X}\phi_\varepsilon )\bar{X}
    \leq
    \sum_{k\geq S}
    \BV{\perron_{\varepsilon}^{k/\varepsilon}\bar{X}\phi_\varepsilon }\int\bar{X}
    \\
    \leq&
    \sum_{k\geq S}
    \BVa{\perron_{\varepsilon}^{k/\varepsilon}\bar{X}\phi_\varepsilon }\int\bar{X}
\end{align*}

Hence, by Proposition \ref{prop: contraction} for each $\delta$ we can find $S$ such that
\[\left|\varepsilon\sum_{n\geq S/\varepsilon}\mu_{\varepsilon}(\bar{X}\bar{X}\circ T_{\varepsilon}^{n})\right|<\delta\]
so it suffices to get the asymptotics of $\mu_{\varepsilon}(\bar{X}\bar{X}\circ T_{\varepsilon}^{n})$ for $n\approx t/\varepsilon$ where $t\leq S$.

The rest of the proof is the same as in \cite[Theorem 2]{difcoef10}, so we only recall the main ideas.

One needs to estimate
\begin{align*}
    D_{n}
    =
    \int \bar{X}(T_{\varepsilon}^{n_{0}}y)\left[\mathcal{L}_{\varepsilon}^{n-2n_{0}}(\mathcal{L}_{\varepsilon}^{n_{0}}(\bar{X}\phi_{\varepsilon})\right](y)dy.
\end{align*}

By Lemma \ref{lemma: continuity},
\begin{align*}
    \int \bar{X}(T_{\varepsilon}^{n_{0}}y)\left[\mathcal{L}_{\varepsilon}^{n-2n_{0}}(\mathcal{L}_{\varepsilon}^{n_{0}}(\bar{X}\phi_{\varepsilon})\right](y)dy
    =&
    \int \bar{X}(T_{\varepsilon}^{n_{0}}y)\left[\mathcal{L}_{\varepsilon}^{n-2n_{0}}(P(\bar{X}\phi_{\varepsilon})\right](y)dy + O(\theta^{n_0})  + o_{\varepsilon\to 0}(1).
\end{align*}
and the mixing of $T$ can be used on both $I_L$ and $I_R$ to see that
\begin{align*}
    \int \bar{X}(T_{\varepsilon}^{n_{0}}y)\left[\mathcal{L}_{\varepsilon}^{n-2n_{0}}(P(\bar{X}\phi_{\varepsilon}))\right](y)dy
    =&
    \sum_{k=0,1} \int_{I_{\xi(k)}} \bar{X}\phi_k
    \int_{I_{\xi(k)}}\mathcal{L}_{\varepsilon}^{n-2n_{0}}(P (\bar{X}\phi_{\varepsilon}))
    \\&+
    (O(\theta^{n_0})  + o_{\varepsilon\to 0}(1))\BV{X}\BVa{P (\bar{X}\phi_{\varepsilon})}.
\end{align*}
In the second integral, from the deconstruction of the projection $P$, an expression
$ f_{j, k, \varepsilon}(n) = \int_{I_{\xi(k)}}\perron_\varepsilon^{n-2n_{0}}\phi_j =  \mu_{j}(T^{n-2n_{0}}_{\varepsilon}x\in I_{k})$ appears.
By Theorem \ref{theo: markov}, $f_{j, k, \varepsilon}(t/\varepsilon)\to p_{j,k}$ locally uniformly.
Using this,
\begin{align*}
    \lim_{\varepsilon\to 0}\varepsilon\sum_{n\leq S/\varepsilon}\mu_{\varepsilon}(\bar{X}\bar{X}\circ T^{n}_{\varepsilon})
    =
\sum_{jk}\int_{0}^{S}\left[p_{j}\mathbf{X}(j)\mathbf{X}(k)p_{jk}(t)-\left(\sum_{j}p_{j}\mathbf{X}(j)\right)^{2}\right]dt.
\end{align*}
as a Riemann integral. Since $S$ is arbitrary, one can obtain Theorem \ref{theo: diffusion} by taking  $S\to\infty$.

\section{Examples}\label{section: examples}

\begin{figure}[H]
\centering
\begin{tabular}{c c}
\begin{tikzpicture}[scale=6]
\draw (0,0) rectangle (1,1);
\draw[gray] (0,0) -- (1,1);
\draw[thick] (0,1) -- (0.15,0) -- (0.3, 0.3) -- (0.45, 0) -- (1, 1);
\node at (0.3, -0.04) {\small{$a$}};
\node at (0.15, -0.04) {\small{$a/2$}};
\node at (0.45, -0.04) {\small{$3a/2$}};

\draw[] (0.3, 0.01) -- (0.3, -0.01);

\node at (0.15, 0.9) {$\T$};

\end{tikzpicture}

&
\begin{tikzpicture}[scale=6]
\draw (0,0) rectangle (1,1);
\draw[gray] (0,0) -- (1,1);
\draw[thick] (0,1) -- (0.15,0) -- (0.3, 0.35) -- (0.45, 0) -- (1, 1);
\node at (0.3, -0.04) {\small{$a$}};
\node at (0.15, -0.04) {\small{$a/2$}};
\node at (0.45, -0.04) {\small{$3a/2$}};

\draw[] (0.3, 0.01) -- (0.3, -0.01);

\draw[dotted] (0, 0.35) -- (0.3, 0.35);
\node at (-0.08, 0.35) {\small{$\frac{2+\varepsilon}{2}a$}};

\node at (0.15, 0.9) {$\T_\varepsilon$};

\end{tikzpicture}
\end{tabular}
\caption{The limit map $\T$ and the perturbed map $\T_\varepsilon$.}\label{fig: example}
\end{figure}
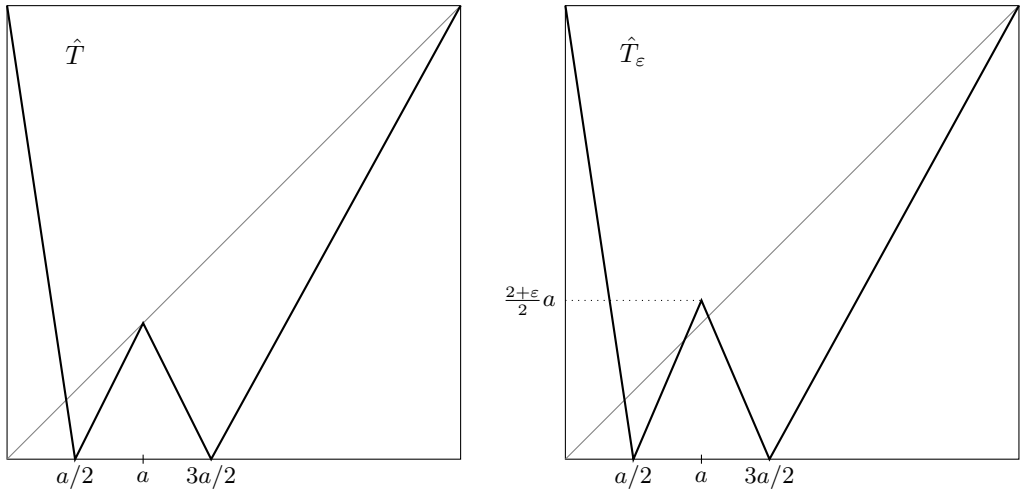

In honour of Keller's W-map and the other families defined in \cite{singular limits, Ulams method, harmonic averages, various instabilities, spectrum instability, family}, we would like to define a family of W-shapes maps, where the fixed point $a$ is not necessarily the point $1/2$.
We fix a parameter $a\in (0,\frac23)$, and define the following family for $\varepsilon \geq 0$, see Fig \ref{fig: example}.

\begin{align*}
    \T_\varepsilon(x)
    =
    \begin{cases}
        -\frac{2}{a}x + 1 & \text{ if } x\leq a/2,
        \\
        (2+\varepsilon)x - \frac{2+\varepsilon}{2}a & \text{ if } a/2 < x\leq a,
        \\
        -(2+\varepsilon)x + \frac{6 + 3\varepsilon}{2}a & \text{ if } a < x\leq 3a/2,
        \\
        \frac{x}{1-3a/2} -\frac{3a/2}{1-3a/2} & \text{ if } 3a/2 < x\leq 1.
    \end{cases}
\end{align*}


Then $h_\varepsilon = 1+ \varepsilon/2$, which means that the family satisfies \ref{C7}. Furthermore, $Leb(\mathcal{B})=a\cdot \varepsilon+o(\varepsilon)$, hence one can tailor the $l.h.r.$ to an extent. Since $\T_0$ is Markovian for all $a \in (0, \frac{2}{3})$, it is easy to see that its invariant density is
\begin{align*}
    \hat \phi_0(x)
    =
    \begin{cases}
    \frac{3}{1+2a}, & 0 \leq x \leq a,
    \\
    \frac{1}{1+2a}, & a < x \leq 1.
    \end{cases}
\end{align*}
Then, by Theorem \ref{theo: convex}, $l.h.r. = (a^2(2-a) + 2a(1-3a/2)(1-a))^{-1}$.

The full generality of our framework and its usefulness depends on the restrictiveness of our conditions. The most limiting detail of \ref{C0} and \ref{C5} is that the family can only have finite number of branches, and that the perturbations can not introduce new critical points. Conditions \ref{C1} and \ref{C6} introduce the complication of the turning fixed point. Recall, that in Remark \ref{remark: C6 opposite} the case opposite to \ref{C6} is discussed. Finally, the technical conditions \ref{C2}, \ref{C3} and \ref{C4} are generic. If the limit map $\T$ is Markovian, then
\ref{C2} and \ref{C4} are simple to check.

\section{Conclusions and questions}\label{section: conclusions}

The original aim of this paper was to introduce a unified framework that accomodates all of the systems discussed in \cite{singular limits, Ulams method, harmonic averages, various instabilities, spectrum instability, family}, which it achieved. Notice, that the perturbed systems in \cite{singular limits} have additional critical points to that of the limit map. However, these are exactly at places where the conjugation map would force critical points in our framework, hence the methods of the proofs can be applied.

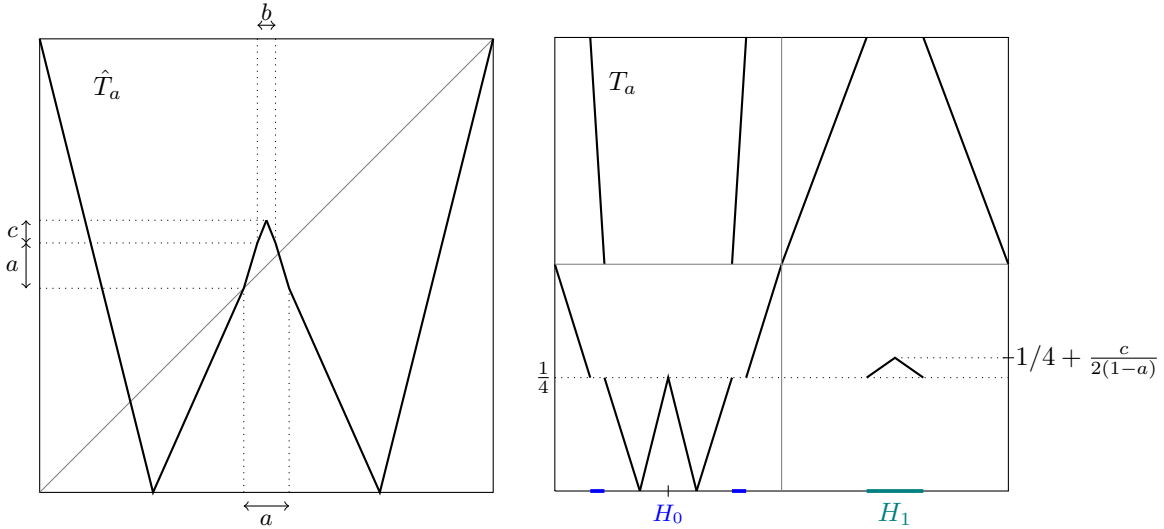
\begin{figure}[H]
\centering
\begin{tabular}{c c}
\begin{tikzpicture}[scale=6]
\draw (0,0) rectangle (1,1);
\draw[gray] (0,0) -- (1,1);
\draw[thick] (0,1) -- (0.25,0) -- (0.45, 0.45) -- (0.48, 0.55) -- (0.5, 0.6);
\draw[thick] (0.5, 0.6) -- (0.52, 0.55) -- (0.55, 0.45) -- (0.75, 0) -- (1, 1);
\draw[dotted] (0.45, 0) -- (0.45, 0.45);
\draw[dotted] (0.55, 0) -- (0.55, 0.45);
\draw[<->] (0.45, -0.03) -- (0.55, -0.03);
\node at (0.5, -0.06) {\small{$a$}};

\draw[dotted] (0, 0.45) -- (0.45, 0.45);
\draw[dotted] (0, 0.55) -- (0.48, 0.55);
\draw[<->] (-0.03, 0.45) -- (-0.03, 0.55);
\node at (-0.06, 0.5) {\small{$a$}};

\draw[dotted] (0, 0.6) -- (0.5, 0.6);
\draw[<->] (-0.03, 0.55) -- (-0.03, 0.6);
\node at (-0.06, 0.575) {\small{$c$}};

\draw[dotted] (0.48, 1) -- (0.48, 0.55);
\draw[dotted] (0.52, 1) -- (0.52, 0.55);
\draw[<->] (0.48, 1.03) -- (0.52, 1.03);
\node at (0.5, 1.06) {\small{$b$}};

\node at (0.15, 0.9) {$\T_a$};

\end{tikzpicture}

&

\begin{tikzpicture}[scale=1.5]
    \draw[black] (0, 0) -- (0, 4) -- (4, 4) -- (4, 0) -- (0, 0);

    \draw[thick] (0, 2) -- (2.5/8, 1);

    \draw[thick] (2.5/8, 4) -- (3.5/8, 2);

    \draw[thick] (3.5/8, 1) -- (6/8, 0);

    \draw[thick] (3/4, 0) -- (1, 1) -- (5/4, 0);

    \draw[thick] (5/4, 0) -- (12.5/8, 1);

    \draw[thick] (12.5/8, 2) -- (13.5/8, 4);

    \draw[thick] (13.5/8, 1) -- (2, 2) -- (3-0.25, 4);

    \draw[thick] (3+0.25, 4) -- (4, 2);

    \draw[thick] (3-0.25, 1) -- (3, 1.175) -- (3+ 0.25, 1);

    \draw[gray]
    (0, 2) -- (4, 2);
    \draw[dotted]
    (3, 1.175) -- (4, 1.175);
    \draw[dotted]
    (0, 1) -- (4, 1);
    \draw[gray]
    (2, 0) -- (2, 4);

    \node at (-0.1, 1) {$\frac{1}{4}$};
    \draw[] (1, -0.05) -- (1, 0.05);

    \draw[] (4 - 0.05, 1.175) -- (4 + 0.05, 1.175);
    \node at (4+0.7, 1.14) {$1/4+\frac{c}{2(1-a)}$};

    \draw[blue, line width = 0.5mm] (2.5/8, 0) -- (3.5/8, 0);
    \draw[blue, line width = 0.5mm] (12.5/8, 0) -- (13.5/8, 0);
    \node at (1, -0.2) {\small\textcolor{blue}{$H_0$}};

    \draw[teal, line width = 0.5mm] (3-0.25, 0) -- (3+0.25, 0);
    \node at (3, -0.2) {\textcolor{teal}{$H_1$}};

    \node at (0.6, 3.6) {$T_a$};

\end{tikzpicture}
\end{tabular}
\caption{The map $\T_a$ and its conjugate map $T_a$.}\label{fig: example abc}
\end{figure}

As another example, on Figure \ref{fig: example abc} we can see the map defined in \cite{singular limits}, and its conjugate map.
Instead of the Markovian family assumed in \cite{singular limits}, our setting applies if we just take the conditions $b\ll c$ and $b \ll a^{1+\delta}$ for some positive $\delta$.

Here $\mu_0(H_0)=\frac{a}{2-2a}$ and $\mu_1(H_1)=b/a$. The authors of \cite{singular limits} denote $\beta=\lim_{a\to 0}a^2/b$, hence $l.h.r.=2/\beta$, a result which now extends to a much larger class of (not necessarily Markovian) maps.

Furthermore, our results can be interpreted also from the pragmatic viewpoint of ACIM stability. Given a system $(\T, \hat\mu)$ that we can only measure with errors $(\T_\varepsilon, \hat \mu_\varepsilon)$, can we say anything meaningful about the probability of events with respect to the original system? If the family is ACIM stable, then the answer is yes. However, even in some cases of instability (Theorem \ref{theo: convex} (III)) one can recover adequate information about probabilities w.r.t. the original system.

One of the key novelties of our method is that it shows that the behavior of the localized state of the perturbed maps indeed depends only on the local nature of the maps. Our results apply to a class of systems that are required to satisfy only quite general conditions globally, and special conditions locally. Verifying \ref{C2} and \ref{C4} can be complicated, but if $\T_0$ is Markovian, it is quite simple.

\begin{remark}
    Notice that the purpose of conditions \ref{C2} and \ref{C4} was to guarantee that no critical point, other than the turning point $a$ itself, induces a complication when constructing the uniform Lasota-Yorke estimate in section \ref{section: open systems}. We believe that Condition \ref{C2} could be easily omitted. The complete relaxation of \ref{C4}, however, would result in facing the general case that would make conditions \ref{C1} and \ref{C6} redundant.
\end{remark}

We would like to pose two questions regarding ACIM instability. Both are meant to refine Keller's original conjecture regarding ACIM instability in \cite{Keller82}. The first question approaches the phenomenon by describing the \textit{shape} of the maps in the family, while the second question links ACIM instability to metastable \textit{behavior}, specifically, to the emergence of a jump Markov process on a slow timescale.

 To formulate the first question, let us introduce the following quantity.

\begin{definition}
    Let $\T: I\to I$ be a map satisfying condition \ref{C0}. The local expansion of $\T$ at $y\in I$ is
    \begin{align*}
        Exp(y)=\lim_{\delta\to 0}\left\{\frac{Leb\left(\T V_\delta(x)\right)}{Leb(V_\delta(x))}: V_\delta(x) \text{ is the connected component of } \T^{-1} ]\T x-\delta, \T x+\delta[ \text{ containing } x \right\}.
    \end{align*}
    If $a\in I$ is a point of period $n\in\N^+$, the local expansion of $\T$ along the orbit of $a$ is
    \begin{align*}
        Exp(\mathcal{O}(a))=\Pi_{k=1}^n Exp(\T^k a)
    \end{align*}.
\end{definition}

Notice that if $\D\T$ is continuous at $x\in I$, then $Exp(x)=\D\T(x)$. If $x\in I$ is a turning point of $\T$, then $Exp(x)=h(x)$ as defined in \eqref{eq: h}.

The precise and general analysis of the case of periodic orbits will be a topic of future research.

\begin{question}
    Let family $\T_\varepsilon$, $\varepsilon>0$ satisfy \ref{C0}, \ref{C5} and let $\T_0$ have exactly one finite orbit $\mathcal{O}$ that contains a critical point. Is it true that the family is ACIM unstable only if there is an eventually periodic point $b_\varepsilon$ such that  $\lim_{\varepsilon\to 0} dist(b_\varepsilon, \mathcal{O})=0$ and the local expansion $Exp(\mathcal{O})$ of $\T_0$ along $\mathcal{O}$ is at most $1$?
\end{question}

\begin{question}
     Let $\T_\varepsilon$, $\varepsilon\geq 0$ be a family of piecewise $C^2$ expanding maps satisfying conditions \ref{C0} and \ref{C5}. Furthermore, assume $\hat\phi$ is positive on at least one of the periodic critical points.  Is it true that $\mu_\varepsilon \xrightarrow{w} t\mu_0 + (1-t)\nu$ holds for some singular measure $\nu$ and $t\in ]0, 1[$ if and only for the family $\T_\varepsilon$ on an appropriate slow timescale an asymptotic jump Markov process emerges, analogous to Theorem \ref{theo: multiple markov}, where the local states are finite orbits containing critical points?
\end{question}

Here $\nu$ is a convex combination of singular measures concentrated on finite orbits of $\T$.

\appendix

\section{Multiple turning points}\label{section: multiple}

In this section we discuss the case where there are several turning points of $\T$. We will denote the conditions in this section by (C$k$'), $0 \leq k \leq 6$. Most of them will be the same as in Section \ref{section: conditions}, but we need to generalize some of them:

\begin{enumerate}
    \item [(C1')] We have turning points $A = \{a^{(1)},\dots, a^{(l)} \}\subset Int(I)$, meaning
    \begin{align*}
        \T(a^{(k)})=a^{(k)}, \quad \partial_\pm\T a^{(k)}=\lim_{x \to a^{(k)}\pm}\D\T (x), \quad sign(\partial_-\T a^{(k)}) \neq sign(\partial_+\T a^{(k)}),
    \end{align*}
    where $sign(y)$ denotes the sign of a value $y\neq 0$, and $\T$ is continuous at every point of $A$.
\end{enumerate}
We define the infinitesimal holes as $H = T^{-1}A\setminus A$.
\begin{enumerate}
    \item [(C2')] The critical points avoid all infinitesimal holes, hence $A$ as well: $A\cap T^k C_0 = \emptyset$
    \item [(C4')] There are no periodic critical points except for all $a\in A$ and possibly $\partial I$.
    \item [(C6')] We assume that $\T_\varepsilon$ are continuous at all $a^{(k)}_\varepsilon \in A_\varepsilon$, where $A_\varepsilon$ is the perturbation of $A$. Furthermore, we assume $sign(\partial_-\T a^{(k)})\T_\varepsilon (a^{(k)}_\varepsilon) > sign(\partial_-\T a^{(k)})a^{(k)}_\varepsilon$ for all $a^{(k)}_\varepsilon \in A_\varepsilon$.
\end{enumerate}

Then by the same logic as in the case of a single turning point, there are disjoint boxes $\B^{(k)}$ around all $a^{(k)}_\varepsilon\in A_\varepsilon$. We denote $\hat I = I\setminus\cup_{k\leq l} \B^{(k)}$.

We can generalize the conjugation map $U_\varepsilon$ such that $U_\varepsilon$ is piecewise linear, and it sends $\B^{(k)}$ to $]k/l, (k+1)/l]$, and $\hat I$ to $[0, 1/l]$.
This way we can define $T_\varepsilon = U_\varepsilon \circ \T_\varepsilon \circ U_\varepsilon^{-1}$.

\subsection{Limit measure}

We define
\begin{align*}
    h^{(k)}_\varepsilon = (|\partial_-\T (a^{(k)}_\varepsilon)|^{-1} + |\partial_+\T (a^{(k)}_\varepsilon)|^{-1})^{-1}
\end{align*}

\begin{align}\label{eq: cpm general}
    c_\pm^{(k)} = \lim_{\varepsilon \to 0} \frac{\partial_{\pm}\T_{\varepsilon}(a_\varepsilon^{(k)}) - \partial_{\pm}\T(a^{(k)})}{\varepsilon}
    =
    \partial_\varepsilon(\partial_{\pm}\T_{\varepsilon}(a_\varepsilon^{(k)}))|_{\varepsilon=0}
\end{align}

\begin{align*}
    c^{(k)}
    =
    \lim_{\varepsilon\to 0}\frac{|\T_\varepsilon(a^{(k)}_\varepsilon)-a^{(k)\varepsilon}|}{\varepsilon}.
\end{align*}

The holes are
 $\hat H_k = \T^{-1}(\hat I)\cap \B^{(k)}$ and
 $\hat H_{0, k} = \T^{-1}(\B^{(k)})\cap \hat I$.
 The limit hole ratios are

\begin{align*}
    l.h.r._k = \lim_{\varepsilon\to 0}\frac{m_{\B^{(k)}}(\hat H_k)}{\hat \mu (\hat H_{0, k} )}, \quad \forall k\leq l.
\end{align*}

\begin{theorem}\label{theo: multiple convex}
    Assume that the family $\T_\varepsilon$, $\varepsilon\geq 0$ satisfies conditions (C0')-(C6').
    \begin{itemize}
        \item [(I')] If there exists exactly one $k\leq l$ such that $h^{(k)}_\varepsilon\leq 1$ for all $\varepsilon >0$, then $\hat \mu_\varepsilon \xrightarrow{w}\delta_{a^{(k)}}$.
        \item [(II')] If  $h^{(k)}_\varepsilon > 1$ holds for all $k\leq l$ and $\varepsilon \geq 0$,
        then $\hat \phi_\varepsilon \xrightarrow{L^1} \hat \phi$.
        \item [(III')] If  $h^{(k)}_\varepsilon > 1$ holds for all $k\leq l$ and $\varepsilon > 0$,
        then $\hat \mu_\varepsilon \xrightarrow{w} \alpha \hat\mu + \sum_{k=1}^l\alpha_k\delta_{a^{(k)}}$, where

        \begin{align*}
            \frac{\alpha}{\alpha_k}
            =
            l.h.r._k
        \end{align*}
    \end{itemize}
\end{theorem}

\begin{remark}
    In Theorem \ref{theo: multiple convex} (III'), $l.h.r._k=0$ for any $k$ with $h^{(k)}>1$. If $h^{(k)}=1$, then
    \begin{align*}
        l.h.r._k
        &=
        \frac{|c^{(k)}_-|\partial_-\T(a^{(k)})^{-2} + |c^{(k)}_+|\partial_+\T(a^{(k)})^{-2}}
        {
        \sum_{x\in H_{0, k}}\hat\phi(x)|\D \T(x)|^{-1}
        \lim_{\varepsilon\to 0}\varepsilon^{-1}Leb(\B^{(k)})
        },\text{ where}
        \\
        \lim_{\varepsilon\to 0}\varepsilon^{-1}Leb(\B^{(k)})
        &=
        -c^{(k)}\frac{|\partial_+\T(a^{(k)})|
        +
        |\partial_-\T(a^{(k)})|}
        {
        (\partial_+\T(a^{(k)}) - \mathbb{1}_{\partial_+\T(a^{(k)})  > 0})(\partial_-\T(a^{(k)}) - \mathbb{1}_{\partial_-\T(a^{(k)})  > 0})
        }
    \end{align*}
\end{remark}

\begin{remark}\label{remark: muliple convex cases}
    Notice, that these three are the possible cases our framework allows.
    Indeed, if for some $\varepsilon_0>0$ we have $i< j \leq l$ such that $h^{(i)}_\varepsilon, h^{(j)}_\varepsilon \leq 1$, then both $\B^{(i)}$ and $\B^{(j)}$ are invariant intervals, which allows for more than one ACIM for $\T_{\varepsilon_0}$, contradicting (C5').
\end{remark}

 By Remark \ref{remark: muliple convex cases}, case (I') of Theorem \ref{theo: multiple convex} implies $h^{(i)}_\varepsilon > 1$ for all $i\neq k$ and $\varepsilon>0$. Then we can see that by the same argument as in case (I) of Theorem \ref{theo: convex} we have $\hat\mu_\varepsilon \xrightarrow{w} \delta_{a^{(k)}}$. Furthermore, if we fix $\varepsilon>0$, then  $\perron_\varepsilon^n f \to \f_\varepsilon$ for all $f \in BV$.

For case (II') we can use the same argument as in case (II) of Theorem \ref{theo: convex} to see that for every $\varepsilon>0$ there exists a unique density $\f_\varepsilon \in BV_\varepsilon$ such that $\perron_\varepsilon \f_\varepsilon = \f_\varepsilon$. Furthermore, by using the limit hole ratio argument, as at the end of the proof of case (II) of Theorem $\ref{theo: convex}$, we can see that $\f_\varepsilon \circ U_\varepsilon = \hat\phi_\varepsilon \xrightarrow{L^1}\hat\phi$.
Furthermore, if we fix $\varepsilon>0$, then $\hat\perron_\varepsilon^n f \to \hat \phi_\varepsilon$ for all $f \in BV$.

In case (III') we separate the set of indices of the turning fixed points $A$. We assume that $h^{(i)} > 1$ for all $i\leq l_1<l$ and $h^{(i)} = 1$ for all $i>l_1$. Then we know that the family $\T_\varepsilon$ behaves well on $\hat I \bigcup \cup_{i\leq l_1}\B^{(i)}$ as an open map. Then, after conjugation we define $I_L=[0, l_1/l]$ and $I_{R_j} = ]j/l, (j+1)/l]$ for all $j\geq l_1$. Notice that we can simultaneously decompose any density $f\in BV$  into components on the intervals $I_{R_j}$, which results in the norm defined below. As in \ref{section: lasota yorke}, we define the $\varepsilon$-variance

\begin{align*}
    \Va{f}= \V{f|_{[0, l_1/l]}} + \varepsilon^{-1}\sum_{j=l_1}^{l-1} \V{f|_{[j/l, (j+1)/l]}}
    +
    2\sum_{k=0}^{l-1}\BV{f|_{[k/l, (k+1)/l]}}
    \geq \V{f},
\end{align*}
and the $BV_\varepsilon$ norm
\begin{align*}
    \BVa{f} = \Va{f} + \Lnorm{f} \geq \BV{f}.
\end{align*}

Then for a $q<1$, for all $l_1 \leq j\leq l$, we can define
\begin{align*}
    H_j
    =&
     I_{R_j} \cap T^{-1}I_L
    \\
    I^{(j)}_0
    =&
    T H_j
    \\
    I^{(j)}_k
    =&
    T^{k} I^{(j)}_0 \setminus T^{k-1}I^{(j)}_0 \quad 1\leq i \leq N_\varepsilon^{(j)},
\end{align*}
where $N_\varepsilon^{(j)}$ is defined as in \eqref{eq: Na}.
Notice that $I^{(i)}_{k_1} \cap I^{(j)}_{k_2} = \emptyset$ unless $(i, k_1) = (j, k_2)$. Then for any $f\in BV$,

\begin{align*}
    f_{L}
    &=
    f|_{I_L},
    \\
    f_{R_j}
    &=
    f|_{I_{R_j}}, \quad l_1\leq j \leq l
    \\
    0 = f^{(j)}_k & : I^{(j)}_k \to \R, \quad l_1\leq j \leq l, \quad 0\leq k \leq q N^{(j)}_\varepsilon.
\end{align*}
We can define the evolution of a density $f\in BV$ as

\begin{align*}
    (\perron_\varepsilon f)_{L}
    =&
    (\perron_\varepsilon f)|_{I_L}
    =
    \perron_{I_L, \varepsilon} f_L
    +
    \sum_{j=l_1}^l \perron_{\varepsilon} f_{qN^{(j)}}^{(j)}
    ,
    \\
    (\perron_\varepsilon f)_{R_j}
    =&
    (\perron_\varepsilon f)|_{I_{R_j}}
    =
    \perron_{R_j, \varepsilon}f_{R_j}
    +
    \perron_{I_L \to I_{R_j}, \varepsilon}f_L
    ,
    \\
    (\perron_\varepsilon f)^{(j)}_0
    =&
    \perron_{H_j \to I^{(j)}_1, \varepsilon}f_{R_j}
    \\
    (\perron_\varepsilon f)^{(j)}_1
    =&
    \perron_{I^{(j)}_0 \to I^{(j)}_1, \varepsilon}f_{R_j}
    +
    \perron_{I^{(j)}_1 \to I^{(j)}_1, \varepsilon}f^{(j)}_{0},
    \\
    (\perron_\varepsilon f)^{(j)}_k
    =&
    \perron_{I^{(j)}_{k-1}\to I^{(j)}_k, \varepsilon} f^{(j)}_{k-1}
\end{align*}

Then by the same arguments as in case (III) of Theorem \ref{theo: convex}, we can prove case (III') of Theorem \ref{theo: multiple convex}.

\subsection{Markov jump process}\label{section: multiple process}

Throughout the section we assume that the following condition holds:
\begin{enumerate}
    \item [(C7')] We assume that case (III') if Theorem \ref{theo: multiple convex} holds for the family $\T_\varepsilon $, $ \varepsilon \geq 0$.
\end{enumerate}

We define a Markov process on the indices $\mathcal{S}=\{0\}\cup \{k\leq l \; | \; h^{(k)} = 1, \; \hat\phi(a^{(k)})>0 \}$. Let the rates be
\begin{align*}
    \beta_{0, k}
    =&
    \lim_{\varepsilon \to 0}\frac{\hat\mu(\hat H_{0, k})}{\varepsilon}
    =
    \lim_{\varepsilon\to 0}\varepsilon^{-1}Leb(\B^{(k)})
    \sum_{x\in H_{0, k}}\hat\phi(x)|\D \T(x)|^{-1},
    \\
    \beta_{0}
    =&
    \lim_{\varepsilon \to 0}\frac{\hat\mu(\hat H_0)}{\varepsilon}
    =
    \sum_{k\in\mathcal{S}\setminus\{0\}} \beta_{0, k}
    \\
    \beta_{k, 0}
    =&
    \lim_{\varepsilon \to 0}\frac{m_{\B^{(k)}}(\hat H_{k})}{\varepsilon}
    =
    |c^{(k)}_-|\partial_-\T(a^{(k)})^{-2}
    +
    |c^{(k)}_+|\partial_+\T(a^{(k)})^{-2}.
\end{align*}

Then the Markov process is a random walk on a star graph with the state $0$ at its center. For any initial state $r\in \mathcal{S}$ (that is, $\mathbb{P}^r(z_0^M=j)=\delta_{rj}$), the transition probabilities for any $j\in \mathcal{S}\setminus \{0\}$ are
\begin{align*}
    d\mathbb{P}^r(\mathcal{T}_i^M = t \;| \; z_{i-1}^M= j)
    =&
    \beta_{j,0} e^{-\beta_{j,0} t} dt,
    &&\mathbb{P}^r(z_i^M = 0 \;| \; z_{i-1}^M= j) = 1
    \\
    d\mathbb{P}^r(\mathcal{T}_i^M = t \;| \; z_{i-1}^M= 0)
    =&
    \beta_{0} e^{-\beta_{0} t} dt,
    &&\mathbb{P}^r(z_i^M = j \;| \; z_{i-1}^M= 0) = \frac{\beta_{0,j}}{\beta_0},
\end{align*}
for any $i\ge 1$.

\begin{theorem}\label{theo: multiple markov}
    Let $\T_\varepsilon$ satisfy conditions (C0')-(C7'). Fix $p, s$ and $S$. For any intervals $\Delta_k=[a_k, b_k]$, and numbers $r_k \in \mathcal{S}$, $k = 1, \dots, p$
    \begin{align*}
        \hat\mu(\varepsilon\Time^{\varepsilon}_{k}\in\Delta_k,\; z(t_k^{\varepsilon})
        =&
        r_k \; \forall k\leq p)
        \to \mathbb{P}^0(\Time_{k}\in\Delta_k,\; z_k^M=r_k \; \forall k\leq p),
        \\
        m_{\B^{(j)}}(\varepsilon\Time^{\varepsilon}_{k}\in\Delta_k,\; z(t_k^{\varepsilon})
        =&
        r_k \; \forall k\leq p)
        \to \mathbb{P}^j(\Time_{k}\in\Delta_k,\; z_k^M=r_k \; \forall k\leq p), \; \forall 0\neq j \in \mathcal{S},
    \end{align*}

and the convergence is uniform for $\max_k b_k \leq S$ and $\min_k a_k\geq s$.
\end{theorem}

Since for every element $a_j\in A$ we have intervals $I_k^{(j)}$, $0\leq N_\varepsilon^{(j)}$, we define
\begin{align*}
    B_\varepsilon = \bigcup_{j\in \mathcal{S}\setminus\{0\}} \bigcup_{k=0}^{N_\varepsilon^{(j)}} I_k^{(j)}
\end{align*}
and
\begin{align*}
    \overline J = J\setminus B_\varepsilon.
\end{align*}

Define $N_\varepsilon$ as $\min_j N_\varepsilon^{(j)}$.
By the same arguments, as in the proof of Lemma \ref{prop: growth lemma}, we have

\begin{lemma}[Growth lemma]\label{prop: growth lemma multiple} There exists a $\overline\Lambda>1$, and a $\overline{C}, \overline{c} >0$ such that for all $\varepsilon, \varepsilon'>0$ small enough, any curve $J \subset I$,  and all $n\in\N^+$ we have the following estimates.

\begin{itemize}
\item[(a)]
For $ J^{(L)} =  J \cap I_L \cap T_\varepsilon^{-n}I_L$ and $ J^{(R)} = \cup_{i,j=1}^l J \cap I_{R_i} \cap T_\varepsilon^{-n}I_{R_j}$,
\begin{align*}
    m_{J}\left(x\in T_\varepsilon^{-n}\overline {J^{(L)}_n} \cup T_\varepsilon^{-n}\overline {J^{(R)}_n}: r_0(x)<\varepsilon'\right)
    \leq&
    \overline C m_{J}\left(x\in J: r_0(x)<\varepsilon'/\overline\Lambda^n\right)
    +
    \overline c \varepsilon'.
\end{align*}

\item[(b)]
For $ J^{(L\to R)} = \cup_{i=1}^l  J \cap I_L \cap T_\varepsilon^{-n}I_{R_i}$,
\begin{align*}
    m_{J}\left(x\in T_\varepsilon^{-n}\overline {J^{(L \to R)}_n}: r_n(x)<\varepsilon'\right)
    \leq&
    \overline C m_{J}\left(x\in J: r_0(x)<\varepsilon\varepsilon'/\overline\Lambda^n\right)
    +
    \overline c \varepsilon'
\end{align*}

\item[(c)]
For $i\in \mathcal S\setminus\{0\}$ set $ J^{(R_i \to L)} =  J \cap I_{R_i} \cap T_\varepsilon^{-n}I_L$. Then we have
\begin{align*}
    m_{J}\left(x\in T_\varepsilon^{-n}\overline {J^{(R_i \to L)}_n}: r_n(x)<\varepsilon'\right)
    \leq&
    \overline C m_{J}\left(x\in J: r_0(x)<\varepsilon'/\overline\Lambda^{n-N_\varepsilon^{(i)}}\right)
    +
    \overline c \varepsilon'.
\end{align*}

\end{itemize}
\end{lemma}

For a centered observable $X\in BV$, which is continuous at every point of $A$, we can define an observable $\X: S\to \R$ on the state space of the Markov-process described above by setting $\X(0)=\hat\mu(X)$ and $\X(k) = \delta_{a^{(k)}}(X)=X(a^{(k)})$. Let $\textbf{D}(\X)$ denote the diffusion coefficient of $\X$. Recall that $D^\varepsilon(X)$ denotes the diffusion coefficient of $X$ w.r.t. the system $\T_\varepsilon$.

\begin{theorem}\label{theo: multiple diffusion}
    If (C0')-(C7') holds for the family $\T_\varepsilon$, $\varepsilon \geq 0$, then
    \begin{align*}
        \varepsilon D^\varepsilon(X) \to \textbf{D}(\X)
    \end{align*}
    for all $X \in BV$ which is continuous at every point of $A$.
\end{theorem}

\section*{Acknowledgements} The authors thank Joshua Peters for stimulating discussions. The support of Hungarian National Research, Development and Innovation Office (NKFIH); grants  142169 and 144059, is thankfully acknowledged. The research of KA has been implemented with the support provided by the Ministry of Culture and Innovation of Hungary from the National Research, Development and Innovation Fund (Project no. 294), financed under the DKÖP-26-1-BME-8 funding scheme.

\end{document}